\documentclass{stavanger}

\usepackage[utf8]{inputenc}
\usepackage{amssymb,amsmath,amsfonts,amsthm}
\usepackage{ dsfont }
\usepackage{hyperref}
\usepackage{natbib}
\usepackage{graphicx}
\usepackage{tikz}
\usepackage{cleveref}
\usepackage{mathtools}
\usepackage{placeins}
\usepackage{bm}
\usepackage{animate}
\usepackage[export]{adjustbox} % Aligning figures left and right
\usepackage{fancybox}
\usepackage{mathbbol}
\usepackage{setspace}
\usepackage{enumitem}

\usepackage{tikz}
\usetikzlibrary{matrix,decorations.pathreplacing,calc,fit,backgrounds}

\startlocaldefs

\newcommand{\td}{\dot t}
\newcommand{\tp}{t^\prime}
\newcommand{\xd}{\dot x}
\newcommand{\xp}{x^\prime}
\newcommand{\pd}{\dot \phi}
\newcommand{\pp}{\phi^\prime}

\pgfmathsetmacro{\myscale}{1}
\pgfkeys{tikz/mymatrixenv/.style={decoration={brace},every left delimiter/.style={xshift=8pt},every right delimiter/.style={xshift=-8pt}}}
\pgfkeys{tikz/mymatrix/.style={matrix of math nodes,nodes in empty cells,
left delimiter={(},right delimiter={)},inner sep=1pt,outer sep=1.5pt,
column sep=8pt,row sep=8pt,nodes={minimum width=20pt,minimum height=10pt,
anchor=center,inner sep=0pt,outer sep=0pt,scale=\myscale,transform shape}}}
\pgfkeys{tikz/mymatrixbrace/.style={decorate,thick}}

\tikzset{greenish/.style={
    fill=green!50!lime!60,draw opacity=0.4,
    draw=green!50!lime!60,fill opacity=0.1,
  },
  orangeish/.style={
    fill=orange!90, draw opacity=0.8,
    draw=orange!90, fill opacity=0.3,
  },
  purpleish/.style={
    fill=violet!90!pink!20, draw opacity=0.5,
    draw=violet, fill opacity=0.3,    
  }}

\usepackage[skins,theorems]{tcolorbox}
\tcbset{highlight math style={enhanced,
  colframe=red,colback=white,arc=0pt,boxrule=1pt}}

\newtcolorbox{mybox}[2][]{enhanced,
attach boxed title to top left={xshift=5mm,yshift=-2mm},
  title={#2},#1}

\endlocaldefs

\begin{document}

%%%%%%%%%%%%%%%%%%%%%%%%%%%%%%%%%%%%%%%%%
%%                                     									       %%
%%  Start of the title page information section						       %%
%%                                     									       %%
%%%%%%%%%%%%%%%%%%%%%%%%%%%%%%%%%%%%%%%%%

\begin{frontmatter}

\begin{fmbox}

%%%%%%%%%%%%%%%%%%%%%%%%%%%%%%%%%%%%%%%%%
%%                                     									       %%
%%  Set the header on the first page, default is "Research Article"		       %%
%%                                     									       %%
%%  Define which division of the Institute of Mathematics and Physics	       %%
%%  this work is associated with								       %%
%%                                     									       %%
%%        PM:	Pure Mathematics                         					       %%
%%        ST:	Mathematical Statistics                     					       %%
%%        MP:	Mathematical Physics                        					       %%
%%        MS:	Material Physics                         					       %%
%%        FP:	Theoretical Physics							       %%
%%                                     									       %%
%%%%%%%%%%%%%%%%%%%%%%%%%%%%%%%%%%%%%%%%%

\dochead{Preprint}{FP}

%%%%%%%%%%%%%%%%%%%%%%%%%%%%%%%%%%%%%%%%%
%%                                     									       %%
%%  Title of the manuscript									       %%
%%                                     									       %%
%%%%%%%%%%%%%%%%%%%%%%%%%%%%%%%%%%%%%%%%%

\title{Exact symmetry conservation and automatic mesh refinement in discrete initial boundary value problems}

%%%%%%%%%%%%%%%%%%%%%%%%%%%%%%%%%%%%%%%%%
%%                                     									       %%
%%  Author information										       %%
%%                                     									       %%
%%%%%%%%%%%%%%%%%%%%%%%%%%%%%%%%%%%%%%%%%
\author[
   addressref={aff0},                   	  % id's of addresses
   corref={aff0},                     		  % id of corresponding address, if any
   %noteref={n1},                        		  % id's of article notes, if any
   email={alexander.rothkopf@uis.no}   		  % email address
]{\inits{AR}\fnm{A.} \snm{Rothkopf}}
\author[
   addressref={aff1},                   	  % id's of addresses
   email={wa.horowitz@uct.ac.za}   		  % email address
]{\inits{WAH}\fnm{W.\ A.} \snm{Horowitz}}
\author[
   addressref={aff2,aff3},                   	  % id's of addresses
   email={jan.nordstrom@liu.se}   		  % email address
]{\inits{JN}\fnm{J.} \snm{Nordstr{\"o}m}}

%%%%%%%%%%%%%%%%%%%%%%%%%%%%%%%%%%%%%%%%%
%%                                     									       %%
%%  Affiliations of the authors									       %%
%%                                     									       %%
%%%%%%%%%%%%%%%%%%%%%%%%%%%%%%%%%%%%%%%%%

\address[id=aff0]{%                          			 % unique id
  \orgname{Faculty of Science and Technology}, 	 % faculty
  \street{University of Stavanger},                     		 % university
  \postcode{4021},                               			 % post or zip code
  \city{Stavanger},                              				 % city
  \cny{Norway}                                   				 % country
}

\address[id=aff1]{%                          			 % unique id
  \orgname{Department of Physics}, 	 % faculty
  \street{University of Cape Town},                     		 % university
  \postcode{7701},                               			 % post or zip code
  \city{Rondenbosch},                              				 % city
  \cny{South Africa}                                   				 % country
}

\address[id=aff2]{%                          			 % unique id
  \orgname{Department of Mathematics}, 	 % faculty
  \street{Link{\"o}ping University},                     		 % university
  \postcode{SE-581 83},                               			 % post or zip code
  \city{Link{\"o}ping},                              				 % city
  \cny{Sweden}                                   				 % country
}

\address[id=aff3]{%                          			 % unique id
  \orgname{Department of Mathematics and Applied Mathematics}, 	 % faculty
  \street{University of Johannesburg},                     		 % university
  \postcode{P.O. Box 524, Auckland Park 2006},                           % post or zip code
  \city{Johannesburg},                              				 % city
  \cny{South Africa}                                   				 % country
}

%%%%%%%%%%%%%%%%%%%%%%%%%%%%%%%%%%%%%%%%%
%%                                     									       %%
%%  Additional notes on the authors								       %%
%%                                     									       %%
%%%%%%%%%%%%%%%%%%%%%%%%%%%%%%%%%%%%%%%%%

%\begin{artnotes}
%\note[id=n1]{Only contributor} % note, connected to author
%\end{artnotes}

\end{fmbox}% comment this for two column layout

%%%%%%%%%%%%%%%%%%%%%%%%%%%%%%%%%%%%%%%%%
%%                                     									       %%
%%  Abstract of the manuscript									       %%
%%                                     									       %%
%%%%%%%%%%%%%%%%%%%%%%%%%%%%%%%%%%%%%%%%%

\begin{abstractbox}
\begin{abstract} % abstract
We present a novel solution procedure for initial boundary value problems. The procedure is based on an action principle, in which coordinate maps are included as dynamical degrees of freedom. This reparametrization invariant action is formulated in an abstract parameter space and an energy density scale associated with the space-time coordinates separates the dynamics of the coordinate maps and of the propagating fields. Treating coordinates as dependent, i.e. dynamical quantities, offers the opportunity to discretize the action while retaining all space-time symmetries and also provides the basis for automatic adaptive mesh refinement (AMR). The presence of unbroken space-time symmetries after discretization also ensures that the associated continuum Noether charges remain exactly conserved. The presence of coordinate maps in addition provides new freedom in the choice of boundary conditions. An explicit numerical example for wave propagation in $1+1$ dimensions is provided, using recently developed regularized summation-by-parts finite difference operators.
\end{abstract}

%%%%%%%%%%%%%%%%%%%%%%%%%%%%%%%%%%%%%%%%%
%%                                     									       %%
%%  Keywords for the article, each one in its separate \kwd{}			       %%
%%                                     									       %%
%%%%%%%%%%%%%%%%%%%%%%%%%%%%%%%%%%%%%%%%%

\begin{keyword}
\kwd{Initial Boundary Value Problems, Space-Time Symmetry, Automatic Mesh Refinement, Noether Charge }
\end{keyword}

\end{abstractbox}

%\end{fmbox}% uncomment this for twcolumn layout

\end{frontmatter}

%%%%%%%%%%%%%%%%%%%%%%%%%%%%%%%%%%%%%%%%%
%%                                     									       %%
%%  Main text starts here										       %%
%%                                     									       %%
%%%%%%%%%%%%%%%%%%%%%%%%%%%%%%%%%%%%%%%%%

\section{Motivation}

\subsection{Executive Summary}

In this study we set out to address three central challenges in the treatment of discretized initial boundary value problems (IBVPs):
\begin{enumerate}[label*=\arabic*.]
\item the breaking of space-time symmetries related to the finite grid discretization and the associated lack of conservation of continuum Noether charges,
\item the need to construct appropriate meshes to accurately resolve the simulated dynamics, and
\item the need for more flexible and less costly implementation  of (non-reflecting) boundary conditions. %(e.g. non-reflecting boundary conditions)
%without \textit{a priori} knowledge of the field dynamics.
\end{enumerate}
We will work directly on the level of the action of the system, from which governing equations can be derived. Our novel action is formulated in a set of abstract temporal and spatial parameters $(\tau,\vec{\sigma})$. While in the usual IBVP treatment conventional space-time coordinates $(t,\vec{x})$ represent the independent parameters of the theory, here space-time coordinates are introduced as dynamical mappings $(t(\tau,\vec{\sigma}),\vec{x}(\tau,\vec{\sigma}))$, which themselves depend on abstract parameters. When discretizing the abstract parameters, the coordinate maps remain continuous and our action retains its continuum space-time symmetries after discretization. As a result, this novel discretization technique preserves the continuum Noether charges \textit{exactly}. Concurrently we show that the presence of dynamical coordinate maps leads to new contributions in boundary terms, offering novel freedom to implement boundary conditions. In addition, the dynamical interplay of the coordinate maps and the fields allows the coordinates to adapt and to form a non-trivial space-time mesh guided by the symmetries of the system, realizing a form of automatic adaptive mesh refinement.

\begin{figure}
\centering\includegraphics[scale=0.33]{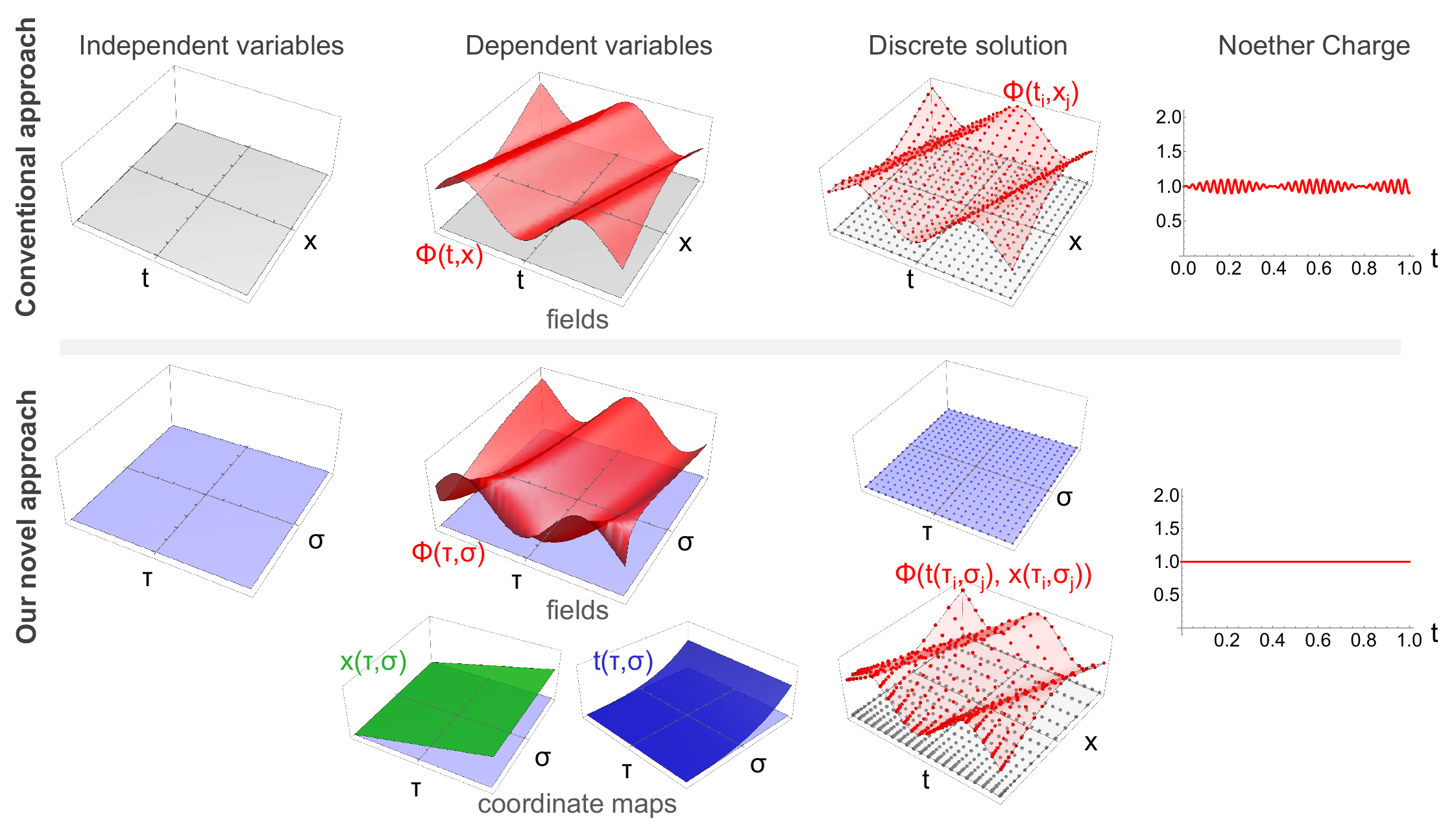}
    \caption{(Top) Sketch of the conventional approach to IBVPs: space-time coordinates are designated as independent variables and the field $\phi(t,x)$ propagates on the background of this space-time scaffold. When discretizing space-time as a hypercubic grid, the field $\phi(t_i,x_j)$ is resolved in a similarly regular fashion. Since discretized space-time is unable to accommodate infinitesimal symmetry transformations, the continuum Noether charge is not preserved. (Bottom) Sketch of our novel approach to IBVPs: A set of abstract parameters $(\tau,\sigma)$ is designated as independent parameters. Both the field $\phi(\tau,\sigma)$ as well as dynamical coordinate maps $t(\tau,\sigma)$ and $x(\tau,\sigma)$ evolve on the background of the $(\tau,\sigma)$ parameters. The evolution of the coordinate maps can be highly non-linear depending on the field dynamics. Discretizing the $(\tau,\sigma)$ parameters leaves the values of the coordinate maps continuous. In turn when expressing the physical field solution in terms of space-time coordinates $\phi(t(\tau_i,\sigma_j),x(\tau_i,\sigma_j))$, one finds in general that a non equidistant space-time grid emerges, which automatically adapts in resolution to the dynamics of the field. Since the discrete action retains its continuum symmetries the continuum Noether charge remains exactly conserved.}
    \label{fig:DifferentApproaches}
\end{figure}

A sketch of the differences between the conventional approach which discretizes space-time coordinates directly and our novel approach, which discretizes the underlying abstract parameters is provided in \cref{fig:DifferentApproaches}.

As a preview to the reader, we list the main achievements of our novel approach to IBVPs in the inset box on the next page. It showcases concrete examples from a numerical investigation of wave propagation in $(1+1)$d dimensions, treated via our novel approach, as presented in detail in \cref{sec:proofofprinciple11d}.
\pagebreak 

\centerline{\begin{minipage}[t]{1.25\linewidth}
\vspace{-1.3cm}
\begin{mybox}[width=1\linewidth]{Executive Summary - with results from (1+1)d wave propagation (c.f. \cref{sec:proofofprinciple11d})}
\begin{itemize}[leftmargin=*]
\vspace{0.3cm}
\item Space-time coordinates are promoted to dynamical degrees of freedom, which propagate together with the fields in an abstract parameter space.\\[0.1cm]
\begin{minipage}[t]{\linewidth}
\centering
 \includegraphics[scale=0.22]{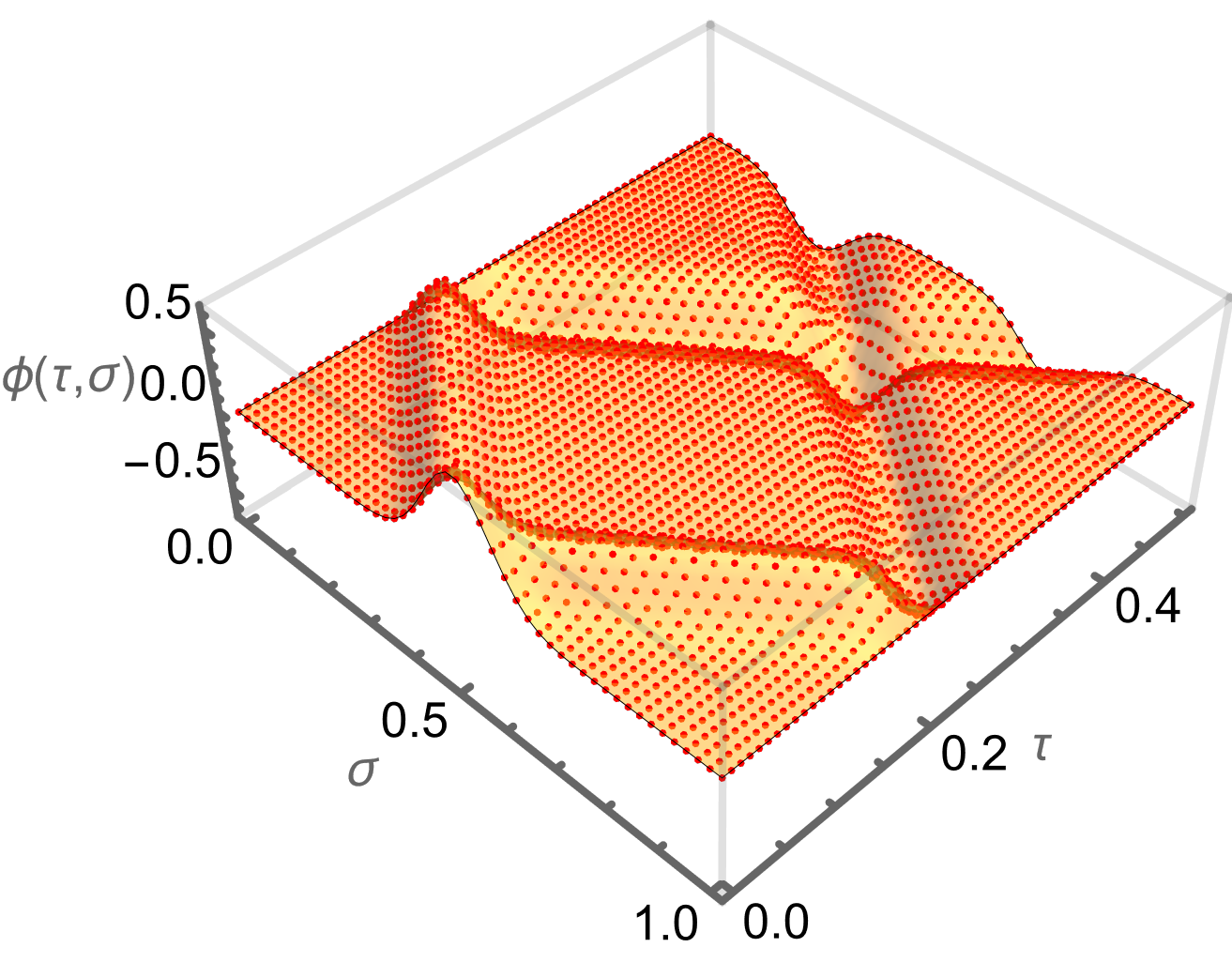}
 \includegraphics[scale=0.2]{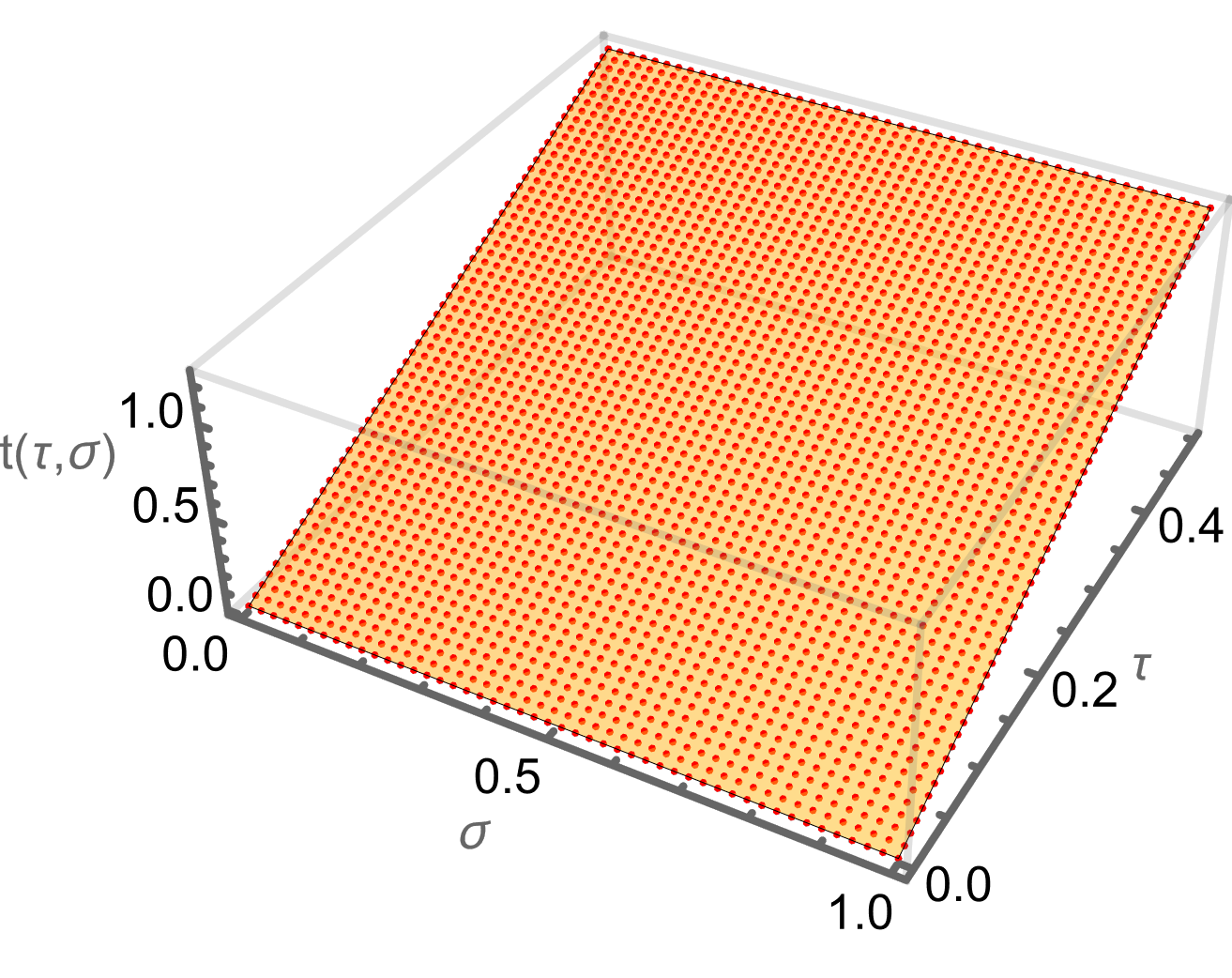}
 \includegraphics[scale=0.2]{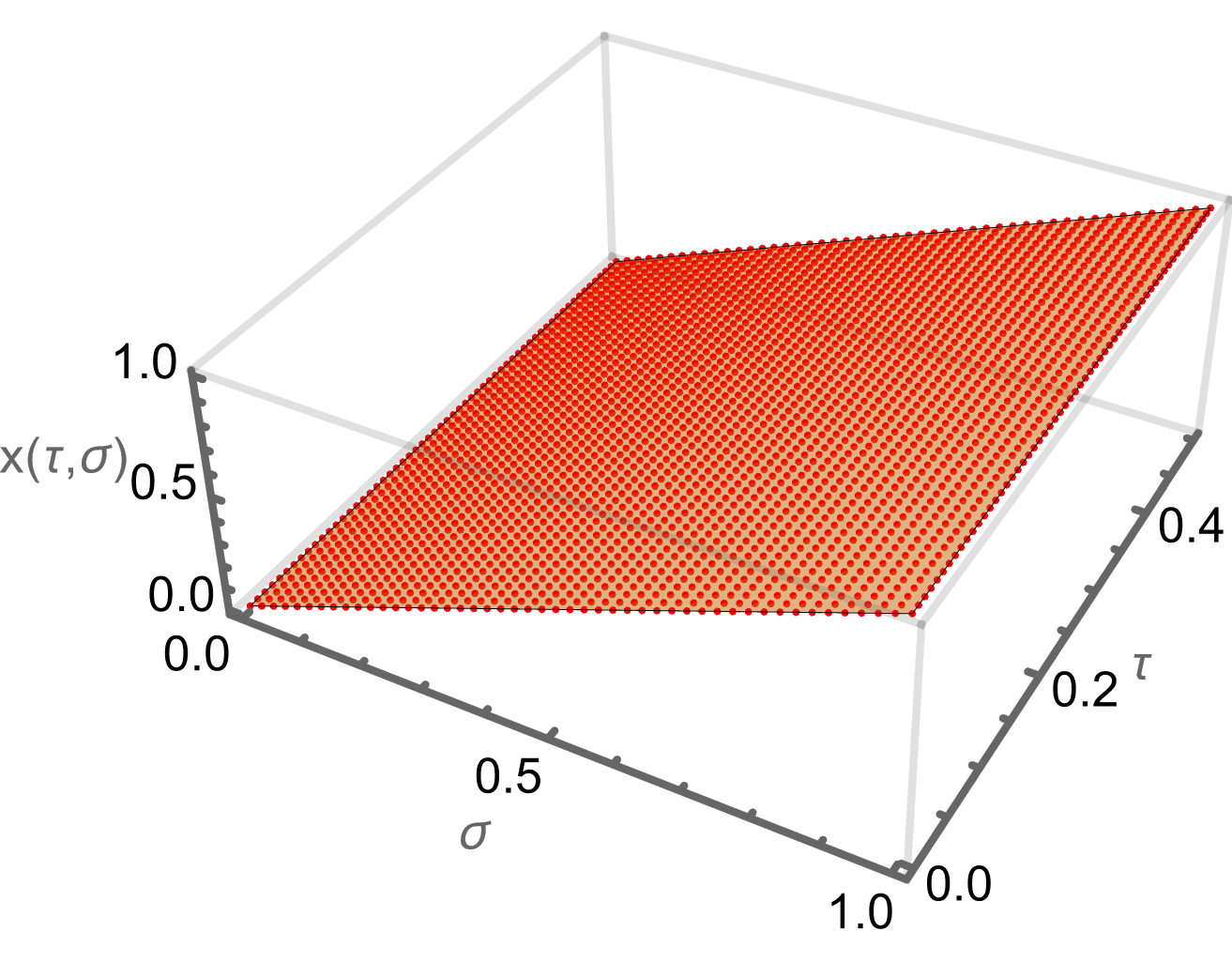}\vspace{0.1cm}
{\begin{spacing}{0.75}{\footnotesize\sffamily\raggedright  {(left) Propagating wave packages reflecting from a fixed spatial boundary. (center) Evolution of the dynamical time mapping. (right) The corresponding spatial coordinate mapping. (c.f. \cref{fig:SolutionPhiAndt} in \cref{sec:proofofprinciple11d})}} \end{spacing}}
\end{minipage}\vspace{0.05cm}
\item Interplay between fields and coordinate maps leads to a coarser or finer space-time resolution, depending on where relevant changes occur in the field configuration. This constitutes a form of dynamic resolution of the space-time coordinates which \textit{realizes automatic adaptive mesh refinement}. This mechanism is independent of the specific discretization scheme and whether the IBVP is solved on the level of the action or governing equations.\\
\begin{minipage}[t]{\linewidth}
\centering
\includegraphics[scale=0.25]{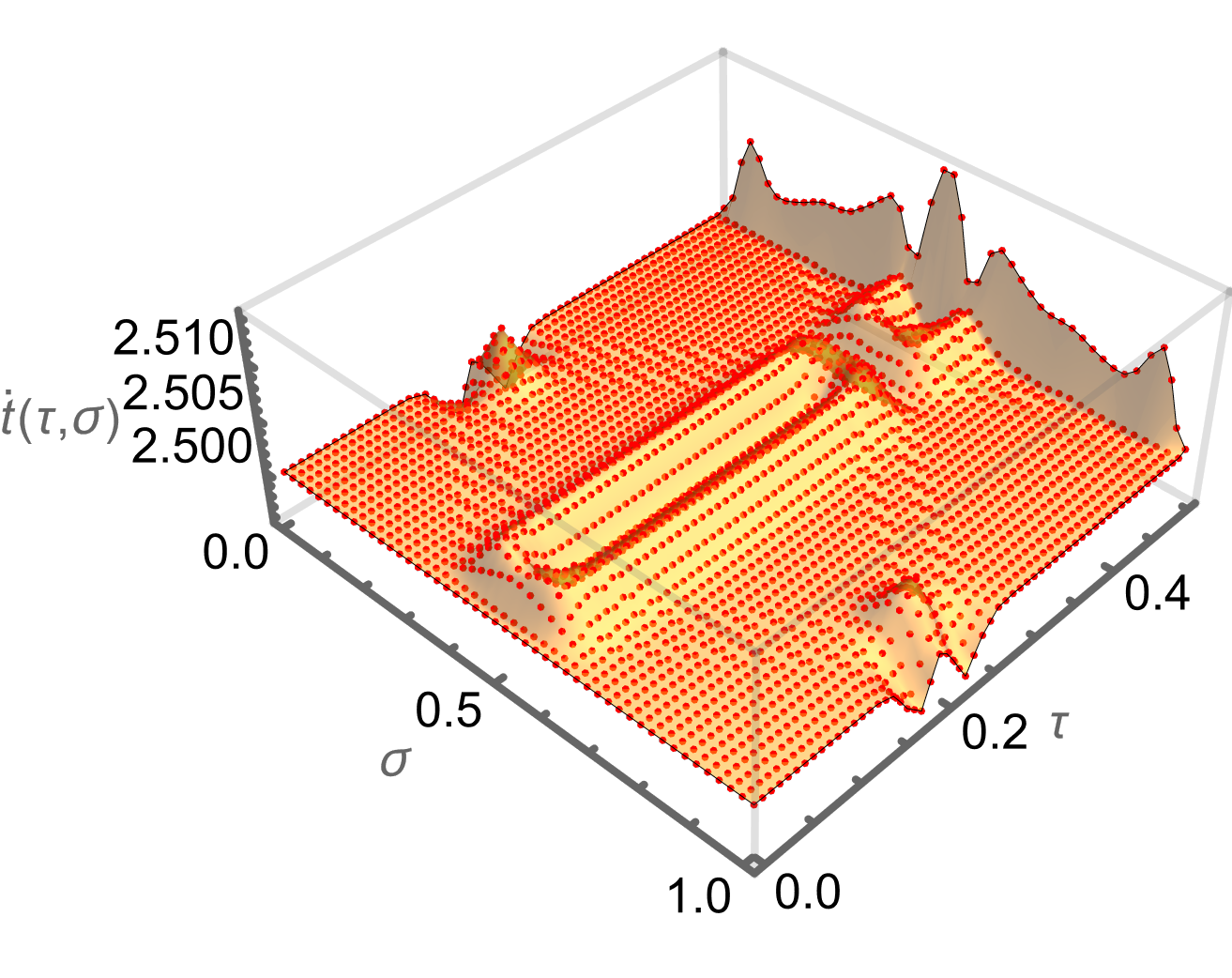}
\vspace{0.1cm}
{\begin{spacing}{0.75}{\footnotesize\sffamily\raggedright  {Non-trivial evolution of the time-mapping visualized via its temporal derivative. Larger values indicate a coarser time grid, smaller values a finer time grid. (c.f. \cref{fig:derivst} in \cref{sec:proofofprinciple11d})}} \end{spacing}}
\end{minipage}\vspace{0.05cm}
 \item Discretization in the abstract parameters instead of in space-time allows us to keep the coordinate maps continuous and thus the \textit{discretized action retains its space-time symmetries}.
\item If a discretization scheme is used that exactly mimics integration by parts, the presence of continuum space-time symmetries in the discretized action allows us to establish a discrete Noether theorem and the associated \textit{Noether charges are exactly preserved}.\\[0.2cm] 
\begin{minipage}[t]{\linewidth}
\centering
      \includegraphics[scale=0.21]{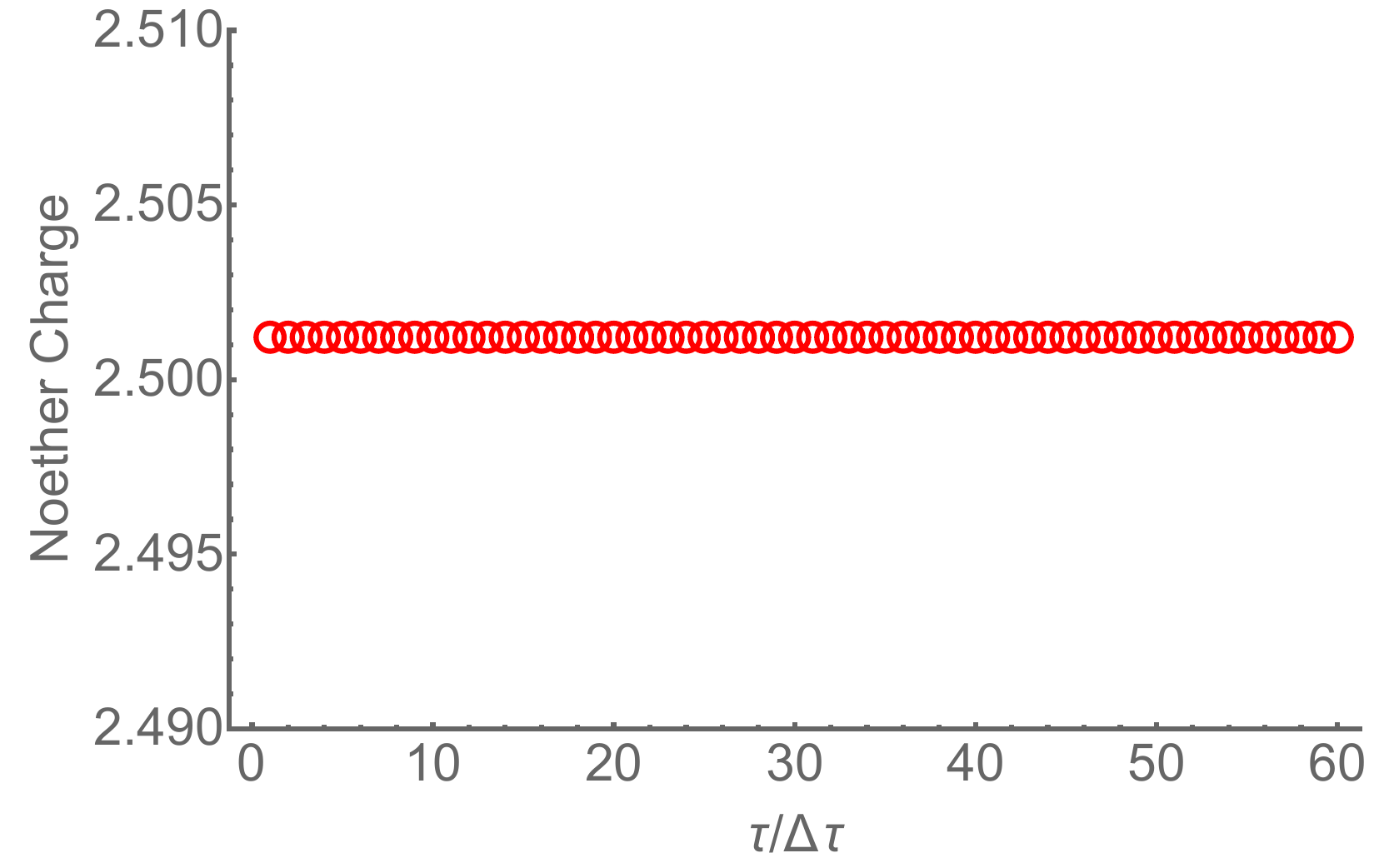}
      {\begin{spacing}{0.75}{\footnotesize\sffamily\raggedright  {The exactly preserved Noether charge associate with time translation symmetry. (c.f. \cref{fig:NoetherChargeNt60} in \cref{sec:proofofprinciple11d})}} \end{spacing}}
\end{minipage}
\item With both fields and coordinates entering as dynamical degrees of freedom, the boundary terms we encounter contain novel contributions, which provide \textit{more flexibility in the construction of boundary conditions}.
\end{itemize}
\end{mybox}
 \end{minipage}}

\subsection{Introduction}

Progress in both fundamental and applied science relies to a large extent on the solution of initial boundary value problems (IBVPs). Wave propagation is a prominent and important example. In the form of electromagnetic waves (for a textbook see e.g. \cite{taflove2005computational}) accurate numerical simulations of wave propagation are a vital ingredient in wireless telecommunication and stealth defence technology. Wave propagation is also indispensable for understanding the minute ripples in space-time emanating from the merger of two black holes in the form of gravitational waves \cite{PhysRevLett.111.241104}. In the context of mechanical waves, the propagation of sound in various media and environments has received great attention (see e.g. \cite{renterghem14}). The behavior of fluids is yet another area where IBVPs take center stage, be it non-relativistic flow (for a textbook see e.g. \cite{Lohner20081}) over airfoils and through intricate tubing or the relativistic flow (see e.g. \cite{DERRADIDESOUZA201635}) of the almost perfect fluid of hot nuclear matter created in relativistic heavy-ion collisions. 

In fact, all of the fundamental forces: gravity, electromagnetism, weak- and strong nuclear force are described in the language of a field theory \cite{landau2000classical,jackson2012classical,weinberg1995quantum, carroll2019spacetime}. Their classical dynamics are conventionally formulated in the form of IBVPs, i.e. a set of coupled partial differential equations supplied with initial and boundary data related to a specific experimental or observational scenario. From a theory point of view, the fundamental mathematical object to summarize the classical mechanics of these fields is their action \cite{goldstein2011classical}, defined as the time integral over the Lagrange functional. Dirac \cite{dirac1933lagrangian} and Feynman's \cite{RevModPhys.20.367} vital contribution to quantum field theory showed that also in the microscopic realm, the dynamics of fields is governed by this classical action. The presence  of quantum fluctuations leads to deviations from the classical configuration, but it is still the classical action that steers those quantum excursions.

The treatment of IBVPs on the level of their governing equations and through their action is intimately related. Take point mechanics as an example (see \cite{goldstein2011classical} for an in-depth introduction). We may understand the motion of particles either by using Newton's second law to analyze forces or we may gain insight by analyzing the energy budget of the system at hand. Both approaches must and do in the end lead to the same results. Governing equations are often motivated by considering the processes in a system, which lead to a change of its state, similar to applying Newton's second law in mechanics. To construct an action on the other hand, requires us to inspect the imbalance between the kinetic and potential energy contributions to our system. It is Hamilton's variational principle that relates the two approaches \cite{goldstein2011classical}. Hamilton's principle states that the classical trajectory or field configuration corresponds to a critical point of the classical action. By performing variational analysis (for a mathematical perspective see \cite{jost1998calculus}) on the action, one can show that under certain boundary conditions, the critical point of the action can also be obtained by solving a set of differential equations, the Euler-Lagrange equations. I.e. we may equivalently either solve the Euler-Lagrange equations or find the critical point of the system action. For multi-component systems, it turns out that it often requires less effort to formulate the action than to construct the corresponding governing equations.

Working with the action has a clear theoretical advantage, related to the symmetries of the system. The reason is that the action is a scalar. I.e. the degrees of freedom are combined in the action such that their non-trivial behavior under space-time symmetries exactly cancels, leaving the action invariant. In the presence of internal symmetries, the d.o.f. are combined such that the action also remains invariant under these symmetry transformations. These invariance properties are not straight forward to ascertain on the level of the governing equations, where an individual degree of freedom may change in a non-trivial fashion under space-time or internal symmetry transformations. One pertinent example are electric and magnetic fields, which transform into each other when changing between inertial frames with different relative velocities (Lorentz boosts). In contrast, the classical action of electromagnetism remains invariant under all types of Lorentz transformations. 

In this study we show that the action formulation offers additional benefits also when it comes to discretization, which is the main reason why we set out to formulate and solve classical IBVPs directly on the level of the system action. 

It is important to note that the action from introductory classical field theory is derived in the context of a boundary value problem, where the field configuration is known at both the initial and final times. Such a BVP setting is non-causal, as the final state must be known before the intermediate field configurations can be computed. In order to exploit the benefits of the action formalism for genuine initial boundary value problems, where only initial and boundary data is known, it has been shown in \cite{galley_classical_2013} and put in practice in \cite{Rothkopf:2022zfb} that one must construct an action with a doubled number of degrees of freedom. (This construction will be discussed in more detail in \cref{sec:IBVPaction} and \cref{sec:appconcond})

As described in the opening motivation, our aim is to address three central challenges associated with the discrete treatment of IBVPs: 1) loss of space-time symmetries, 2) appropriate choices of meshes and 3) the construction of boundary conditions.

The breaking of translation symmetry is one key artifact in discretizing space-time and representative of the first challenge. In the presence of a finite grid spacing, only translations of that size (or integer multiples thereof) can be implemented, prohibiting arbitrary infinitesimal translations. In turn the concept of a continuous symmetry, the basis for Noether's theorem \cite{Noether1918}, is absent and the associated Noether charges, i.e. energy, linear and angular momentum, of the system are no longer conserved. 

In some systems it may be possible to partially salvage the situation with regard to time translations, by going over to a Hamiltonian formulation, where time at first is kept continuous. But the symplectic solvers deployed in those cases only conserve energy on average (see e.g. \cite{PhysRev.159.98}). In addition, the Hamiltonian approach relies crucially on the availability of well-defined canonical momenta. Especially in systems with inherent constraints (see the discussion in \cite{dirac2001lectures}), these may be difficult to define or may require the choice of a particular gauge, as is the case with the Gauss constraint in classical electrodynamics. 

Energy conservation is intimately tied to the stability of the system dynamics (see e.g. \cite{regan2002neumann} and also \cite{nordstrom2023nonlinear}), hence non-conservation adversely affect the stability of numerical schemes. There is research ongoing to define a Noether theorem for finite grids in the context of so-called calculus of variations on time scales (see e.g. \cite{10.1063/1.5140201}). It however requires one to introduce a specific discretization scheme, which does not possess the same beneficial mimetic\footnote{A mimetic discretization exactly mimics certain continuum properties in the discrete setting.} properties as e.g. summation-by-parts (SBP) operators (for reviews see e.g. \cite{svard2014review,fernandez2014review,lundquist2014sbp}), and which appears difficult to extend to higher orders. To give Noether's theorem inherent meaning requires realizing a discretization in which infinitesimal space-time symmetry transformations can be accommodated.

The second challenge is related to the appropriate construction of meshes on which the solution of the discretized IBVP is obtained. It is intuitively clear that regions in space-time where the field evolution is more rapid require a higher density of grid points to maintain accuracy of the solution. The systematic construction of these grids is known as adaptive mesh refinement (AMR)  \cite{berger1984adaptive,lohner1987adaptive,berger1989local}, a well established field of research with a rich history. Over the past decades different strategies for the refinement of meshes have been put forward, based on three central concepts. The first is feature detection (see e.g. \cite{persson2006sub}), which deploys so-called sensors (often for shocks) to determine where violent dynamics require increased resolution. The second follows the evolution of so-called adjoint variables (see e.g. \cite{nemec2008adjoint,offermans2023error}) that can be used to improve the accuracy of the solution itself, as well as to provide an error bound, informing one in which space-time region accuracy of the solution is lost. The third approach is based on \textit{aposteriori} error estimates (see e.g. \cite{mavriplis1994adaptive,henderson1999adaptive,kompenhans2016adaptation}) which can be used to direct refinement of meshes to maintain a desired error threshold.

Each of the above strategies thus proposes a different set of criteria to identify space-time regions in which the resolution of the simulation must be modified, given a certain user defined goal for accuracy. It would be highly desirable to uncover a theoretical guiding principle behind mesh refinement, which  produces an inherent AMR criterion directly from the formulation of the IBVP.

The third challenge arises from difficulties in constructing boundary conditions for the simulation domain. In particular in the handling of propagating features, such as localized wave packets or extended wave fronts, unphysical reflections at the boundaries may and most often do contaminate the interior of the domain. Several approaches to so-called non-reflecting boundary conditions have been developed, among them complete radiation boundary conditions (CRBC) \cite{Complete,CRBCel} and perfectly matched layers (PML) \cite{johnson2021notes}. %True transparency at the boundaries may require that the solution must already be known or an ill-posed inverse problem ensues. 
Partially non-reflecting boundaries may be constructed, such that they allow wave-fronts traveling in a certain direction to pass. While it is already possible to treat systems of waves with different dispersion relations (see e.g. \cite{AppCol09}), so far no conclusive solution to the challenge of non-reflecting boundaries has been obtained. In order to support a comprehensive resolution of this challenge, formulations that provide additional flexibility of boundary treatment are called for.

In order to address the first two of the three above mentioned challenges, our study puts forward two central contributions.

The \textit{first central contribution} of this work is the construction of a novel continuum action for second order IBVPs in $(d+1)$ dimensions. Instead of being formulated directly in space-time, we introduce a $(d+1)$ dimensional abstract parameter space in which all degrees of freedom, including space-time coordinate maps, propagate. These coordinate maps translate parameter coordinates into space-time coordinates and enter the action together with the propagating fields. 

The presence of dynamical coordinate maps leads to an interplay with the fields, which determines how fast space-time progresses in terms of the abstract coordinates. The fact that coordinate maps are now dynamical means that they too appear in the boundary terms one encounters, when deriving the governing equations from the action. Their presence thus offers new flexibility in how to implement consistent boundary conditions and thus offers new possibilities for developing non-reflecting boundary conditions. 

The \textit{second central contribution} is a strategy to discretize the novel action in the abstract parameters and not in space-time coordinates. To this end we deploy appropriately regularized multidimensional summation-by-parts finite difference operators. This discretization takes place on the level of the abstract parameters with a naive hypercubic mesh\footnote{The discretization performed here constitutes a multidimensional generalization of the discretization strategy developed for the world-line approach in \cite{Rothkopf:2023ljz}.}. Since the values of the coordinate maps remain continuous, the discretized action also retains its invariance under all infinitesimal continuum space-time symmetries. Explicitly preserving the space-time symmetries in this way ensures that Noether's theorem remains fully valid in the discrete setting and that the associated Noether charges remain conserved.  

The interplay between fields and coordinate maps in the novel action leads to a self-consistent regulation of how space-time flows. This furnishes the basis for achieving automatic AMR after discretization. We will demonstrate in \cref{sec:proofofprinciple11d}, using scalar wave propagation in $(1+1)$ dimensions as example, that a discretization on the level of the abstract parameter space leads to a non-trivial evolution of the coordinate maps, coupled to the field dynamics. This in turn defines a non-trivial mesh of space-time coordinates. It turns out that our procedure leads to coordinate maps that evolve such that space-time is resolved more finely in regions of rapid dynamics, while space-time resolution is coarser where the dynamics of the fields are less rapid. 

The rest of the paper is structured as follows. In \cref{sec:contaction} we derive the novel reparameterization invariant continuum action for fields in $(d+1)$ dimensions. We start out in \cref{sec:BVPaction} in the context of the conventional boundary value setting of classical field theory, to highlight the main novelty of our action, which is the presence of dynamical coordinate maps. In \cref{sec:IBVPaction} we proceed to formulate the novel action with doubled degrees of freedom, necessary to describe a genuine causal IBVP. In preparation for discretizing the action, we discuss in \cref{sec:lagmultcont} how initial and boundary data can be explicitly included in the action via Lagrange multipliers and how this affects the Noether charges of the system. \Cref{sec:discrIBVP} lays out the details of our discretization strategy for the novel action based on appropriately regularized summation-by-parts finite difference operators. All theoretical development in the preceding sections will be put to work in \cref{sec:proofofprinciple11d}, addressing the concrete example of scalar wave propagation in $(1+1)$ dimensions. We provide the explicit form of the continuum action, discretize this action and solve for the classical solution by numerically locating the critical point of the IBVP action, i.e. without solving the governing equations. The convergence properties of the discretization strategy will be assessed.

\section{Continuum action with dynamical coordinate maps}
\label{sec:contaction}
In this section we construct a reparameterization invariant action for classical field theory with dynamical coordinate maps, inspired by the world-line formalism of relativistic point mechanics (see, e.g. \cite{landau2000classical,carroll2019spacetime}). 

For the world-line formalism, the theory of general relativity and its guiding principle of reparameterization invariance provides us with an action that naturally contains dynamical coordinate maps. %(see \cref{eq:geodesic,eq:geodesicincl} in \cref{sec:BVPaction}). 
Based on this starting point, it is possible to explore (as was done in \cite{Rothkopf:2023ljz}) the consequences of the presence of coordinate maps, which even after discretization allows for preservation of space-time symmetries. While in string theory actions with explicit coordinate maps have been studied for a long time, e.g. in the context of Born-Infeld electromagnetism on D-branes \cite{Taylor:2002uv}, field theory actions are predominantly formulated with reference to space-time coordinates. Our goal in this section is thus to reverse engineer an action for a wide variety of IBVPs in the presence of dynamical coordinate maps, using the world-line formalism as guide. This constitutes a genuinely novel approach to IBVPs not explored in the literature so far.

Conventionally, the variational principle in classical mechanics and field theory is formulated as a boundary value problem (BVP). One assumes that the values of the propagating degrees of freedom are known at initial time and final time and the critical point of the classical action provides us with the classical trajectory or field configuration spanning between these time slices. As discussed in \cite{galley_classical_2013} and more recently in \cite{Rothkopf:2022zfb}, such a BVP is not causal, as the final state of the system must be known before the trajectory can be determined. In order to solve a genuine initial boundary value problem, a different variational treatment is called for, which involves a doubling of the degrees of freedom of the system. 

In the following subsection we begin the construction of our novel action with dynamical coordinate maps in the context of classical field theory as a boundary value problem. This allows us to present to the reader the novel ingredients in our approach, before implementing a genuine IBVP formulation in the subsequent subsection. 

\subsection{Boundary Value Problem formulation}
\label{sec:BVPaction}
The starting point for our construction of a reparametrization invariant action with dynamical coordinate maps is an analogy with the relativistic dynamics of point particles in general relativity. 

Let us briefly review the motion of a point particle under the influence of arbitrary forces, encoded in the potential function $V(x)$ (see \cite{Rothkopf:2022zfb} for a more detailed discussion). The standard non-relativistic action, applicable for motion at velocities much slower than the speed of light reads
\begin{align}
    S_{\rm nr}=\int \, dt\, L_{\rm nr}[\vec{x},\dot{\vec{x}}] = \int \, dt\, \Big\{ \frac{1}{2} m \dot \vec{x}^2(t) - V\big(\vec{x}(t)\big) \Big\}, \label{eq:nonrelac}
\end{align}
where $L_{\rm nr}$ refers to the classical Lagrangian, which consists of the difference between kinetic and potential energy.
In contrast, the world-line formalism of general relativity posits that in order to correctly describe the physics also at velocities close to those at the speed of light, one must consider the following action
\begin{align}
    S_{\rm wl}&=\int \, d\gamma\,\Big\{  (-mc)\sqrt{ G_{\mu\nu}\frac{d X^\mu}{d\gamma} \frac{dX^\nu}{d\gamma} } - V(X) \Big\}. \label{eq:geodesic}
\end{align}
Summation over repeated indices is implied. The kinetic term of this action is derived following the guiding principle of reparameterization invariance, central to the general theory of relativity\footnote{In order to show that the potential term too obeys reparameterization invariance, one has to start from the so-called Polyakov action with an auxiliary field $\xi$, which allows one to couple mass and potential terms in a manifestly reparameterization invariant fashion. The Polyakov action can be shown to be classically equivalent to the world-line action \cite{zwiebach2004first}.}. The action $ S_{\rm wl}$ describes the motion of a point particle in a space-time characterized by the metric tensor $G$\footnote{\label{sec:fnG}The metric tensor $G$ describes how length and time intervals, as well as angles are defined in a given space-time. In flat $(3+1)$ dimensional space-time $G\in\mathbb{R}^{4 \times 4}$ is often denoted as $\eta_{\mu\nu}\equiv G_{\mu\nu}$, has negative determinant (we choose here the signature convention $(+,-,-,-)$ \cite{goldstein2011classical}) and is diagonal with the entries $G_{\mu\nu}={\rm diag}[G_{00},G_{11},G_{22},G_{33}]={\rm diag}[c^2,-1,-1,-1]$. The inverse of $G_{\mu\nu}$ is denoted by $G^{\mu\nu}$.}, but instead of time $t$, uses an abstract world-line parameter $\gamma$ to parameterize the one-dimensional trajectory $X^\mu(\gamma)=(t(\gamma),\vec{x}(\gamma))^\mu=(t(\gamma),x_1(\gamma),\ldots,x_d(\gamma))^\mu$ that the particle traces out in $(d+1)$ dimensional space-time. Here the superscript $\mu=0,\ldots,d$ refers to the $\mu$-th components of the vector $X$ and is used in the following to implement the Einstein summation convention. As discussed in \cite{Rothkopf:2022zfb}, for a diagonal metric and a potential, which is small compared to the rest energy of the point particle $V(\vec{x})/2mc^2\ll1$, one may rewrite the world-line action \cref{eq:geodesic} as a genuine geometric, i.e. a geodesic action
\begin{align}
    S_{\rm wl}\approx S_{\rm geo} = \int \, d\gamma\,(-mc)\Big\{  \sqrt{ \Big(G_{00}+\frac{V(\vec{x})}{2mc^2}\Big)\frac{d X^0}{d\gamma} \frac{dX^0}{d\gamma}+ G_{ii}\frac{d X^i}{d\gamma} \frac{dX^i}{d\gamma} } \Big\}.\label{eq:geodesicincl}
\end{align}

When taking the non-relativistic limit in flat space-time $G=\eta$ (see \cref{sec:fnG}), which builds on the separation of two scales $V(\vec{x})/2mc^2\ll 1$ and $\frac{d|\vec{x}|/d\gamma}{dt/d\gamma}/c = v/c \ll 1$, one is able to recover a non-relativistic action from \cref{eq:geodesicincl}. This non-relativistic action features one decisive difference with respect to the usual \cref{eq:nonrelac}
\begin{align}
    S_{\rm nr} \overset{\rm from \; \cref{eq:geodesicincl}}{=} \int \, dt\, \Big\{-mc^2 +  \frac{1}{2} m \dot{\vec{x}}^2(t) - V\big(\vec{x}(t)\big) \Big\},\label{eq:nonrelacderiv}
\end{align}
which is the presence of the constant rest mass term $-mc^2$. In addition one also learns that the factor $1/2$ in front of $m\dot{\vec{x}}^2$ is directly related to the expansion of the square root in the world-line action.

Note that a constant term in the action does not affect the classical trajectory. Indeed, a constant term vanishes under variation of the action, and does not contribute to the stationary condition $\delta S_{\rm nr}=0$ defining the trajectory. This explains why one cannot guess the ``correct value" of the constant term without an additional guiding principle. In the case of general relativity, the guiding principle is reparametrization invariance, which leads to the construction of the world-line action, and the non-relativistic limit $v\ll c$ is a systematic approximation of that world-line action. In particular, this expansion in $v/c$ leads to a series of correction terms to the non-relativistic action \cref{eq:nonrelac}. Since the effect of these corrections only becomes relevant in the genuine relativistic regime, it was historically difficult in experiment to spot that \cref{eq:nonrelac} is missing relevant physics. %For example, a machine learning algorithm applied to experimentally available data up to the beginning of the 20th century would not have been able to infer the inadequacy of the non-relativistic action.

Keeping the above line of arguments in mind, let us turn to genuine field theory. The fully relativistic field theory action found in standard textbooks (see e.g. chapter 4 of \cite{carroll2019spacetime}) includes both \textit{kinetic} terms, referring to the space-time derivatives of the field, as well as (self-)\textit{interaction} terms in the form of an in general non-linear potential function $V(\phi)$. This action too is reparameterization invariant and is conventionally written as
\begin{align}
    S=\int d^{(d+1)}X \, \sqrt{-{\rm det}[G]} \frac{1}{2}\Big( G^{\mu\nu} \partial_\mu\phi(X)\partial_\nu\phi(X) - V(\phi) \Big). \label{eq:stdFTaction}
\end{align}
For notational convenience, we took the liberty to define the self-interaction function $V(\phi)$ such that an overall factor $1/2$ is taken to the front (i.e. a mass term quadratic in the fields would now read $V_{\rm mass}=m\phi^2$). The above action describes second order systems, which feature a kinetic term quadratic in the velocities $\partial_\mu\phi$, as is the case for e.g. scalar wave propagation \cite{Lohner20081} (for which $V(\phi)=0$) or self interacting Klein-Gordon fields with $V(\phi)\neq0$.

To connect the action of \cref{eq:stdFTaction}, written in an abstract relativistic form, to the familiar physics of, e.g., scalar wave propagation, let us consider a flat space-time, where ${\rm det}[G]=-c^2$, and let us set $V(\phi)=0$. We then recover the standard action for scalar wave propagation
\begin{align}
    S^{\rm wp}=\int d^{(d+1)}X\, c\,\underbracket{\Big( \frac{1}{c^2}\,\frac{1}{2}  \partial_0\phi(X)\partial_0\phi(X) - \frac{1}{2}\partial_i\phi(X)\partial_i\phi(X)  \Big)}_{{\cal L}^{\rm wp}}. \label{eq:stdWPaction}
\end{align}
The factor $c$ up front in the integrand ensures that the units of the integration measure match those of the derivatives in the parentheses. The variational principle tells us that the classical solution is obtained at the critical point of this action. After computing the stationary condition 
\begin{align}
  &0=\delta S^{\rm wp}= \int d^{(d+1)}Xc \Big\{ \frac{1}{c^2} (\partial_0\phi)\delta(\partial_0\phi) - (\partial_i\phi)\delta(\partial_i\phi)\Big\}\\
  \nonumber &\overset{\rm IBP}{=}\int d^{(d+1)}Xc \Big\{ -\Big( \frac{1}{c^2}  \partial^2_0\phi(X) - \partial^2_i\phi(X)  \Big)\delta\phi \Big\} + \frac{1}{c^2}(\partial_0\phi \,\delta\phi)|_{t^{\rm i}}^{t^{\rm f}} + \sum_j  (\partial_j\phi\, \delta \phi)|_{x_j^{\rm i}}^{x_j^{\rm f}}
\end{align}
one finds that only if the spatial boundary terms vanish (the temporal ones in the BVP setting do, since $\delta \phi$ vanishes), one is led to the following governing (Euler-Lagrange) equation
\begin{align}
    \frac{\partial {\cal L}^{\rm wp}}{\partial  \phi(X)} - \partial_\mu \Big[ \frac{\partial {\cal L}^{\rm wp}}{\partial (\partial_\mu \phi(X))}\Big]= -\Big( \frac{1}{c^2}  \partial^2_0\phi(X) - \partial^2_i\phi(X)  \Big)=0, \label{eq:stdWPgoveq}
\end{align}
i.e. the well-known scalar wave equation with wave propagation speed $c$.

The standard action in \cref{eq:stdFTaction} remains invariant under arbitrary differentiable reparameterizations of the space-time $X\to X^\prime$ due to the presence of the term $\sqrt{-{\rm det}[G]}$. This term absorbs the Jacobian $J=\partial X/\partial X^\prime$ of the coordinate transformation to yield $\sqrt{-{\rm det}[J^TGJ]}=\sqrt{-{\rm det}[G^\prime]}$, which defines the new \textit{induced metric} $G^\prime$ in the new coordinate system\footnote{The metric $G$ is associated with the inner product of four-dimensional space-time vectors, which in turn defines distances between events taking place in the space-time described by the coordinates $X$. The induced metric $G^\prime$ implements the corresponding inner product in the coordinate system defined by $X^\prime$. For more details see, e.g., \cite{zwiebach2004first}.}. The action is also invariant under the full Poincare group of space and time translations, spatial rotations and relativistic rotations involving space and time, called boosts. 

After inspecting the properties of the standard action, let us take the first step towards \textit{reverse engineering} a field theory action with dynamical coordinate maps. Taking insight from the relation between the non-relativistic point-particle action \cref{eq:nonrelacderiv} and the geodesic action \cref{eq:geodesicincl}, we add to \cref{eq:stdFTaction} a constant term $-T$ akin to the rest mass term $-mc^2$ (or tension $T$ of a D-brane, see discussion on page \pageref{par:tension}) . Since the extremum of $S_{\rm ft}$ is unchanged by the constant, the classical physics too remains unchanged
\begin{align}
    S=\int d^{(d+1)}X \, \sqrt{-{\rm det}[G]} \Big\{ -T + \frac{1}{2}\Big(G^{\mu\nu} \partial_\mu\phi(X)\partial_\nu\phi(X) + V(\phi) \Big)\Big\}.\label{eq:stdFTactionpconst}
\end{align}
Note that in order for the constant $T$ to be accommodated in this action, it must carry the units of an energy density.

For the next step, we interpret the presence of the constant term and the factor $1/2$ as an indication that we are actually dealing with the low energy limit of a more general action (akin to the world-line action) which contains a square root term. I.e. we assume that in \cref{eq:stdFTactionpconst} small terms at higher order in inverse powers of some dimensionless ratio $\kappa\sim{\rm energy\, density}/T$ have been discarded. As part of our reverse engineering strategy, we undo this expansion by introducing a square root as
\begin{align}
    \nonumber&{\cal S}_{\rm BVP}\equiv\\
    &\int d^{(d+1)}X \, \sqrt{-{\rm det}[G]} \big(-T\big) \Big\{ 1 - \frac{1}{2T}\Big(G^{\mu\nu} \partial_\mu\phi(X)\partial_\nu\phi(X) + V(\phi) \Big)\Big\} + {\cal O}(\kappa^2),\\
    &=\int d^{(d+1)}X \, \sqrt{-{\rm det}[G]} \big(-T\big) \sqrt{ 1 - \frac{1}{T}\Big(G^{\mu\nu} \partial_\mu\phi(X)\partial_\nu\phi(X) + V(\phi) \Big)}.  \label{eq:stdFTactionbacksqrt}
\end{align}
In order for the previous step to be justified the ratio $\kappa$ must indeed be small; conversely, $T$ must correspond to a large energy density compared to the energy density represented by the field $\phi$. This situation mirrors the world-line formalism, where the rest-mass of a point particle is much larger than the kinetic and potential energy encountered on everyday scales.

Now we take the next vital step, which consists of elevating the space-time coordinates $X^\mu$ to dynamical degrees of freedom. This too mirrors the world-line formalism, where one elevates $t\to t(\gamma)$ and it provides the mechanism that will subsequently enable the preservation of space-time symmetries after discretization.

To this end we change coordinates in the integral from space-time $X^\mu=(t, \vec{x})^\mu=(t,x_1,\ldots,x_d)^\mu$ (in \cref{eq:stdFTactionbacksqrt}) to a set of arbitrary parameters $\Sigma^a=(\tau,\vec{\sigma})^a=(\tau,\sigma_1,\ldots,\sigma_d)^a$ with $a=0,\ldots,d$ (in \cref{eq:stdFTactionbacksqrt2}). The two coordinate systems are related by an invertible and differentiable mapping $X^\mu(\Sigma)$, which is explicitly encoded via the Jacobian $J^\mu_a=\partial X^\mu/\partial \Sigma^a$ that appears in the measure as we carry out the change of variables
\begin{align}
    {\cal S}_{\rm BVP}=\int d^{(d+1)}\Sigma \, |{\rm det}[J]|&\sqrt{-{\rm det}[G]} \big(-T\big)\label{eq:stdFTactionbacksqrt2}\\
    \nonumber&\times \sqrt{ 1 - \frac{1}{T}\Big(G^{\mu\nu} \partial_\mu\phi(X(\Sigma))\partial_\nu\phi(X(\Sigma)) + V(\phi) \Big)}.  
\end{align}
We will take the range of the parameters $\Sigma^a$ to be finite and to span a hypercube interval from their lowest $(\tau^{\rm i},\sigma_1^{\rm i},\ldots,\sigma_d^{\rm i})$ to highest values $(\tau^{\rm f},\sigma_1^{\rm f},\ldots,\sigma_d^{\rm f})$. I.e. we will always refer to initial and final values of the parameters with superscript letters ${\rm i}$ and ${\rm f}$ respectively.

Let us inspect the \textit{induced} metric $g_{ab}$ on the parameter manifold, defined as the quantity $g_{ab}\equiv G_{\mu\nu} J^\mu_a J^\nu_b$. The induced metric naturally arises from combining the determinant of the Jacobian with the metric dependent contribution of the measure 
\begin{align}
&\sqrt{-{\rm det}[J]{\rm det}[G]{\rm det}[J]}\\
\nonumber &\qquad\qquad=\sqrt{-{\rm det}[J^T]{\rm det}[G]{\rm det}[J]}=\sqrt{-{\rm det}[J^T G J]}= \sqrt{-{\rm det}[g]},    
\end{align}
so that the action reads
\begin{align}
     {\cal S}_{\rm BVP}=\int d^{(d+1)}\Sigma \, \sqrt{-{\rm det}[g]}& \big(-T\big)\label{eq:stdFTactionbacksqrtmod1}\\
     \nonumber&\times \sqrt{ 1 - \frac{1}{T}\Big(G^{\mu\nu} \partial_\mu\phi(X(\Sigma))\partial_\nu\phi(X(\Sigma)) + V(\phi) \Big)}.  
\end{align}

If we also absorb the determinant ${\rm det}[g]$ into the main square root expression in the action, we obtain
\begin{align}
    {\cal S}_{\rm BVP}&=\int d^{(d+1)}\Sigma \, \big(-T\big) \sqrt{ -{\rm det}[g] + \frac{1}{T}\Big(G^{\mu\nu} \partial_\mu\phi(\Sigma)\partial_\nu\phi(\Sigma) + V(\phi) \Big){\rm det}[g]}.  \label{eq:stdFTactionbacksqrtmod2}
\end{align}

What we have achieved so far is to formulate an action for fields (\cref{eq:stdFTactionbacksqrtmod2}) in an arbitrary parameter space, described by the coordinates $\Sigma^a$. The maps $X(\Sigma)$ that relate coordinates in parameter space to space-time coordinates enter in the ${\rm det}[g]$ terms in the square root. Note that the field dependent term on the right in \cref{eq:stdFTactionbacksqrtmod2} however still contains reference to the metric $G$ and derivatives in terms of the old coordinates $\partial_\mu=\partial/\partial X^\mu$. In order to complete the transition from $X$ to $\Sigma$ coordinates, we have to rewrite the RHS of \cref{eq:stdFTactionbacksqrtmod2}  solely in terms of derivatives in the new coordinates $\partial_a=\partial/\partial \Sigma^a$.

To this end let us recognize that the contraction of the partial derivatives $G^{\mu\nu}\partial_\mu\partial_\nu$ is invariant under reparameterization. It thus follows that we can simply replace it by a contraction in terms of the new derivatives  $\partial_\mu\phi \partial_\nu\phi G^{\mu\nu}=\partial_a\phi \partial_b\phi g^{ab}$. We can confirm this replacement by explicit calculation, writing in component notation $(\partial^a\phi J_a^\mu) (\partial^b\phi J_b^\nu) G_{\mu\nu}$ with $J^a_\mu J^b_\nu G^{\mu\nu}=g^{ab}$ by definition of $g$ and the fact that $J^a_\mu=(J^{-1})^\mu_a$. Let us also move the potential term to the left under the square root, then
\begin{align}
    {\cal S}_{\rm BVP}&=\int d^{(d+1)}\Sigma \, \big(-T\big) \sqrt{ \Big(\frac{1}{T}V(\phi)-1\Big){\rm det}[g] + \frac{1}{T} \partial_a\phi(\Sigma)\partial_b\phi(\Sigma) g^{ab}  {\rm det}[g]}.  \label{eq:stdFTactionbacksqrtmod3}
\end{align}

Note that both the inverse (denoted by $g^{-1}=g^{ab}$) and the determinant ${\rm det}[g]$ of the induced metric $g_{ab}$ appear in the $\phi$ dependent term on the right under the square root. As long as $g$ is invertible, we can further simplify the action by applying a well-known matrix identity, involving the adjugate of the metric, ${\rm adj}[g]=g^{-1}{\rm det}[g]$. This leads us to the final expression for our novel reparameterization invariant action with dynamical coordinate maps
\begin{align}
\hspace{-0.8cm}\tcbhighmath[boxrule=2pt,arc=1pt,colback=gray!10!white,colframe=gray,
  drop fuzzy shadow=gray]{{\cal S}_{\rm BVP}=\int d^{(d+1)}\Sigma \, \big(-T\big) \sqrt{ \Big(\frac{1}{T}V(\phi)-1\Big){\rm det}[g] + \frac{1}{T} \partial_a\phi(\Sigma)\partial_b\phi(\Sigma) {\rm adj}[g]_{ab}}.} \label{eq:novelaction}
\end{align}
Summation of the indices $a$ and $b$ is implied. Since no additional factors of $g$ enter this summation, we write both indices as subscripts. We emphasize that the induced metric $g$ makes reference to $X^\mu$ and $\partial_aX^\mu$. The corresponding Lagrangian ${\cal L}_{\rm BVP}$ is defined as ${\cal S}_{\rm BVP}\equiv\int d^{(d+1)}\Sigma\,{\cal L}_{\rm BVP}[X,\partial_a X,\phi,\partial_a \phi]$.

This novel action \cref{eq:novelaction} constitutes a generalisation of the world-line action for a point particle to field dynamics. Thus we may expect some similarities but also differences between \cref{eq:novelaction} and \cref{eq:geodesicincl}. 

Note first that our new action is manifestly invariant under Poincar\'e transformations. These transformations only affect the coordinate maps $X^\mu$, which appear in the form of the determinant or adjugate of the matrix $g$. As the coordinate maps enter the induced metric $g$ by contraction of their derivatives via the metric tensor $G$, the entries of $g$ by themselves are invariant under Lorentz transformations. And since it is only the derivatives of $X^\mu$ that appear, any constant shift in the coordinate maps is annihilated. In turn the entries of $g$ and hence ${\cal S}$ are invariant under the full Poincar\'e group of transformations.

Note also the overall multiplicative factor of $-T$ in \cref{eq:novelaction}, the analog to $-mc$ for the world-line formalism. This quantity in both cases denotes the scale at which the full dynamical nature of the coordinates becomes relevant. For a point particle, at kinetic energies close to its rest energy, we must take into account that its motion through time can no longer be considered independent from its motion through space, contrary to the non-relativistic limit. Similarly here, once the energy density carried by the field becomes of the order of $T$, the dynamics of the coordinate mappings $X^\mu$ encoded in $g$ are inescapably intertwined with those of the field $\phi$. 

In the geodesic action \cref{eq:geodesicincl}, the effects of a potential were approximately absorbed into a term related to the dynamics of the coordinate $t$, promoted to a dynamical degree of freedom. In the case of field theory in \cref{eq:novelaction}, the potential function $V(\phi)$ is naturally associated with the ${\rm det}[g]$ term under the square root, which governs the newly introduced dynamical coordinate maps.

Let us consider the physics of the coordinate maps in the absence of a field, by taking the $T\to\infty$ limit. Interestingly, in that case, we recover an action given solely by the square root of the determinant of the induced metric. When $d=1$, this expression plays a central role in string theory, where it is known as the Nambu-Goto action (for details see \cite{zwiebach2004first}) and describes a spatially extended string that propagates, tracing out a so-called world-sheet. In an otherwise flat space-time, the string world-sheet, as the name suggest, constitutes a flat sheet with constant slope. For $d>1$, this expression is known as a D-brane action (see, e.g., \cite{Green:1996bh}).

\label{par:tension} In our case the action includes reference to all $(d+1)$ space-time coordinates $X$ and we consider a $(d+1)$-dimensional world-volume that is traced out, as the underlying parameters $\Sigma$ vary. The constant $T$ within string theory has a clear physics interpretation: it denotes the tension $T/c$ of the string, and constitutes a fundamental property of nature. Similarly we may interpret $T$ here as a reference energy density that tells us how much energy is required in order to produce a stretching or compression of the coordinate mappings.

The D-brane action is invariant under coordinate reparameterizations and by construction, so is our novel action, including the field degrees of freedom.

The variational principle of classical field theory tells us that the classical field configuration is found as the critical point of the system action. Furthermore this critical point can either be found by direct variation of the action itself or equivalently by solving a set of partial differential equations of motion, the so-called Euler-Lagrange equations. 

To obtain the Euler-Lagrange equations of motion, let us proceed by carrying out the variation of \cref{eq:novelaction} in all the dependent variables $X^\mu$, $\partial_a X^\mu$, $\phi$ and $\partial_a \phi$. 
\begin{align}
\nonumber&\delta {\cal S}_{\rm BVP}[X, \partial_a X, \phi, \partial_a \phi]\\
\nonumber &=\int d^{(d+1)}\Sigma \Big\{ 
\frac{\partial {\cal L}_{\rm BVP}}{\partial X^\mu}\delta X^\mu + \frac{\partial {\cal L}_{\rm BVP}}{\partial (\partial_a X^\mu)}\delta (\partial_a X^\mu) + \frac{\partial {\cal L}_{\rm BVP}}{\partial \phi}\delta \phi + \frac{\partial {\cal L}_{\rm BVP}}{\partial (\partial_a \phi)}\delta (\partial_a \phi) \Big\}\\
&=\int d^{(d+1)}\Sigma \Big\{ 
\Big( \frac{\partial {\cal L}_{\rm BVP}}{\partial X^\mu} -\partial_a \frac{\partial {\cal L}_{\rm BVP}}{\partial (\partial_a X^\mu)}\Big)\delta X^\mu + \Big(  \frac{\partial {\cal L}_{\rm BVP}}{\partial \phi} - \partial_a \frac{\partial {\cal L}_{\rm BVP}}{\partial (\partial_a \phi)} \Big)\delta \phi \Big\} \label{eq:varSeom}\\
&\qquad + \sum_{n=0}^{d} \int \prod_{\tiny\begin{array}{c}a=0\\a\neq n\end{array}}^{d}d\Sigma_a \Big[ \frac{\partial {\cal L}_{\rm BVP}}{\partial (\partial_n X^\mu)} \delta X^\mu \Big]_{\Sigma_n^{\rm i}}^{\Sigma_n^{\rm f}}\label{eq:varSbnd1}\\ 
&\qquad + \sum_{n=0}^{d} \int \prod_{\tiny\begin{array}{c}a=0\\a\neq n\end{array}}^{d}d\Sigma_a \Big[ \frac{\partial {\cal L}_{\rm BVP}}{\partial (\partial_n \phi)} \delta \phi \Big]_{\Sigma_n^{\rm i}}^{\Sigma_n^{\rm f}}
\label{eq:varSbnd2}\end{align}
The equivalence between the critical point, defined by $\delta {\cal S}_{\rm BVP}[X, \partial_a X, \phi, \partial_a \phi]=0$ and the equations of motion in the parentheses in (\ref{eq:varSeom}) requires the boundary terms in  (\ref{eq:varSbnd1}) and (\ref{eq:varSbnd2}) to vanish. In the BVP setting we consider here, the case $n=0$ makes reference to boundaries in the temporal direction $\Sigma_0=\tau$. By construction, we keep the values of the d.o.f. fixed at initial and final temporal parameter. I.e. the variations $\delta X^\mu$ and $\delta\phi$ have to vanish on the temporal boundaries at initial $\tau^{\rm i}=\Sigma_0^{\rm f}$ and final $\tau^{\rm f}=\Sigma_0^{\rm f}$. For the spatial boundaries related to the lowest and highest values of $\sigma_1,\ldots,\sigma_d$ there exist different possibilities to make the corresponding terms in (\ref{eq:varSbnd1}) and (\ref{eq:varSbnd2}) vanish. One possibility is to impose Dirichlet boundary conditions, which entails that the variations themselves vanish. On the other hand, another possibility is to impose von-Neumann or mixed boundary conditions, which involve specification of the derivatives of the degrees of freedom. We will discuss the treatment of boundary conditions in more detail in the explicit example of wave propagation in 1+1 dimensions in \cref{sec:proofofprinciple11d}.

%A central novelty of our action with dynamical coordinate mappings is the ability to discretize its parameters $\Sigma$ instead of introducing a discretization on the level of the space-time coordinates $X$. This fact, as we will show in detail in the next section, allows us to retain all space-time symmetries of the theory after discretization. The classical solution of the combined dynamical coordinate maps and fields will be obtained from a numerical optimization procedure based on the system action. 

In preparation for the next steps, recall that the location of the critical point of a functional does not change if we apply a monotonic function to its integrand. Also, for a reparameterization invariant action there exist an infinite number of ways to express the classical dynamics of the system it describes. In the world-line formalism, which describes the motion of a single point particle along its world-line trajectory, reparameterization invariance refers to the fact that the same world-line trajectory in $(d+1)$ dimensions can be traversed with different speeds relative to the world-line parameter $\gamma$. Such freedom leads to difficulties for a numerical minimizer in locating the classical field configuration. At the same time, the invariance of the physics under different parameterizations means that we are free to choose a parameterization at will.

These difficulties were avoided in the world-line formalism by squaring the integrand of the action. Motivated by this choice, we consider an action ${\cal E}$, where the integrand has been squared and the pre-factors have been modified for notational convenience
\begin{align}
\nonumber{\cal E}_{\rm BVP}&= \int d^{(d+1)}\Sigma\; E_{\rm BVP}[X,\partial_a X,\phi, \partial_a\phi]\\ 
&=\int d^{(d+1)}\Sigma \; \frac{1}{2} \Big\{ \Big(\frac{1}{T}V(\phi)-1\Big){\rm det}[g] + \frac{1}{T} \partial_a\phi(\Sigma)\partial_b\phi(\Sigma) {\rm adj}[g]_{ab}\Big\}. \label{eq:novelactionE}
\end{align}
Note that this action, as intended, is no longer invariant under general reparameterizations of $\Sigma$. However the critical point of ${\cal E}_{\rm BVP}$ and ${\cal S}_{\rm BVP}$ are equivalent, as the integrand in the action differs by only a monotonic transformation.

The equations of motion associated with the functional are obtained from its variation as
\begin{align}
\nonumber &\delta {\cal E}_{\rm BVP}[X, \partial_a X, \phi, \partial_a \phi]\\
&=\int d^{(d+1)}\Sigma \Big\{ 
\frac{\partial E_{\rm BVP}}{\partial X^\mu}\delta X^\mu + \frac{\partial E_{\rm BVP}}{\partial (\partial_a X^\mu)}\delta (\partial_a X^\mu) + \frac{\partial E_{\rm BVP}}{\partial \phi}\delta \phi + \frac{\partial E_{\rm BVP}}{\partial (\partial_a \phi)}\delta (\partial_a \phi) \Big\}\label{eq:varE}\\
&=\int d^{(d+1)}\Sigma \Big\{ 
\Big( \frac{\partial E_{\rm BVP}}{\partial X^\mu} -\partial_a \frac{\partial E_{\rm BVP}}{\partial (\partial_a X^\mu)}\Big)\delta X^\mu + \Big(  \frac{\partial E_{\rm BVP}}{\partial \phi} - \partial_a \frac{\partial E_{\rm BVP}}{\partial (\partial_a \phi)} \Big)\delta \phi \Big\} \label{eq:varEeom}\\
&\qquad + \sum_{n=0}^{d} \int \prod_{\tiny\begin{array}{c}a=0\\a\neq n\end{array}}^{d}d\Sigma_a \Big[ \frac{\partial E_{\rm BVP}}{\partial (\partial_n X^\mu)} \delta X^\mu \Big]_{\Sigma_n^{\rm i}}^{\Sigma_n^{\rm f}}\label{eq:varEbnd1}\\ 
&\qquad + \sum_{n=0}^{d} \int \prod_{\tiny\begin{array}{c}a=0\\a\neq n\end{array}}^{d}d\Sigma_a \Big[ \frac{\partial E_{\rm BVP}}{\partial (\partial_n \phi)} \delta \phi \Big]_{\Sigma_n^{\rm i}}^{\Sigma_n^{\rm f}}
\label{eq:varEbnd2}\end{align}
under the condition that the boundary terms (\ref{eq:varEbnd1}) and (\ref{eq:varEbnd2}) vanish.

For a manifestly Poincar\'e invariant action, such as \cref{eq:novelactionE}, Noether's theorem guarantees that each of its global continuous symmetries leads to one conserved current, which in turn defines a conserved charge. We may use the fact that $\delta_X {\cal E}\equiv{\cal E}[X+\delta X]-{\cal E}[X]=0$ and since the transformation only affects changes in the coordinate maps one ends up with the first two terms in the curly brackets in (\ref{eq:varE}). Carrying out integration by parts but retaining the total derivative we obtain
\begin{align}
0=\delta_X {\cal E}_{\rm BVP}&=\int d^{(d+1)}\Sigma \Big\{ 
\underbracket{\Big( \frac{\partial E_{\rm BVP}} {\partial X^\mu} -\partial_a \frac{\partial E_{\rm BVP}}{\partial (\partial_a X^\mu)}\Big)}_{\rm equations\,of\,motion\, for \, X}\delta X^\mu  +  \partial_a\Big[ \underbracket{\frac{\partial E_{\rm BVP}}{\partial (\partial_a X^\mu)} \delta X^\mu}_{J^a} \Big]\Big\}\label{eq:derivN1}\end{align}
Noether's theorem only holds if the equations of motion are fulfilled, i.e. if the term in the parenthesis in (\ref{eq:derivN1}) vanishes. We see that the conserved Noether current $J^a$ then is given by the expression
\begin{align}
    \partial_aJ^a= \partial_0J^0 + \sum_{n=1}^{d} \partial_n J^n =   \partial_a\Big[ \frac{\partial E_{\rm BVP}}{\partial (\partial_a X^\mu)} \delta X^\mu \Big] =0 \label{eq:NoetherCur}
\end{align}
The corresponding conserved charge is also defined as in regular field theory,
\begin{align}
    Q= \int \prod_{a>0}^d d\Sigma^a J^0 = \int d\sigma_1 \ldots d\sigma_d J^0, \label{eq:defQ}
\end{align}
as the spatial volume average over the temporal component of the Noether current $J^0$. Note that reference to temporal and spatial here is understood to refer to the parameters $\tau$ and $\vec{\sigma}$ respectively. In order for $Q$ to be constant in time, its temporal derivative must vanish. To evaluate $\partial_\tau Q$ we first exploit the definition \cref{eq:NoetherCur} and then apply the fundamental theorem of calculus
\begin{align}
    \partial_\tau Q&\overset{def.\,\cref{eq:defQ}}{=} \int   d\sigma_1 \ldots d\sigma_d \big(\partial_\tau J^0\big) \overset{\cref{eq:NoetherCur}}{=} - \int   d\sigma_1 \ldots d\sigma_d \big( \sum_{n=1}^{d} \partial_{\sigma_n}J^n\big) \\
    &=-\sum_{n=1}^{d} \int \prod_{\tiny\begin{array}{c}a=0\\a\neq n\end{array}}^{d}\left.d\sigma_a J^n \right|_{\sigma_n^{\rm i}}^{\sigma_n^{\rm f}}.
\end{align}
One has to inspect how the spatial components $J^n$ behave on the spatial boundaries. Only if they vanish there, does \cref{eq:defQ} imply $\partial Q/\partial t=0$ (c.f. (\ref{eq:varEbnd1}) and (\ref{eq:varEbnd2})). We will give an explicit expression for the relevant Noether charges when exploring 1+1 dimensional wave motion in \cref{sec:proofofprinciple11d}.

\subsection{Initial value formulation}
\label{sec:IBVPaction}

The action functional of \cref{eq:novelactionE} does not yet allow us to formulate a causal initial boundary value problem. This becomes apparent, when we consider the derivation of the equations of motion from it in \cref{eq:varE}. The boundary terms only vanish if the variation on the final $\tau$ slice are set to zero, i.e. the values of the d.o.f. must be known at final time.

The way to accommodate the natural absence of information about the final state of the system is to introduce a cleverly arranged doubling of degrees of freedom. This approach to causal initial value problems is known in the quantum physics community as the Schwinger-Keldysh construction (see e.g. \cite{sieberer_keldysh_2016}) and has been independently derived in the classical context by Galley in \cite{galley_classical_2013}. After deploying this construction to conventional point particle motion in \cite{Rothkopf:2022zfb} and to the world-line formalism in \cite{Rothkopf:2023ljz}, we will now exploit it to formulate a novel action with dynamical coordinate maps for initial boundary value problems involving fields.

\begin{figure}
\centering
    \includegraphics[scale=0.4]{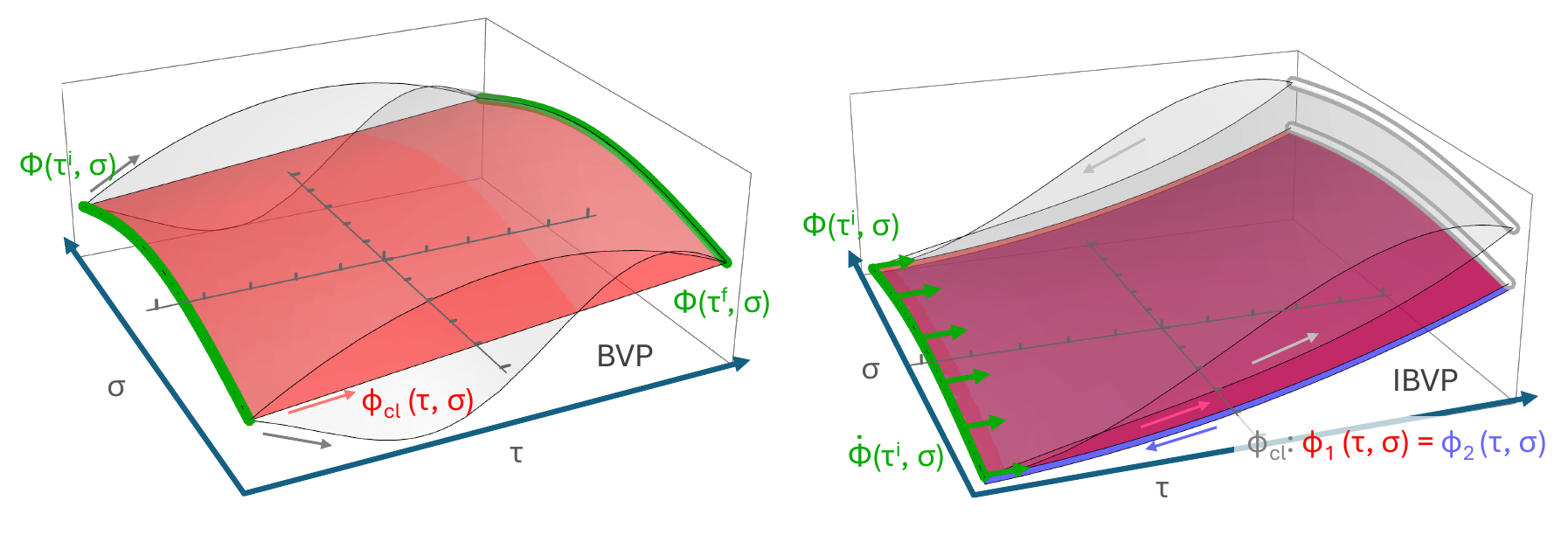}
    \caption{(left) The conventional BVP formulation of classical field theory, where both initial data $\phi(\tau^{\rm i},\vec{\sigma})$ and final data $\phi(\tau^{\rm f},\vec{\sigma})$ (in green) are provided a priori. While different field configurations exist that connect these two datasets (gray sheets), there exists a unique field configuration $\phi_{\rm cl}$ (in red) that constitutes the critical point of the action. It is this configuration that represents the configuration realized in nature. (right) The doubled d.o.f. construction of \cite{galley_classical_2013}, necessary to formulate causal IBVPs on the action level. Provided initial data on values $\phi(\tau^{\rm i},\vec{\sigma})$ and derivatives $\dot \phi(\tau^{\rm i},\vec{\sigma})$ (in green), one constructs a double shooting method. One copy of the field evolves forward (arrows to the right) and one copy evolves backward from the final state of the forward branch (arrows to the left). While different field configuration pairs exists that accommodate the initial data on the forward branch (e.g. gray sheets), only a single one (red and blue sheets) fulfills the requirement that it constitutes the critical point of the action with doubled degrees of freedom and that the forward and backward branch solutions agree (for more details see e.g. \cite{Rothkopf:2022zfb}).}
    \label{fig:sketchGalley}
\end{figure}

The doubling of degrees of freedom can be understood (see sketches and caption in \cref{fig:sketchGalley}) to furnish a double shooting method. One considers one copy of the degrees of freedom to live on the so-called forward temporal branch, a second copy on the backward temporal branch. The forward branch is initialized by the physical initial conditions while the backward branch receives information from the final state of the forward branch. As shown and discussed in detail in \cite{galley_classical_2013}, the solution of the IBVP corresponds to the trajectory, where the forward and backward path agree and the backward path reproduces the correct values and derivatives at initial time.

In the following we will thus refer to coordinate mappings and fields on the two branches as $X_1,X_2$ and $\phi_1,\phi_2$, each as functions of the parameters $\Sigma$ (which are not doubled).

Remember that it is the boundary terms in the final $\tau$ slice which prevent a causal formulation in \cref{eq:varE}. The approach of \cite{galley_classical_2013} resolves this issue by assigning the forward branch and the backward branch the same Lagrangian with opposite sign in the action
\begin{align}
{\cal E}_{\rm IBVP}&= \int d^{(d+1)}\Sigma\; E_{\rm IBVP}[X_1,X_2,\partial_a X_1,\partial_a X_2,\phi_1, \phi_2, \partial_a\phi_1,\partial_a\phi_2] ,\label{eq:novelactionEIVP}
%\nonumber &=\int d^{(d+1)}\Sigma\; \Big\{ E_{\rm BVP}[X_1,\partial_a X_1,\phi_1, \partial_a\phi_1] -E_{\rm BVP}[X_2,\partial_a X_2,\phi_2, \partial_a\phi_2] \Big\},
\end{align}
where 
\begin{align}
    \nonumber E_{\rm IBVP}=E_{\rm BVP}[X_1,\partial_a X_1,\phi_1, \partial_a\phi_1] -E_{\rm BVP}[X_2,\partial_a X_2,\phi_2, \partial_a\phi_2],
\end{align}
and by enforcing the following \textit{connecting conditions} on the values and derivatives of the degrees of freedom on the final time slice
\begin{align}
\nonumber &X_1^\mu(\tau=\tau^{\rm f},\vec{\sigma})=X_2^\mu(\tau=\tau^{\rm f},\vec{\sigma}), &&\phi_1(\tau=\tau^{\rm f},\vec{\sigma})=\phi_2(\tau=\tau^{\rm f},\vec{\sigma})\\
&\partial_0 X_1^\mu|_{\tau=\tau^{\rm f}}= \partial_0 X_2^\mu|_{\tau=\tau^{\rm f}}, &&  \partial_0 \phi_1|_{\tau=\tau^{\rm f}}= \partial_0 \phi_2|_{\tau=\tau^{\rm f}}.\label{eq:SKconnecting}
\end{align}

As derived in detail in \cref{sec:appconcond}, these conditions establish a causal variational principle, connecting the equations of motion, with the critical point of the functional ${\cal E}_{\rm IBVP}$ (given that  the spatial boundary terms vanish). The variational principle requires us to provide physical initial conditions about the values and derivatives of the degrees of freedom at initial time as well as connecting conditions between the forward and backward branches, which in turn avoid any non-causal temporal boundary terms.

In setting up the causal variational approach we introduced a second copy of the degrees of freedom. This proliferation of dynamical variables needs to be undone, for which, as argued in \cite{galley_classical_2013}, one must take the so-called physical limit. This limiting procedure consists of requiring that $X^\mu_1-X^\mu_2=X^\mu_-=0$ and $\phi_1-\phi_2=\phi_-=0$. In turn the dynamic variables $X^\mu_+$, $\phi_+$ converge towards the classical solution. In the equations of motion 
\begin{align}
\Big( \frac{\partial E_{\rm IBVP}}{\partial X_\pm^\mu} -\partial_a \frac{\partial E_{\rm IBVP}}{\partial (\partial_a X_\pm^\mu)}\Big)=0, \quad
\Big(  \frac{\partial E_{\rm IBVP}}{\partial \phi_\pm} - \partial_a \frac{\partial E_{\rm IBVP}}{\partial (\partial_a \phi_\pm)} \Big)=0,
\end{align}
the physical limit finds another interpretation. As $E_{\rm IBVP}= E_{\rm BVP}[X_1,\phi_1]- E_{\rm BVP}[X_2,\phi_2]$ arises from the difference of the BVP functionals it contains at least a linear dependence on the relative variables. In the physical limit, thus only the Euler-Lagrange equation where the derivative w.r.t. the relative coordinates is taken, has a chance to survive. Terms in $E_{\rm IBVP}$ that feature at most a linear dependence on $X^\mu_-$ and $\phi_-$ reduce to finite expressions containing $X^\mu_+$ and $\phi_+$. A more rigorous derivation of the physical limit can be made in the Schwinger-Keldysh formalism of quantum field theory, where the physical limit of \cite{galley_classical_2013} amounts to taking the classical limit of Feynman's path integral, as discussed in detail in \cite{Berges:2007ym}.

Let us take a brief look at the type of initial data we must provide to our initial boundary value problem. Since all dynamical quantities are dependent on the parameters $\Sigma^a$ and the kinetic terms in the action \cref{eq:novelaction} contain derivatives quadratic in $\tau=\Sigma_0$, both the values $X^\mu(\tau^{\rm i},\vec{\sigma}),\phi(\tau^{\rm i},\vec{\sigma})$ and the first derivative in $\tau$, i.e. $\partial_0 X^\mu(\tau^{\rm i},\vec{\sigma}),\partial_0\phi(\tau^{\rm i},\vec{\sigma})$ must be provided. From physical experiments we may obtain access to the values of the field at initial time and we are aware of the space-time coordinates at which this initial field has been measured. We can also determine the temporal derivative $\partial_{X^0}\phi$. On the other hand, the speed by which the coordinate maps and field evolve along $\Sigma_0$ are not fixed a priori. Thus there remains a freedom, as it is only the ratio $\frac{\partial\phi}{\partial X^0} = \frac{\partial \phi}{\partial \Sigma_0}/\frac{\partial X^0}{\partial \Sigma_0}$ that needs to be fixed. Since we consider time $X^0$ as a dynamic variable, its final value at $\tau^{\rm f}$ will depend on the speed assigned to it initially and it is an emergent property of the dynamics between the $X^\mu$'s and $\phi$.

\subsection{Lagrange multiplier formulation of initial, boundary and connecting conditions}
\label{sec:lagmultcont}

The construction of our causal variational approach to initial value problems, relied on the enforcement of different types of conditions on the dynamical degrees of freedom, be it in the form of initial data, spatial boundary data or the connecting condition of the forward and backward branches. So far we enforced these conditions implicitly, which is valid in an analytic setting. In preparation for a numerical treatment, we however have to incorporate these conditions explicitly in the system action. This is possible by use of Lagrange multipliers. I.e we amend \cref{eq:novelactionEIVP} by additional sets of multiplier functions $\lambda,\tilde \lambda$ for initial conditions and $\gamma,\tilde\gamma$ for the connecting conditions. The tilde in $\tilde \lambda$ and $\tilde \gamma$ indicates that these Lagrange multiplier functions are associated with constraints involving temporal derivatives, while the constraints enforced by $\lambda$ and $\gamma$ make reference only to the values. To enforce spatial boundary conditions we introduce the Lagrange multiplier functions $\kappa,\tilde \kappa, \xi,\tilde \xi$, where the tilde refers to those implementing the boundary conditions at $\vec{\sigma}^{\rm f}$ and those without at $\vec{\sigma}^{\rm i}$. Whenever we refer to quantities for which explicit Lagrange multipliers have been taken into account, a superscript letter $\rm L$ will be added. For notational simplicity we here refer only to Dirichlet boundary conditions in spatial direction, which can be generalized to Neumann or mixed boundary conditions

\begin{align}
&\nonumber{\cal E}^{\rm L}_{\rm IBVP}=\int d^{(d+1)}\Sigma\; \Big\{ E_{\rm BVP}[X_1,\partial_a X_1,\phi_1, \partial_a\phi_1] -E_{\rm BVP}[X_2,\partial_a X_2,\phi_2, \partial_a\phi_2]\Big\} \\
\nonumber+&\int \prod_{a=1}^{d}d\Sigma_a\Big\{ \lambda_\mu\big( X_1^\mu(\tau^{\rm i},\vec{\sigma}) - X^\mu_{\rm IC}\big) +\lambda_\phi \big( \phi_1(\tau^{\rm i},\vec{\sigma}) - \phi_{\rm IC}\big)\\
\nonumber+&\tilde \lambda_\mu\big( \partial_0 X_1^\mu(\tau^{\rm i},\vec{\sigma}) - \dot X^\mu_{\rm IC}\big)  +\tilde \lambda_\phi \big( \partial_0\phi_1(\tau^{\rm i},\vec{\sigma}) - \dot \phi_{\rm IC}\big)\\
\nonumber+&\gamma_\mu\big( X_1^\mu(\tau^{\rm f},\vec{\sigma}) - X_2^\mu(\tau^{\rm f},\vec{\sigma})\big) +\gamma_\phi \big( \phi_1(\tau^{\rm f},\vec{\sigma}) - \phi_2(\tau^{\rm f},\vec{\sigma}) \big)\\
+&\tilde \gamma_\mu\big( \partial_0 X_2^\mu(\tau^{\rm f},\vec{\sigma}) - \partial_0 X_2^\mu(\tau^{\rm f},\vec{\sigma})\big)  +\tilde \gamma_\phi \big( \partial_0\phi_1(\tau^{\rm f},\vec{\sigma}) - \partial_0\phi_2(\tau^{\rm f},\vec{\sigma})\big)\Big\}\\
\nonumber+&\sum_{j=1}^{d} \int \prod_{\tiny\begin{array}{c}a=0\\a\neq j\end{array}}^{d} d\Sigma_a\Big\{\kappa^j_{\mu} \big( X_1^\mu(\sigma_j^{\rm i}) - X^\mu_{\rm sBCL}(\sigma_j^{\rm i})\big)
\nonumber+\xi^j_{\mu} \big( X_2^\mu(\sigma_j^{\rm i}) - X^\mu_{\rm sBCL}(\sigma_j^{\rm i})\big)\\
\nonumber+&\tilde \kappa^j_{\mu} \big( X_1^\mu(\sigma_j^{\rm f}) - X^\mu_{\rm sBCR}(\sigma_j^{\rm f})\big)
\nonumber+\tilde \xi^j_{\mu} \big( X_2^\mu(\sigma_j^{\rm f}) - X^\mu_{\rm sBCR}(\sigma_j^{\rm f})\big)\\
\nonumber+&\kappa^j_{\phi} \big( \phi_1(\sigma_j^{\rm i}) - \phi_{\rm sBCR}(\sigma_j^{\rm i})\big)
\nonumber+\xi^j_{\phi} \big( \phi_2(\sigma_j^{\rm i}) - \phi_{\rm sBCR}(\sigma_j^{\rm i})\big)\\
+&\tilde \kappa^j_{\phi} \big( \phi_1(\sigma_j^{\rm f}) - \phi_{\rm sBCL}(\sigma_j^{\rm f})\big)
+\tilde \xi^j_{\phi} \big( \phi_2(\sigma_j^{\rm f}) - \phi_{\rm sBCL}(\sigma_j^{\rm f})\big)\Big\},\label{eq:novelactionEIVPl1}\\
=&\int d^{(d+1)}\Sigma\; \Big\{ E^{\rm L}_{\rm BVP}[X_1,\partial_a X_1,\phi_1, \partial_a\phi_1] - E^{\rm L}_{\rm BVP}[X_2,\partial_a X_2,\phi_2, \partial_a\phi_2]\Big\}.\label{eq:novelactionEIVPlabs}
\end{align}
Above in the second and third line, we enforce the initial conditions on the forward branch of $X^\mu$ and the field $\phi$. I.e. we specify their value and derivative at initial temporal parameter $\tau^{\rm i}$. In the fourth and fifth line we connect the two branches at final temporal parameter $\tau^{\rm f}$. The different quantities $\lambda,\tilde\lambda,\gamma,\tilde \gamma$ are functions of all the $d$ spatial parameters $\vec{\sigma}$ but do not depend on the temporal parameter $\tau$. 

Enforcement of the spatial boundary conditions constitutes the rest of the terms, where separate Lagrange multiplier functions are applied to fix left and right spatial boundaries, corresponding to spatial parameters $\vec{\sigma}^{\rm i}$ and $\vec{\sigma}^{\rm f}$ respectively. The spatial boundary conditions are enforced on each of the branches separately with different Lagrange multiplier functions. Here the quantities $\kappa^j,\tilde \kappa^j,\xi^j,\tilde\xi^j$ are functions that depend on the temporal parameter $\tau$ and $(d-1)$ spatial parameters, i.e. they depend on all components of $\vec{\sigma}$, except for the $j$-th parameter $\sigma_j$, indicate by the superscript.

As discussed in \cite{Rothkopf:2023vki} for the world-line formalism, the Lagrange multiplier contributions can be absorbed into the action integral by amending each of the boundary terms with delta-functions, placing the multiplier functions at the correct initial or final parameter values. This is indicated in \cref{eq:novelactionEIVPlabs} by the redefinition of  $E_{\rm BVP}\to E^{\rm L}_{\rm BVP}$. The presence of the Lagrange multiplier terms affects the expressions one obtains for the equations of motion, and for the Noether charge in (\ref{eq:varEeom}) and (\ref{eq:NoetherCur}) we have to take derivatives of $E^{\rm L}_{\rm BVP}$ instead of the boundary value Lagrangian $E_{\rm BVP}$. This leads to a modified Noether current
\begin{align}
    \partial_a(J^{\rm L})^a=  \partial_a\Big[ \frac{\partial  E^{\rm L}_{\rm BVP}}{\partial (\partial_a X^\mu)} \delta X^\mu \Big] =0,\label{eq:NoetherCurmodL}
\end{align}
where the $a$-th component is denoted $(J^{\rm L})^a$.
In turn we find the following modification of the Noether charge density
\begin{align}
    (J^{\rm L})^0=  \Big(\frac{\partial E_{\rm BVP}}{\partial (\partial_0 X^\mu)} + \tilde\lambda_\mu\delta(\tau-\tau^{\rm i}) +\tilde \gamma_\mu\delta(\tau-\tau^{\rm f}) \Big) \delta X^\mu. \label{eq:NoetherChrgDenmodL}
\end{align}
Note that only those Lagrange multipliers $\tilde\lambda_\mu,\tilde\gamma_\mu$ appear, which in the action are associated with constraints involving $\tau$ derivatives of the coordinate maps. We will provide a concrete example of the corrections to the conserved charges when discussing wave propagation in $1+1$ dimensions in \cref{sec:1p1dcontform}.

At this point we have furnished the full continuum variational formalism in $(d+1)$ dimensions. Elevating coordinate maps to dynamical degrees of freedom led us to the action ${\cal E}_{\rm BVP}$ in \cref{eq:novelactionE}, on which we base our genuine IBVP formulation ${\cal E}_{\rm IVP}$ in \cref{eq:novelactionEIVP} and whose underlying requirements in the form of initial, connecting, and boundary conditions have been made explicit in ${\cal E}^{\rm L}_{\rm IBVP}$ in \cref{eq:novelactionEIVPlabs}. It now remains to discretize this variational approach to IBVPs, which is the focus of the next section.

\section{Discretization of the IBVP}
\label{sec:discrIBVP}

What distinguishes our novel variational approach to IBVPs from the conventional formulation is the presence of dynamical coordinate maps. As we will show below, these maps allow us to discretize the system on the level of the underlying parameters $\Sigma$ instead of the space-time coordinates $X$. In turn we obtain a discretized action  which \textit{retains all space-time symmetries} and thus \textit{preserves the continuum Noether charges}. The fact that the coordinate maps adapt to the field dynamics in a way that preserves the continuum symmetries provides a mechanism, where the coarseness of the dynamically emerging space-time grid is automatically adjusted, implementing in effect an \textit{automatic adaptive mesh refinement} (AMR) procedure.

Retention of symmetries and the preservation of the associated Noether charges in the discrete setting is the central goal of our discretization strategy. The derivation of Noether currents and charges requires the application of integration-by-parts (IBP), hence, our discretization will be based on summation-by-parts (SBP) finite difference operators \cite{svard2014review,fernandez2014review,lundquist2014sbp}, which mimic IBP exactly in the discrete setting. %After discretizing the action, the classical field configuration will be obtained by numerical optimization.

Our notation in the discrete setting is as follows: discrete vectors will be referred to as bold letters, such as $\bm X$, $\bm \phi$ or $\bm f$. Matrix operators that act on these vectors are denoted by capital double-struck letters, e.g. $\mathds{H}$ or $\mathds{D}$. In case we wish to refer to operations affecting functions only along a single parameter (either $\tau$ or a $\vec{\sigma}$), we will use instead lowercase double-struck letters, such as $\mathbb{h}_a$ or $\mathbb{d}_a$, where subscripts denote which of the $(d+1)$ space-time dimension the matrix operator is associated with. The identity matrix in $(d+1)$ dimensions is thus referred to as $\mathds{I}$ whereas its one-dimensional counterpart is $\mathbb{i}$.

%Below we derive the discrete formalism in the general (d+1) dimensional setting with the goal of providing the reader with the details necessary for implementation in practice. An illustrative example of the discrete formalism in (1+1) dimensions can be found in subsection \ref{sec:discr2d}, within our proof of principle \cref{sec:proofofprinciple11d}.

The underlying $(d+1)$ dimensional parameter space described by $\Sigma^a=(\tau,\vec{\sigma})^a$ is discretized on a hypercubic grid with $N_a$ points and grid spacing $\Delta_a$, not necessarily equal for each dimension. Once we select the intervals $\Sigma_a\in [\Sigma_a^{\rm i},\Sigma_a^{\rm f}]$ in which the parameters lie, the grid spacings follow as $\Delta_a=(\Sigma_a^{\rm f}-\Sigma_a^{\rm i})/(N_a-1)$. Any function  $f(\Sigma)$, which depends on the $\Sigma_a$ parameters, is discretized as an array $ \bm f$ with \texttt{TotVol}$= \prod_{a=0}^{d} N_a$ entries. We refer to the discretized dynamical degrees of freedom, the coordinate maps and field on the forward and backward branch, as ${\bm X}^\mu_{1,2}$ and ${\bm \phi}_{1,2}$ respectively. The order of the individual entries in these arrays determines the form of tensor products (implemented as Kronecker products) deployed below. We choose to let $\Sigma_0=\tau$ run slowest and $\Sigma_d=\sigma_d$ fastest. 
We need to specify which entry of a discrete array ${\bm f}$ refers to a specific position on the abstract parameter grid. Take e.g. the coordinates $(\tau^{\rm i}+\Delta\tau n_0,\sigma_1^{\rm i}+\Delta\sigma_1 n_1,\ldots,\sigma_d^{\rm i}+\Delta\sigma_d n_d)$, which correspond to taking $n_0$ steps with $\Delta\tau$ and $n_i$ steps with $\Delta\sigma_i$.  Correspondingly we introduce an index function \texttt{index}$(\tau^{\rm i}+\Delta\tau n_0,\sigma_1^{\rm i}+\Delta\sigma_1 n_1,\ldots,\sigma_d^{\rm i}+\Delta\sigma_d n_d)= (\ldots(n_0 N_1+n_1)N_2 \ldots + n_d)$ that computes the integer position along the discrete function array corresponding to that specific position. The counter variables $n_0\in[0,N_\tau-1]$ and $n_a\in[0,N_a-1]$ take on only integer values. In $(1+1)$ dimensions, the index is computed explicitly as \texttt{index}$(n_\tau,n_\sigma)=(n_\tau N_\sigma +n_\sigma)$.

Integration-by-parts rests on the interplay of differentiation and integration. We start by approximating the $(d+1)$ dimensional integral in the action functional. Note that discretization here refers to the abstract parameters. Let us introduce an inner product on discretized functions ${\bm f}$ and ${\bm g}$, which is characterized by a positive definite diagonal matrix $\mathds{H}\in\mathbb{R}^{\texttt{TotVol}\times\texttt{TotVol}}$, so that $\int d^{(d+1)}\Sigma \; f(\Sigma)g(\Sigma) \approx ( {\bm f}, {\bm g}) = {\bm f}^T \mathds{H} {\bm g}$. 
Since in our action, the boundaries of the integrals in each direction $\Sigma_a$ are independent from the values of the other integration variables and we use quadrature rules, which are local to each dimension, we can consider the multidimensional integral described by $\mathds{H}$ as consecutive and independent applications of one-dimensional integration. Hence the matrix $\mathds{H}$ can be written as a tensor product of one-dimensional quadrature rules $\mathbb{h}_a\in\mathbb{R}^{N_a\times N_a}$
\begin{align}
    \mathds{H}=\mathbb{h}_0\otimes\mathbb{h}_1\otimes\ldots \otimes\mathbb{h}_d,
\end{align}
where $\mathbb{h}_0$ implements integration along $\tau$ and $\mathbb{h}_i$ along the $\sigma_i$ direction.

Similarly, we define a finite difference operator, which implements the SBP property in one dimension and use the tensor product to  construct the corresponding operator in multiple dimensions. To obtain a finite difference operator $\mathbb{d}_a\in\mathbb{R}^{N_a\times N_a}$ of order $r$ that approximates the continuum derivative $d/d\Sigma_a$ and which mimics IBP exactly, it must fulfill the following SBP properties
\begin{align}
    &  \mathbb{d}_a{\bm \Sigma}_a^n = n {\bm \Sigma}_a^{n-1} \quad {\rm for\,n \leq r, \; no\, summation\,in\,a\, implied,\;} {\bm \Sigma}_a{\rm \in\mathbb{R}^{N_a}}\label{eq:derivapprx}\\
    & \mathbb{d}_a=\mathbb{h}_a^{-1}\mathbb{q}_a, \quad \mathbb{q}_a^{\rm T}+\mathbb{q}_a=\mathbb{e}_a^N-\mathbb{e}_a^0={\rm diag}[-1,0,\ldots,0,1].\label{eq:SBP}
\end{align}
Here ${\bf \Sigma}_a$ denotes a vector with $N_a$ entries, which linearly increase in magnitude, e.g., ${\bf \Sigma}_a=(0,\Delta_a,2\Delta_a,\ldots)$. The matrix $\mathbb{q}_a$ encodes the so-called SBP property and \cref{eq:derivapprx} establishes that $\mathbb{d}_a$ acts as a derivative, using polynomials as the underlying function space (for a definition of SBP operators on general function spaces see \cite{glaubitz2023summation,glaubitz2024summation,glaubitz2024energy}). The matrices $\mathbb{e}^0_a$ and $\mathbb{e}^N_a$ project out the values at the appropriate boundaries of the $\Sigma_a$ direction. To be able to consistently act on our arrays ${\bm X}^\mu_{1,2}$ and ${\bm \phi}_{1,2}$ we must define $\mathds{D}_a\in\mathbb{R}^{\texttt{TotVol}\times\texttt{TotVol}}$ via the tensor product with unit matrices $\mathbb{i}$ as
\begin{align}
    \mathds{D}_a= \overbracket{\underbracket{\mathbb{i} \otimes \ldots \otimes \mathbb{i}}_{\rm a\, times} \otimes \mathbb{d}_a \otimes \mathbb{i} \ldots \otimes \mathbb{i}}^{\rm d+1\, times}
\end{align}

As concrete realizations of the $\mathbb{h}_a$ and $\mathbb{d}_a$ we will use the two lowest order diagonal norm SBP schemes referred to as \texttt{SBP121} and \texttt{SBP242}. While the former exhibits second order accuracy in the interior, the latter is of fourth order. The accuracy on the boundary is reduced by one order for \texttt{SBP121} and by two orders for \texttt{SBP242}\footnote{For the \texttt{SBP121} operator this is evident as the second order symmetric stencil in the interior in \cref{eq:SBP121} reduces to a first order forward or backward stencil on the boundary.}. The quadrature rule for the \texttt{SBP121} scheme is the well-known trapezoidal rule  
\begin{equation}
 \mathbb{h}_a^{[1,2,1]}=\Delta \Sigma_a \left[ \begin{array}{ccccc} 1/2 & & & & \\ &1 & & &\\ & &\ddots && \\ &&&1&\\ &&&&1/2 \end{array} \right],
\quad 
\mathbb{d}_a^{[1,2,1]}=
\frac{1}{2 \Delta \Sigma_a}
\left[ \begin{array}{ccccc} -2 &2 & & &\\ -1& 0& 1& &\\ & &\ddots && \\ &&-1&0&1\\ &&&-2&2 \end{array} \right].\label{eq:SBP121}
\end{equation}
The uniquely defined \texttt{SBP121} operator features the lowest order central symmetric stencil in the interior and the forward and backward stencils at the left and right boundary of the parameter interval respectively. %Thus the SBP property leads to a collection of $2\times2$ entries associated with the boundary behavior in the top left and lower right corners of the above matrices, different from the interior stencil region. 
To achieve fourth order accuracy in the interior and second order on the boundary, the corresponding symmetric stencil is modified by a collection of $4\times 4$ entries in the corners
\begin{align}
\nonumber &\mathbb{h}_a^{[2,4,2]}=\Delta \Sigma_a \left[ \begin{array}{ccccccc} 
\frac{17}{48} & & & & & & \\
 & \frac{59}{48} & & & & & \\
 & &\frac{43}{48} & & & & \\
 & & &\frac{49}{48} & & & \\
 & & & & &1 & \\
 & & & & & &\ddots \\
 \end{array} \right],\end{align}

 \begin{align}
&\mathbb{d}_a^{[2,4,2]}=
\frac{1}{\Delta \Sigma_a}
\left[ \begin{array}{ccccccccc} 
-\frac{24}{17}&\frac{59}{34} & -\frac{4}{17} & -\frac{3}{34} & && & & \\
-\frac{1}{2}& 0 & \frac{1}{2} & 0 & && & & \\
 \frac{4}{43}& -\frac{59}{86} & 0 & \frac{59}{86}&-\frac{4}{43} && & & \\
\frac{3}{98}& 0& -\frac{59}{86}  & 0&\frac{32}{49}&-\frac{4}{49}& & & \\
& &\frac{1}{12}  & -\frac{2}{3}&0&\frac{2}{3}& -\frac{1}{12}& & \\
&&&&&&&\ddots&
 \end{array} \right].\label{eq:SBP242}
\end{align}

%The matrices $\mathbb{d}_a$ by construction are not anti-symmetric and their 
Note that the left and right eigenvectors of the matrices $\mathbb{d}_a$ are not necessarily the same. In fact it has been shown in \cite{Rothkopf:2022zfb} that while the right eigenvectors contain the physical zero mode, i.e. the constant function, the zero modes among the left eigenvectors are so-called $\pi$-modes (a manifestation of the doubler problem of lattice field theory \cite{Wilson:1974sk}). When SBP operators are used in the setting of discretized partial differential equations, only their right eigenvectors play a role. Since in that case only the constant function is annihilated by the derivative operator, one calls $\mathds{D}_a$ null-space consistent (see ref.~\cite{svard2019convergence}). In the discretization of an action functional, which is quadratic in the derivatives, additional care must be taken. When locating the critical point of the action, the $\pi$-mode will, if not tamed, contaminate the solution \cite{Rothkopf:2022zfb}. In this sense the matrices $\mathbb{d}_a$ of \cref{eq:SBP121,eq:SBP242} are not yet null-space consistent.

In prior publications \cite{Rothkopf:2022zfb,Rothkopf:2023ljz} it has been established that in a one-dimensional setting, the $\pi$-mode can be avoided by including boundary information into $\mathds{D}_a$, taking inspiration from the technique of simultaneous-approximation terms (SAT) \cite{carpenter1994time}. The strategy relies on incorporating penalty terms related to boundary information into a novel derivative $\bar{\mathds{D}}_a$, without unphysical zero modes\footnote{Since the functions we apply the derivative operator on are purely real, the Wilson-term regularization \cite{Wilson:1974sk} often deployed in lattice field theory is not applicable here.}
\begin{align}
    &\bar{\mathds{D}}^f_a {\bm f} = \mathds{D}_a {\bm f} + \sigma_0\mathds{H}^{-1}\mathds{E}^0_a({\bm f} - {\bm f}_{\rm bnd})= \mathds{D}_a {\bm f} + \mathds{S}_a({\bm f} - {\bm f}_{\rm bnd}),\label{eq:regSBP}
\end{align}
where the projector matrices $\mathds{E}^0_a$ and $\mathds{S}_a$ are defined as 
\begin{align}
    \mathds{E}^0_a= \overbracket{\underbracket{\mathbb{i}\otimes \ldots \otimes \mathbb{i}}_{\rm a\, times} \otimes \mathbb{e}^0_a \otimes \mathbb{i} \ldots \otimes \mathbb{i}}^{\rm d+1\, times},\quad  \mathds{S}_a=\overbracket{\underbracket{\mathbb{i}\otimes \ldots \otimes \mathbb{i}}_{\rm a\, times} \otimes \sigma_0\mathbb{h}_a^{-1}\mathbb{e}^0_a \otimes \mathbb{i} \ldots \otimes \mathbb{i}}^{\rm d+1\, times}
\end{align}
We choose $\sigma_0=1$ to obtain a regularized multidimensional $\bar{\mathds{D}}^f_a$ with eigenvalues featuring a positive definite real part, in agreement with the one-dimensional setting in \cite{Rothkopf:2022zfb,Rothkopf:2023ljz}.

As we only consider Dirichlet boundary conditions in spatial directions, the array ${\bm f}_{\rm bnd}$ contains values of the function ${\bm f}$ only on the boundary itself and may be set to zero otherwise. Note that a different regularized SBP operator ensues for each degree of freedom ($\bar{\mathds{D}}^\mu_a$ to act on ${\bm X}^\mu_{1,2}$ and $\bar{\mathds{D}}^\phi_a$ to act on ${\bm \phi}_{1,2}$), as each may be assigned different boundary behavior. The above definition entails that the derivative operator in the temporal direction $\tau$ will include the initial conditions of the corresponding degree of freedom, while for the spatial derivative operator we will use the spatial boundary conditions at $\sigma^{\rm i}$. 

The shift operation involving ${\bm f}_{\rm bnd}$ in \cref{eq:regSBP} can be efficiently implemented by introducing affine coordinates. The introduction of affine coordinates entails that all discrete function arrays are amended by one more entry of value one and each operator matrix is endowed with one more row and column with the value one placed in the lower right corner. The new column on the right contains the shift. The new row, except for the corner is filled with the values zero. The resulting structure of the regularized SBP operator is sketched below with the conventional structures shaded in purple and the affine coordinate additions in orange
\begin{equation}
\bar{\mathds{D}}_a^f \bar{{\bf f}} =\begin{tikzpicture}[baseline={-0.5ex},mymatrixenv]
        \matrix [mymatrix,inner sep=4pt] (m)  
        {
      &    &   &\\
      & \mathds{D}_a+\mathds{S}_a & & -\mathds{S}_a{\bf f}_{\rm bnd}   \\
      &    &   &\\
      & 0   &   & 1\\
    };
    \begin{scope}[on background layer,rounded corners]
     \node [fit=(m-1-1) (m-3-3),purpleish,inner xsep=1.5pt,inner ysep=2.5pt]{};
     \node [fit=(m-1-4) (m-3-4),orangeish,inner xsep=1.5pt,inner ysep=2.5pt]{};
    \node [fit=(m-4-1) (m-4-3),orangeish,inner xsep=1.5pt,inner ysep=2.5pt]{};
    \node [fit=(m-4-4) (m-4-4),orangeish,inner xsep=1.5pt,inner ysep=2.5pt]{};
    \end{scope}
\end{tikzpicture}
\begin{tikzpicture}[baseline={-0.5ex},mymatrixenv]
        \matrix [mymatrix,inner sep=4pt] (m)  
        {
       \\
      {\bf f}\\
      \\
     1 \\
    };
    \begin{scope}[on background layer,rounded corners]
     \node [fit=(m-1-1) (m-3-1),purpleish,inner xsep=2.5pt,inner ysep=2.5pt]{};
    \node [fit=(m-4-1) (m-4-1),orangeish,inner xsep=1.5pt,inner ysep=2.5pt]{};
    \end{scope}
\end{tikzpicture}.\label{eq:sketchregSBP}
\end{equation}
For a more comprehensive discussion of the implementation of the affine coordinate regularization we refer the reader to \cite{Rothkopf:2022zfb,Rothkopf:2023ljz}.

With appropriately regularized multidimensional SBP operators in place, we proceed to discretize first the novel action ${\cal E}_{\rm BVP}$, as its Lagrangian $E_{\rm BVP}$ constitutes the building block of the genuine causal IBVP formulation. In addition it allows us to highlight the central property of our approach, the retention of continuum space-time symmetries after discretization. Replacing derivatives by regularized SBP operators we obtain
\begin{align}
     {\bm g}_{ab}=G_{\mu\nu} (\bar{\mathds{D}}^\mu_a {\bm X}^\mu)\circ (\bar{\mathds{D}}^\nu_b {\bm X}^\nu), \quad {\rm det}[{\bm g}]=\sum_{i_0,\ldots,i_{d}}\epsilon_{i_0\cdots i_{d}} {\bm g}_{0,i_0}\circ\cdots\circ{\bm g}_{d,i_{d}}\label{eq:discrg}
\end{align}
where the symbol $\circ$ refers to pointwise multiplication of discrete arrays and $\epsilon_{i_0\cdots i_{d}}$ denotes the $(d+1)$ dimensional totally anti-symmetric tensor. For ease of notation when dealing with integration over products of more than two functions, we introduce the vector quantity ${\bf h}=\mathds{H}{\bf 1}$, where ${\bf 1}$ denotes the vector with all entries being one. ${\bf h}$ contains the entries of the diagonal matrix $\mathds{H}$ allowing us to write the integral of ${\bf f}$ as $({\bf f},{\bf 1})={\bf f}^T{\bf h}$. Applying the discretization described above to the continuum action ${\cal E}_{\rm BVP}$ in \cref{eq:novelactionE} results in the following expression
\begin{align}
    \nonumber \mathds{E}_{\rm BVP}&[{\bm X}_1^\mu,\bar{\mathds{D}}^\mu_a {\bm X}_1^\mu,{\bm \phi}_1, \bar{\mathds{D}}^\phi_a{\bm \phi}_1]=\\
    &\frac{1}{2}\Big\{ \Big(\frac{1}{T}V({\bm \phi}_1)-1\Big)\circ{\rm det}[{\bm g}_1] +\frac{1}{T} (\bar{\mathds{D}}^\phi_a{\bm \phi}_1) \circ  (\bar{\mathds{D}}^\phi_b{\bm \phi}_1)\circ{\rm adj}[{\bm g}_1]_{ab}\Big\}^T {\bm h}.\label{eq:discrEBVP}
\end{align}

Let us start by investigating the discretized induced metric ${\bm g}_{ab}$ in (\ref{eq:discrg}). We see that the discrete derivatives in (\ref{eq:discrg}) are combined with the space-time metric $G$ in the same way as in the continuum. By Einstein summation convention we sum over both $\mu$ and $\nu$ indices. This entails that under a Lorentz transform, a transformation matrix $\Lambda$ is applied to each of the coordinate maps ${\bf X}$. The discretized expression remains invariant by definition of $\Lambda$, i.e $G_{\mu\nu}\Lambda^\mu_\kappa\Lambda^\nu_\gamma=G_{\kappa\gamma}$. Thus each entry of the matrix ${\bm g}$ individually stays invariant under Lorentz transformations in the discrete setting. Furthermore the discrete derivatives $\bar{\mathds{D}}_a^\mu$ annihilate the constant function exactly, such that adding a constant global shift to each entry of the arrays ${\bf X}^\mu$ will not affect the discrete ${\bm g}$, making it invariant under the full Poincar\'e group of transformations. In turn, as the only reference to space-time coordinates in \cref{eq:discrEBVP} is made via ${\bm g}$ we can conclude that the discretized action $\mathds{E}_{\rm BVP}$ retains the full manifest Poincar\'e invariance of the continuum action. 

This situation is crucially different from the conventional discretization in the space-time coordinates themselves. In our novel formulation it is the underlying abstract parameters $\Sigma$, which are discretized, leaving the entries of the ${\bf X}^\mu$ arrays to take on arbitrary values. In particular the values of the entries of ${\bf X}^\mu$ can be changed by arbitrarily infinitesimal amounts, the prerequisite for an application of Noether's theorem. 

In the derivation of the continuum Noether current only integration by parts was deployed and our discretization via SBP operators mimics IBP exactly in the discrete setting. We may thus directly infer that the expression for the Noether current ${\bm J}^a$ and its associated Noether charge density ${\bm q}\equiv{\bm J}^0$ in the discrete setting reads
\begin{align}
    \mathds{D}_a{\bf J}^a=  \mathds{D}_a\Big[ \frac{\partial \mathds{E}_{\rm BVP}}{\partial (\mathds{D}_a {\bf X}^\mu)} \delta {\bf X}^\mu \Big] =0, \quad {\bm q} = \frac{\partial \mathds{E}_{\rm BVP}}{\partial (\mathds{D}_0 {\bf X}^\mu)}\delta {\bf X}^\mu.\label{eq:NoetherCurDiscr}
\end{align}
Note that the conserved Noether charge ${\bm Q}$ is obtained from integration of the associated charge density ${\bm q}$ over all spatial direction (c.f. \cref{eq:defQ}).

The differentiation with respect to the array containing the derivative of the space-time coordinate maps $\partial/\partial (\mathds{D}_a {\bf X}^\mu)$ is to be understood in a point-wise fashion. In the discrete setting it can be explicitly implemented using the discrete Dirac delta function $ \mathfrak{d}_k$ (also known as lifting operator \cite{nordstrom2017roadmap}). When ${\bf f}$ is integrated over $ \mathfrak{d}[k]$ the result must correspond to the $k$-th entry of ${\bf f}$, i.e. ${\bf f}^{\rm T}{\mathds{H}}\mathfrak{d}[k] ={\bf f}[k]$. We can consider both the discrete counterpart to the delta function collapsing integration in all dimensions $\delta^{(d+1)}(\Sigma-\Sigma_{\rm ref})$ and those acting only in a single dimension $\delta(\Sigma^a-\Sigma^a_{\rm ref})$
\begin{align}
    \mathfrak{d}[k] = {\mathds{H}}^{-1} {\bf e}_k, \quad \mathfrak{d}^a[j] =  \mathbb{h}_a^{-1} {\bf e}^a_j \label{eq:defdeltas} 
\end{align}
where ${\bf e}\in\mathbb{R}^{\texttt{TotVol}}$ and ${\bf e}^a\in\mathbb{R}^{N_a}$ refer to appropriate unit vectors filled with zero except for the value one at the position $k$ or $j$ at which the delta function is non-zero. 

Using the discrete delta function defined above, we can mimic continuum functional derivatives $\frac{\delta }{\delta f(z)}\int dx f(x) g(x) = \int dx \delta(x-z)g(x) = g(z)$ in the discrete setting as  $\partial ({\bf f},{\bf g})/\partial {\bf f}[k]= (\mathfrak{d}[k],{\bf g})={\bf g}[k]$.

In preparation for expressing the Noether charge in the discrete setting, we also need to define a matrix operator, which implements spatial integration. I.e. it must transform a discrete array ${\bf f}$ with \texttt{TotVol} entries to an array that only contains $N_\tau$ entries, an operation represented via a $N_\tau\times$\texttt{TotVol} matrix we call $\mathds{H}_\sigma$. For convenient implementation in practice, we define the $\prod_{a>0}^{d} N_a$ component vector $({\bm h}_\sigma)_i=(\mathbb{h}_1\otimes\ldots\otimes\mathbb{h}_d)_{ii}$, which contains the diagonal entries of the spatial quadrature matrix $\mathbb{h}_1\otimes\ldots\otimes\mathbb{h}_d$. Since we decided that our coordinates are ordered in the discrete setting such that spatial coordinates run fastest, the matrix $\mathds{H}_\sigma$ can be written as
\begin{align}
    \mathds{H}_\sigma = \left(\begin{array}{cccc} {\bm h}^T_\sigma & 0 & \cdots &0\\
                                           0 & {\bm h}^T_\sigma& \cdots &0\\
                                           \vdots&&&\\
                                           0 & 0 &\cdots& {\bm h}^T_\sigma\end{array}\right).\label{eq:spatint}
\end{align}

In \cref{sec:lagmultcont} we showed that in order to explicitly include the necessary conditions related to initial and boundary data, as well as the connection of the forward and backward branch, we can use Lagrange multiplier functions absorbed into a redefinition of the BVP functional $E_{\rm BVP}\to E^{\rm L}_{\rm BVP}$. In the discrete setting we similarly incorporate these terms in a new $\mathds{E}_{\rm BVP}\to {\mathds{E}}^{\rm L}_{\rm BVP}$ using the discrete delta function defined above. The Lagrange multiplier terms affect the resulting equations of motion and the Noether charge. As discussed in detail for the one-dimensional case in \cite{Rothkopf:2023vki}, the modification to the Noether charge can be uniquely related to the Lagrange multiplier terms associated with the derivatives of the spatial coordinates. Denoting the discretized Lagrange multiplier functions by bold-face Greek letters, the discretized version of the continuum Noether charge density of \cref{eq:NoetherChrgDenmodL} reads
\begin{align}
    {\bm q}^{\rm L} = \frac{\partial {\mathds{E}}^{\rm L}_{\rm BVP}}{\partial (\mathds{D}_0 {\bf X}^\mu)}\delta {\bf X}^\mu = \Big( \frac{\partial {\mathds{E}}_{\rm BVP}}{\partial (\mathds{D}_0 {\bf X}^\mu)} + \tilde{\bm \lambda}_\mu \circ \mathfrak{d}^0[0] + \tilde{\bm \gamma}_\mu\circ \mathfrak{d}^0[N_0] \Big)\delta {\bf X}^\mu,\label{eq:NoetherCurDiscrLagr}
\end{align}
while the following explicit expression for the corresponding Noether charge ensues
\begin{align}
    {\bm Q}^{\rm L} =  \Big( \mathds{H}_\sigma \frac{\partial {\mathds{E}}_{\rm BVP}}{\partial (\mathds{D}_0 {\bf X}^\mu)} + ({\bm h}_\sigma^T \tilde{\bm \lambda}_\mu) \mathfrak{d}^0[0] + ( {\bm h}_\sigma^T \tilde{\bm \gamma}_\mu) \mathfrak{d}^0[N_0] \Big)\delta {\bf X}^\mu.\label{eq:NoetherChargeDiscrLagr}
\end{align}

When considering the concrete implementation of our discretization strategy for wave propagation in $(1+1)$ dimensions in the following section, we will show that the Noether charge density ${\bm q}^{\rm L}$ indeed leads to a Noether charge ${\bm Q}^{\rm L}$ which is exactly conserved over the whole duration of the simulation.

The Lagrange multiplier terms also affect the form of the equations of motion, leading to deviations from the continuum form at the last and second to last $\tau$ slice. However, we will show that, as in the world-line formalism, these modifications do not spoil the convergence of the solution to the correct continuum limit nor do they interfere with the conservation of the Noether charge.

\section{Proof-of-principle: wave propagation in $(1+1)$ dimensions}
\label{sec:proofofprinciple11d}
Having established the general variational approach to IBVPs with dynamical coordinate maps in $(d+1)$ dimensions in the previous sections, we proceed to present an explicit implementation for a dynamical system of central interest, scalar wave propagation in $(1+1)$ dimensional flat space-time. It serves as a demonstration for how our novel approach addresses two and connects to the third of the three central challenges to discretized IBVPs. In the context of the first challenge, we show how space-time symmetries can be retained after discretization and that the continuum Noether charges stay exactly preserved in the discrete context. Taking aim at the second challenge, we show how a non-constant discretization of the space-time coordinates emerges dynamically, realizing a form of automatic adaptive mesh refinement (AMR). Along the way we will see how the presence of the coordinate maps offers new freedom in the implementation of boundary conditions, connecting to the third challenge.

In order to keep computational cost for this proof-of-principle to a minimum, we consider $(1+1)$ dimensions, where the abstract parameter space is two dimensional $\Sigma^a=(\tau,\sigma)^a$ and the most general coordinate maps read $X^0(\tau,\sigma)=t(\tau,\sigma)$ and $X^1(\tau,\sigma)=x(\tau,\sigma)$. In $(1+1)$ dimensions the continuum theory features three space-time symmetries: translations in time $t$, translations in space $x$, and relativistic rotations involving both the time and space coordinate, so called boosts. I.e. even though spatial rotation are absent in $(1+1)$ dimensions, our proof-of-principle contains a space-time symmetry structure that is similar to the complexity found in higher dimensions. The value of the scale $T$ is taken to be large compared to the energy density stored in the initial conditions of the field degrees of freedom.

In order to avoid complications in numerically locating the critical point of the IBVP action, %due to the invariance of the action under a large class of coordinate transformations, 
we will consider a trivial mapping in space $x(\tau,\sigma)=\sigma$ while leaving the temporal mapping fully dynamic $t=t(\tau,\sigma)$. This choice allows us to focus on the preservation of time translation symmetry after discretization, associated with the concept of a conserved energy. %We however do not maintain spatial translation and boost symmetry. 
Conservation of energy, while not only a central physical tenet of the classical dynamics, is vital to the stability of the discretization.

\subsection{Continuum formulation}
\label{sec:1p1dcontform}
In flat (Minkowski-)space-time in $(1+1)$ dimensions, the metric tensor takes the form $G={\rm diag}[c^2,-1]$. In the presence of dynamical coordinate maps we must determine the corresponding induced metric $g$ on the abstract parameter manifold spanned by the $(\tau,\sigma)$ directions. Let us denote derivatives with respect to $\tau$ as dots, as in $\dot t\equiv\partial t/\partial\tau $ and $\dot x\equiv \partial x/\partial\tau $, and we will indicate derivatives w.r.t. the spatial parameter $\sigma$ as primes $t^\prime\equiv \partial t/\partial\sigma$ and $ x^\prime\equiv \partial x/\partial\sigma$. The explicit expression for the Jacobian of the coordinate maps then is
\begin{align}
J^\mu_a=\partial X^\mu/\partial \Sigma^a=\left(\begin{array}{cc} \dot t& t^\prime \\ \dot x& x^\prime \end{array}\right).
\end{align}
The induced metric $g=J^TGJ$ and its adjugate thus read
\begin{align}
    g=\left( \begin{array}{cc} c^2\td^2-\xd^2 & c^2\td\tp - \xd\xp \\c^2 \tp\td - \xp\xd & c^2(\tp)^2-(\xp)^2\end{array}\right), \; {\rm adj}[g]= \left( \begin{array}{cc} c^2(\tp)^2-(\xp)^2 & \xd\xp-c^2\td\tp\\ \xp\xd-c^2 \tp\td & c^2\td^2-\xd^2\end{array}\right),
\end{align}
and the determinant reduces to
\begin{align}
    {\rm det}[g]= -c^2(\td\xp-\xd\tp)^2\label{eq:fulldetg}.
\end{align}
It is remarkable that the full determinant reduces to such a simple expression.

Setting the potential term to zero as $V(\phi)=0$ in our novel action ${\cal S}_{\rm BVP}$ (introduced in \cref{eq:novelaction}) corresponds to a system that describes the propagation of scalar waves with velocity $c$. Derivatives of the dynamical field $\phi(\tau,\sigma)$ are denoted similarly to those of the coordinate maps as $\dot \phi\equiv\partial \phi/\partial\tau $ and $\phi^\prime\equiv \partial \phi/\partial\sigma$. For the $(1+1)$ dimensional setting we thus have
\begin{align}
 \nonumber {\cal S}_{\rm BVP}=&\int d\tau d\sigma \, \big(-T\big) \sqrt{ \Big(-{\rm det}[g] + \frac{1}{T} \partial_a\phi(\Sigma)\partial_b\phi(\Sigma) {\rm adj}[g]_{ab}}\\
=&\int d\tau d\sigma \, \big(-T\big) \Big\{ c^2(\td\xp-\xd\tp)^2\label{eq:novelaction1p1}\\
\nonumber+&\frac{1}{T} \Big( \pd^2(c^2(\tp)^2-(\xp)^2) + 2\pd\pp(\xd\xp-c^2\td\tp) + (\pp)^2(c^2\td^2-\xd^2) \Big) \Big\}^{1/2}. 
\end{align}
Note how each term under the square root contains exactly two $\tau$ and two $\sigma$ derivatives. This means that the integral is not only invariant under a differentiable redefinition of $(\tau,\sigma)\to f(\tau,\sigma)$ but that at the same time the dimensions of the parameters $(\tau,\sigma)$ are irrelevant; the dimensions of the parameters cancel between the measure and the derivative terms. Whenever temporal and spatial coordinate maps appear in \cref{eq:novelaction1p1} they do so quadratically with a relative minus sign, which encodes the proper invariance of the expressions under Lorentz transformations. As only dotted and primed coordinates contribute, global shifts in time and space too constitute manifest symmetries.
In the following we will, for notational clarity, choose natural units such that $c=1$. Specifying to the trivial spatial mapping $x(\tau,\sigma)=\sigma$ and removing the square root from the integrand we obtain the following action on which we base our proof-of-principle 
 \begin{align}
{\cal E}_{\rm BVP}\overset{x=\sigma}{=}\int &d\tau d\sigma \, \frac{1}{2} \Big\{ (\td)^2+\frac{1}{T} \Big( \pd^2((\tp)^2-1) - 2 \pd\pp\td\tp + (\pp)^2(\td^2) \Big) \Big\}. \label{eq:novelactionE1p1}
\end{align}
Note that even with the trivial mapping $x=\sigma$ we still retain manifest symmetry under global time translation, since each reference to the $t$-mapping in \cref{eq:novelactionE1p1} involves either its $\tau$ or $\sigma$ derivative. Let us derive the equations of motion for this system according to (\ref{eq:varEeom}).

The variation of this action leads to
\begin{align}
\nonumber\delta {\cal E}_{\rm BVP}= \int d\tau d\sigma \Big\{& \td\delta\td + \frac{1}{T} \Big( \pd\delta\pd((\tp)^2-1) +\pd^2\tp\delta\tp - \delta\pd\pp\td\tp - \pd\delta\pp\td\tp\\
& -\pd\pp\delta\td \tp - \pd\pp\td\delta\tp + \pp\delta\pp\td^2+(\pp)^2\td\delta\td\Big)\Big\}.
\end{align}
Carrying out integration by parts results in 
\begin{align}
\delta {\cal E}_{\rm BVP}= &-\int d\tau d\sigma \Big\{\Big[ \ddot t + \frac{1}{T}\Big\{ \frac{\partial}{\partial \sigma}\Big( \pd^2\tp - \pd\pp\td \Big) + \frac{\partial}{\partial \tau} \Big( (\pp)^2\td - \pd\pp\tp \Big)\Big\}\Big] \delta t\label{eq:1p1dvareomt}\\
&+ \Big[ \frac{\partial}{\partial \tau} \Big( \pd ((\tp)^2-1) -\pp\td\tp \Big) + \frac{\partial}{\partial \sigma} \Big( \pp \td^2 - \pd\td\tp \Big) \Big]\delta \phi \Big\}\label{eq:1p1dvareomphi}\\
\nonumber &+\int d\sigma \Big\{ \Big[ \dot t + \Big( (\pp)^2\td - \pd\pp\tp \Big)\Big] \delta t\Big\} \Big|_{\tau^{\rm i}}^{\tau^{\rm f}} + \int d\tau \Big\{ \Big[ \pd^2\tp - \pd\pp\td\Big]\delta t   \Big\}\Big|_{\sigma^{\rm i}}^{\sigma^{\rm f}}\\
\nonumber &+\int d\sigma \Big\{ \Big[ \pd ((\tp)^2-1) -\pp\td\tp \Big] \delta \phi\Big\} \Big|_{\tau^{\rm i}}^{\tau^{\rm f}} + \int d\tau \Big\{ \Big[ \pp \td^2 - \pd\td\tp  \Big]\delta \phi   \Big\}\Big|_{\sigma^{\rm i}}^{\sigma^{\rm f}},
\end{align}
from which we can read off the coupled equations of motion for $t$ and $\phi$ from (\ref{eq:1p1dvareomt}) and (\ref{eq:1p1dvareomphi}) respectively
\begin{align}
 \ddot t + \frac{1}{T}\Big\{ \frac{\partial}{\partial \sigma}\Big( \pd^2\tp - \pd\pp\td \Big) + \frac{\partial}{\partial \tau} \Big( (\pp)^2\td - \pd\pp\tp \Big)\Big\}&=0,\label{eq:1p1deomt}\\
 \frac{\partial}{\partial \tau} \Big( \pd ((\tp)^2-1) -\pp\td\tp \Big) + \frac{\partial}{\partial \sigma} \Big( \pp \td^2 - \pd\td\tp \Big)&=0.\label{eq:1p1deomphi}
\end{align}
We may recover the conventional wave equation from \cref{eq:1p1deomphi} after re-instituting the explicit factors of $c$, by restricting to a trivial time map $t(\tau,\sigma)=\tau$
\begin{align}
 0&\quad=\quad \frac{\partial}{\partial \tau} \Big( \pd (c^2(\tp)^2-1) -c^2\pp\td\tp \Big) + \frac{\partial}{\partial \sigma} \Big( c^2 \pp \td^2 - c^2\pd\td\tp \Big)\\
 &\overset{t=\tau, x=\sigma}{=}-\frac{\partial}{\partial t} \Big( \pd \Big) + \frac{\partial}{\partial x} \Big( c^2 \pp \Big).
\end{align}

In the limit of $T\to\infty$ we would recover from \cref{eq:1p1deomt} the simple relation $\ddot t(\tau,\sigma)=0$ for the time mapping, which is solved by a sheet, constant in the $\sigma$ direction and with constant slope $\dot t_{\rm IC}$ along $\tau$. Note that for $T\to\infty$ the action only contains reference to $\td^2$, meaning that it decomposes into independent temporal maps for each value of $\sigma$. As there are no terms involving spatial derivatives $\tp$ in the action, also no spatial boundary terms arise when deriving the governing equations. In the presence of fields the dynamics of $\phi$ couples to both $\tau$ and $\sigma$ derivatives of $t$, leading to non-trivial dynamics of the time mapping. 

In contrast to the conventional wave equation, the equation of motion (\ref{eq:1p1deomphi}) for the field here not only contains second derivatives in $\tau$ and $\sigma$ but due to the dynamical nature of the $t$ mapping, also mixed terms. 

Let us take a look at the spatial boundary conditions necessary for achieving equivalence between the equations of motion \cref{eq:1p1deomt,eq:1p1deomphi} and the critical point of the action \cref{eq:novelactionE1p1}. Note that the temporal boundary conditions will be taken care of by the forward backward construction, when going over from ${\cal E}_{\rm BVP}\to{\cal E}_{\rm IBVP}$.
We thus find that the following two terms must vanish
\begin{align}
    \int d\tau \Big\{ \Big[ \pd^2\tp - \pd\pp\td\Big]\delta t   \Big\}\Big|_{\sigma^{\rm i}}^{\sigma^{\rm f}}\overset{!}{=}0, \qquad \int d\tau \Big\{ \Big[ \pp \td^2 - \pd\td\tp  \Big]\delta \phi   \Big\}\Big|_{\sigma^{\rm i}}^{\sigma^{\rm f}}\overset{!}{=}0.\label{eq:bndcnd1p1}
\end{align}
In contrast to the conventional formulation of IBVPs we now have both $\phi$ and $t$ appearing in these boundary expressions. This amounts to a \textit{new freedom available} to us in order to construct boundary conditions, such that the boundary terms vanish. 

In our proof-of-principle here, we simply choose vanishing spatial Dirichlet boundary conditions for the field $\phi$. This entails that $\phi$ and $\delta \phi$, as well as $\dot \phi$ (due to the time independence of Dirichlet boundary conditions), are all zero on the spatial boundary. This choice leaves us with vanishing boundary terms in \cref{eq:bndcnd1p1} without the need to specify a boundary condition on $t$. Conversely one could imagine boundary conditions enforced for $t$, which would make \eqref{eq:bndcnd1p1} vanish without the need to provide an explicit boundary condition on $\phi$.
I.e. the presence of the coordinate maps provides a \textit{unique and novel flexibility in the treatment of boundary conditions for the dynamical fields}, especially interesting for non-reflecting boundary conditions.

Having confirmed that for our choice of boundary conditions (\ref{eq:bndcnd1p1}) vanishes, we have established the equations of  motion for our $(1+1)$ dimensional wave propagation.

Following \cref{eq:NoetherCur}, the components of the conserved Noether current $\partial_a J_t^a = 0$ associated with manifest invariance under infinitesimal time translations $\delta t$, i.e. $\delta X^\mu=\delta t \delta^{\mu 0}$ read
\begin{align}
    &J_t^0=\frac{\partial E_{\rm BVP}}{\partial (\td) } = \td + \frac{1}{T}\Big( (\pp)^2\td - \pd\pp\tp \Big)\\
    &J_t^1=\frac{\partial E_{\rm BVP}}{\partial (\tp)} =  \frac{1}{T}\Big( (\pd)^2\tp - \pd\pp\td \Big)
\end{align}
The zeroth component of a conserved current defines the Noether charge according to \cref{eq:defQ}, which yields
\begin{align}
    Q_t(\tau)=\int d\sigma J_t^0 = \int d\sigma \Big\{ \td + \frac{1}{T}\Big( (\pp)^2\td - \pd\pp\tp \Big)\Big\}.\label{eq:1p1dnoethcont}
\end{align}
This quantity in the presence of dynamical coordinate maps constitutes the energy, which depends both on contributions from the fields and from the coordinate maps themselves. 

Let us confirm that our choice of boundary condition conserves the Noether charge
\begin{align}
    \nonumber \partial_\tau Q_t &\overset{{\rm def.}\, \cref{eq:1p1dnoethcont}}{=} \int d\sigma \frac{\partial }{\partial \tau} \Big\{ \td + \frac{1}{T}\Big( (\pp)^2\td - \pd\pp\tp \Big)\Big\}\\
    \nonumber &\overset{ \partial_0J_t^0=-\partial_1J_t^1}{=} -\int d\sigma \frac{\partial }{\partial \sigma} \Big\{ \frac{1}{T}\Big( \pd^2\tp - \pd\pp\td\Big) \Big\}\\
    &\qquad = -\frac{1}{T}\Big( \pd^2\tp - \pd\pp\td\Big)\Big|_{\sigma_{\rm i}}^{\sigma^{\rm f}}=0.
\end{align}
In the second line we have exploited that the Noether current is conserved by construction. And indeed due to the vanishing of $\phi$ on the spatial boundaries and thus $\dot \phi$ being zero there, the Noether charge does not change over time. 

Note that the Noether charge of \cref{eq:1p1dnoethcont}, in contrast to its counterpart in the world-line formalism, makes reference to all dynamical degrees of freedom, i.e. both their spatial and temporal derivatives. This places tight constraints on the growth of derivatives, and while not a proof, suggests that the evolution of the system is stable. This will be considered in future work. 

When we wish to formulate a causal IBVP we have to go over to the forward backward construction and incorporate initial and boundary conditions, as well as the connecting conditions for both the field and temporal mapping via explicit Lagrange multipliers. Note that the evolution of the time mapping too is handled as an IBVP, meaning we do not specify an initial and final time but instead an initial time $t_{\rm IC}$ and the derivative of time $\dot t_{\rm IC}$ with respect to the temporal parameter $\tau$. Similarly we have to supply $\phi_{\rm IC}$ and $\dot \phi_{\rm IC}$ at the initial $\tau$ slice. At this point we also add all possible Lagrange multiplier functions $\kappa,\tilde \kappa,\xi, \tilde \xi$ related to the spatial (Dirichlet) boundary conditions
\begin{align}
{\cal E}^{\rm L}_{\rm IBVP}=&\int d\tau d\sigma \, \frac{1}{2} \Big\{ (\td_1)^2+\frac{1}{T} \Big( \pd^2_1((\tp_1)^2-1) - 2 \pd_1\pp_1\td_1\tp_1 + (\pp_1)^2(\td_1^2) \Big) \Big\}\label{eq:1p1dEIBVPc}\\
\nonumber -&\int d\tau d\sigma \, \frac{1}{2} \Big\{ (\td_2)^2+\frac{1}{T} \Big( \pd_2^2((\tp_2)^2-1) - 2 \pd_2\pp_2\td_2\tp_2 + (\pp_2)^2(\td_2^2) \Big) \Big\}\\
\nonumber +&\int d\sigma \Big\{ \lambda^t(\sigma)\Big(t_1(\tau^{\rm i},\sigma)-t_{\rm IC}(\sigma)\Big)+ \lambda^\phi(\sigma)\Big(\phi_1(\tau^{\rm i},\sigma)-\phi_{\rm IC}(\sigma)\Big)\Big\}\\
\nonumber +&\int d\sigma \Big\{ \tilde \lambda^t(\sigma)\Big(\dot t_1(\tau^{\rm i},\sigma)-\dot t_{\rm IC}(\sigma)\Big)+ \tilde \lambda^\phi(\sigma)\Big(\dot \phi_1(\tau^{\rm i},\sigma)-\dot \phi_{\rm IC}(\sigma)\Big)\Big\}\\
\nonumber +&\int d\sigma \Big\{ \gamma^t(\sigma)\Big(t_1(\tau^{\rm f},\sigma)-t_2(\tau^{\rm f},\sigma)\Big)+ \gamma^\phi(\sigma)\Big(\phi_1(\tau^{\rm f},\sigma)-\phi_2(\tau^{\rm f},\sigma)\Big)\Big\}\\
\nonumber +&\int d\sigma \Big\{ \tilde\gamma^t(\sigma)\Big(\dot t_1(\tau^{\rm f},\sigma)-\dot t_2(\tau^{\rm f},\sigma)\Big)+ \tilde \gamma^\phi(\sigma)\Big(\dot \phi_1(\tau^{\rm f},\sigma)-\dot \phi_2(\tau^{\rm f},\sigma)\Big)\Big\}\\
\nonumber +&\int d\tau \Big\{ \kappa^{t}(\tau)\Big(t_1(\tau,\sigma^{\rm i})-t_{\rm sBCL}(\tau)\Big)+ \kappa^{\phi}(\tau)\Big(\phi_1(\tau,\sigma^{\rm i})-\phi_{\rm sBCL}(\tau)\Big)\Big\}\\
\nonumber +&\int d\tau \Big\{ \tilde\kappa^t(\tau)\Big(t_1(\tau,\sigma^{\rm f})-t_{\rm sBCR}(\tau)\Big)+ \tilde\kappa^\phi(\tau)\Big(\phi_1(\tau,\sigma^{\rm f})-\phi_{\rm sBCR}(\tau)\Big)\Big\}\\
\nonumber +&\int d\tau \Big\{ \xi^{t}(\tau)\Big(t_2(\tau,\sigma^{\rm i})-t_{\rm sBCL}(\tau)\Big)+ \xi^{\phi}(\tau)\Big(\phi_2(\tau,\sigma^{\rm i})-\phi_{\rm sBCL}(\tau)\Big)\Big\}\\
\nonumber +&\int d\tau \Big\{ \tilde\xi^t(\tau)\Big(t_2(\tau,\sigma^{\rm f})-t_{\rm sBCR}(\tau)\Big)+ \tilde\xi^\phi(\tau)\Big(\phi_2(\tau,\sigma^{\rm f})-\phi_{\rm sBCR}(\tau)\Big)\Big\}.
\end{align}

The presence of Lagrange multipliers leads to a modification of the expression for the energy. In the physical limit the values of the degrees of freedom on the forward and backward branch must agree, hence we define the Noether charge simply from the forward time mapping and field
\begin{align}
    Q_t^{\rm L}(\tau)=\int d\sigma \Big\{ \td_1 + \frac{1}{T}\Big( (\pp_1)^2\td_1 - \pd_1\pp_1\tp_1 \Big) + \tilde\lambda^t\delta(\tau-\tau^{\rm i})  + \tilde\gamma^t\delta(\tau-\tau^{\rm f}) \Big\}.\label{eq:1p1dnoethcontl}
\end{align}

With the full continuum IBVP action defined in \cref{eq:1p1dEIBVPc} we are now ready to discretize the action and obtain the classical trajectory by numerically determining its critical point.

\subsection{ Discretized formulation}
\label{sec:discr2d}

In the following we discretize the two-dimensional parameter space $(\tau,\sigma)$ using a grid with $(N_\tau,N_\sigma)$ points. We choose finite intervals for the values of $\tau\in[\tau^{\rm i},\tau^{\rm f}]$ and $\sigma\in[\sigma^{\rm i},\sigma^{\rm f}]$ leading to $\Delta \tau= (\tau^{\rm f}-\tau^{\rm i})/(N_\tau-1)$ and $\Delta \sigma= (\sigma^{\rm f}-\sigma^{\rm i})/(N_\sigma-1)$. The arrays to house the dynamical degrees of freedom are thus ${\bf t}_{1,2}$ and ${\bm \phi}_{1,2}$ each with \texttt{TotVol}$=N_\tau N_\sigma$ entries\footnote{The spatial mapping was chosen to be trivial $x=\sigma$ and thus does not require to be treated as a separate array.}. For reference to subsets of these arrays let us introduce projection operators $\mathds{P}_a^k$. When applied to an array of \texttt{TotVol} entries, they produce the $k$-th slice in the $a$ direction, e.g. $\mathds{P}_\tau^0[{\bm \phi}]=(\phi(0,0), \phi(0,\Delta\sigma),\ldots,\phi(0,\sigma^{\rm f}))\in \mathbb R^{N_\sigma}$.

As discussed in \cref{sec:discrIBVP}, we deploy summation by parts operators $\mathds{D}_\tau=\mathbb{d}_\tau\otimes{\mathbb{i}}$ and $\mathds{D}_\sigma={\mathbb{i}}\otimes\mathbb{d}_\sigma$, as well as approximate integration by the compatible quadrature matrix $\mathds{H}=\mathbb{h}_\tau\otimes\mathbb{h}_\sigma$. The array with the diagonal entries of $\mathds{H}$ we refer to as ${\bm h}$, which is used to express integration of a function $\bm f$ as $({\bm f},{\bm 1})={\bm f}^T \mathds{H} {\bm 1}= {\bm f}^T {\bm h}$. For the temporal direction, where we have a genuine initial value problem, as the solution at the final point $\tau^{\rm f}$ is \textit{a priori} unknown, we introduce the previously discussed regularization in the SBP operators. This leads to $\bar{\mathds{D}}_{\tau}^t$ and $\bar{\mathds{D}}_{\tau}^\phi$, whose regularization is based on penalty terms including reference to the initial values. We found that when there are spatial Dirichlet boundary conditions enforced, we do not need to introduce regularization, as the $\pi$-mode is avoided automatically by requiring the function to take a specific value at both ends of the interval.

The Lagrange multiplier arrays are referred to as ${\bm \lambda}^t,\tilde{{\bm \lambda}}^t,{\bm \lambda}^\phi,\tilde{{\bm \lambda}}^\phi$ for the initial conditions stored in ${\bm t}_{\rm IC},{\bm \phi}_{\rm IC}$ and $\dot {\bm t}_{\rm IC},\dot {\bm \phi}_{\rm IC}$ arrays. Note that these arrays only contain $N_\sigma$ entries, as they make reference solely to the first $\tau$ slice. Similarly we have ${\bm \gamma}^t,\tilde{{\bm \gamma}}^t,{\bm \gamma}^\phi,\tilde{{\bm \gamma}}^\phi$ for the connecting conditions. As we saw in \cref{eq:bndcnd1p1} if we supply vanishing spatial Dirichlet boundary data to the field on the forward and backward branch using $N_\tau$ component arrays ${\bm \kappa}^\phi,\tilde{{\bm \kappa}}^\phi,{\bm \xi}^\phi,\tilde{{\bm \xi}}^\phi$ no spatial boundary conditions need to be supplied for ${\bm t}$. 

We arrive at the following discretized action functional
\begin{align}
    \nonumber \mathds{E}_{\rm IBVP}^{\rm L} =& \frac{1}{2}\Big\{ (\bar{\mathds{D}}_\tau^t {\bm t}_1)^2+\frac{1}{T} \Big( (\bar{\mathds{D}}_\tau^\phi {\bm \phi}_1)^2\circ\big( (\mathds{D}_\sigma {\bm t}_1)^2 - 1 \big)\\
    \nonumber &\qquad -2(\bar{\mathds{D}}_\sigma^\phi {\bm \phi}_1)\circ(\bar{\mathds{D}}_\tau^\phi {\bm \phi}_1)\circ (\bar{\mathds{D}}_\tau^t {\bm t}_1)\circ(\mathds{D}_\sigma^t {\bm t}_1)  + (\bar{\mathds{D}}_\sigma^\phi {\bm \phi}_1)^2\circ (\bar{\mathds{D}}_\tau^t {\bm t}_1)^2 \Big)\Big\}^T {\bm h}\\
    \nonumber-&\frac{1}{2}\Big\{ (\bar{\mathds{D}}_\tau^t {\bm t}_2)^2+\frac{1}{T} \Big( (\bar{\mathds{D}}_\tau^\phi {\bm \phi}_2)^2\circ \big( (\mathds{D}_\sigma {\bm t}_2)^2 - 1 \big)\\
     \nonumber &\qquad -2(\bar{\mathds{D}}_\sigma^\phi {\bm \phi}_2)\circ(\bar{\mathds{D}}_\tau^\phi {\bm \phi}_2)\circ (\bar{\mathds{D}}_\tau^t {\bm t}_2)\circ(\mathds{D}_\sigma^t {\bm t}_2) + (\bar{\mathds{D}}_\sigma^\phi {\bm \phi}_2)^2 \circ(\bar{\mathds{D}}_\tau^t {\bm t}_2)^2\Big)\Big\}^T {\bm h}\\
    \nonumber+& \big( {\bm \lambda}^t \big)^T \mathbb{h_\sigma} \big( \mathds{P}_\tau^0[{\bm t}_1] - {\bm t}_{\rm IC} \big) + \big( {\bm \lambda}^\phi \big)^T \mathbb{h_\sigma} \big( \mathds{P}_\tau^0[{\bm \phi}_1] - {\bm \phi}_{\rm IC} \big)\\
    \nonumber+& \big( \tilde{{\bm \lambda}}^t \big)^T \mathbb{h_\sigma} \big( \mathds{P}_\tau^0[(\mathds{D}_\tau{\bm t}_1)] - \dot{\bm t}_{\rm IC} \big) + \big( \tilde{\bm \lambda}^\phi \big)^T \mathbb{h_\sigma} \big( \mathds{P}_\tau^0[(\mathds{D}_\tau{\bm \phi}_1)] - \dot{\bm \phi}_{\rm IC} \big)\\
    \nonumber+& \big( {\bm \gamma}^t \big)^T \mathbb{h_\sigma} \big( \mathds{P}_\tau^{N_\tau}[{\bm t}_1] -  \mathds{P}_\tau^{N_\tau}[{\bm t}_2] \big) + \big( {\bm \gamma}^\phi \big)^T \mathbb{h_\sigma} \big(  \mathds{P}_\tau^{N_\tau}[{\bm \phi}_1] -  \mathds{P}_\tau^{N_\tau}[{\bm \phi}_2] \big)\\
    \nonumber+& \big( \tilde{\bm \gamma}^t \big)^T \mathbb{h_\sigma} \big(  \mathds{P}_\tau^{N_\tau}[(\mathds{D}_\tau{\bm t}_1)] -  \mathds{P}_\tau^{N_\tau}[(\mathds{D}_\tau{\bm t}_2)] \big)\\
    \nonumber+& \big( \tilde{\bm \gamma}^\phi \big)^T \mathbb{h_\sigma} \big(  \mathds{P}_\tau^{N_\tau}[(\mathds{D}_\tau{\bm \phi}_1)] -  \mathds{P}_\tau^{N_\tau}[(\mathds{D}_\tau{\bm \phi}_2)] \big)\\
    \nonumber+& \big( {\bm \kappa}^\phi \big)^T \mathbb{h_\tau} \big(  \mathds{P}_\sigma^{0}[{\bm \phi}_1] - {\bm 0} \big) + \big( \tilde {\bm \kappa}^\phi \big)^T \mathbb{h_\tau} \big(  \mathds{P}_\sigma^{N_\sigma}[{\bm \phi}_1] - {\bm 0} \big)\\
    +& \big( {\bm \xi}^\phi \big)^T \mathbb{h_\tau} \big(  \mathds{P}_\sigma^{0}[{\bm \phi}_2] - {\bm 0} \big) + \big( \tilde {\bm \xi}^\phi \big)^T \mathbb{h_\tau} \big(  \mathds{P}_\sigma^{N_\sigma}[{\bm \phi}_2] - {\bm 0} \big).\label{eq:discrEsim}
\end{align}

For the discrete Noether charge corresponding to \cref{eq:1p1dnoethcontl}, we must implement spatial integration as described in \cref{eq:spatint}. To this end we construct the $N_\tau\times (N_\tau \cdot N_\sigma)$ matrix $\mathds{H}_\sigma$. Applying $\mathds{H}_\sigma$ to a discretized function $\bm f$ results in an array of $N_\tau$ length. Here the $N_\sigma$ component vector $({\bm h}_\sigma)_i=(\mathbb{h}_\sigma)_{ii}$, which contains the diagonal entries of the spatial quadrature matrix $\mathbb{h}_\sigma$. The matrix $\mathds{H}_\sigma$ then has the form (\ref{eq:spatint}). 

To obtain the array for the Noether charge with $N_\tau$ entries we thus apply
\begin{align}
   \nonumber {\bm Q}_t^{\rm L}= \mathds{H}_\sigma &\underbracket{\Big\{ ( \mathds{D}_\tau{\bm t_1}) + \frac{1}{T} \Big( (\mathds{D}_\sigma {\bm \phi}_1)^2\circ(\mathds{D}_\tau{\bm t_1)}- (\mathds{D}_\tau {\bm \phi}_1)\circ(\mathds{D}_\sigma {\bm \phi}_1)\circ(\mathds{D}_\sigma {\bm t}_1)\Big)\Big\}}_{{\bm J}^0\in {\mathbb{R}^{N_\tau\times N_\sigma}}}\\
  +  & \underbracket{\Big\{ ({\bm h}_\sigma^T \tilde{\bm \lambda}^t)  \mathfrak{d}^\tau[0]+ ( {\bm h}_\sigma^T \tilde{\bm \gamma}^t ) \mathfrak{d}^\tau[N_\tau]\Big\}}_{\rm Lagr.\, mult.\, contrib.},\label{eq:dsicrNoetherl}
\end{align}
where $\mathfrak{d}^\tau[k]$ refers to an array with $N_\tau$ entries, defined in (\ref{eq:defdeltas}).

To confirm quantitatively the conservation of this Noether charge, let us introduce $\Delta {\bm E}= {\bm Q}_t^{\rm L} - {\bm Q}_t^{\rm L}[0]$ the difference between the discrete Noether charge ${\bm Q}_t^{\rm L}$ and its initial value, which, due to its dependence on only the derivatives of the degrees of freedom, is fully specified by the initial conditions.

%In the interior the solutions we exactly fulfill the discretized equations of motion.  The only deviations occur on the boundary due to the ambiguity in assigning the Lagrange multipliers to the either the forward or backward branch.  We will explicitly demonstrate that these deviations on the boundary do not affect the global convergence to the correct solution under grid refinement.

As we saw for the Noether charge, the presence of Lagrange multiplier terms, implementing the initial and connection conditions, affects the values of our system at the final (and potentially initial) time slice. 
The influence of the Lagrange multipliers manifests itself in the solutions of both the coordinates and fields, which show deviations from the discretized equations of motion on the boundary (for a detailed discussion see ref.~\cite{Rothkopf:2023vki}).
Note that in the interior the solutions we obtain exactly fulfill the discretized equations of motion and the only deviations occur on the boundary. We will explicitly demonstrate that these deviations on the boundary do not affect the global convergence to the correct solution under grid refinement.

\subsection{Numerical results}

We next present numerical results that exemplify how our novel approach addresses two of the three central challenges to formulating and discretizing IBVPs
\begin{enumerate}[label*=\arabic*.]
    \item Due to the presence of a manifest continuum symmetry (time translations) after discretization, the discrete Noether charge ${\bf Q}_t^{\rm L}$ remains preserved over the whole simulation period.\\
    Since our SBP operators exactly mimic integration by parts, the value of the Noether charge remains conserved at its initial value set solely by the initial data.
    \item The dynamical coordinate map $t(\tau,\sigma)$ is not a trivial sheet but exhibits different values of $\tau$ and $\sigma$ derivatives throughout the simulation domain. This amounts to regions of finer and coarser temporal resolution, constituting a dynamically emerging adaptive mesh, i.e. an automatic AMR procedure.
\end{enumerate}

Our simulation domain is defined by the parameters $\tau\in[0,1/2]$ and $\sigma\in[0,1]$.%\footnote{The reason for choosing different intervals lies in the observation that for exactly the same intervals the numerical minimization is comparatively slow. We believe this is related to us naively elevating our regularization of the SBP operators from the one-dimensional scenario to a multidimensional setting and will be addressed in future work.}. 
Choosing $N_\tau=60$ and $N_\sigma=48$ for the presentation of the main results below, we have $\Delta \tau=1/118$ and $\Delta \sigma =1/47$. Our choice of grid points and spacings represents a compromise between capturing relevant wave dynamics and computational cost.

The physical extent of our simulation domain is determined by the choice of coordinate mappings. As we use the trivial mapping $x=\sigma$ in spatial direction, the total x-direction extent is equal to that of $\sigma$. On the other hand the temporal mapping is fully dynamic and we control its extent by prescribing an initial value $t(\tau^{\rm i},\sigma)=0$ and velocity of $\td(\tau^{\rm i},\sigma)=5/2$. 

We choose a value of $T=10^4$ such that deviations from the simple sheet behavior of the time mapping will be of at most on the percent level. In realistic systems, one will encounter an even larger value of $T$, since the corrections to the behavior encoded in the conventional field theory action $S_{\rm ft}$ in \cref{eq:stdFTaction} are minute. What is important is that $T$ is present, so that the coordinate maps can be considered as dynamical. From a practical point of view, as long as $T$ is chosen much larger than the energy density of the field, encoded in the initial conditions, the solution obtained for the field will benefit from the automatic AMR procedure and also converge to a continuum result that is indistinguishable from that of the conventional action.

The dynamical adaptation of the time-mapping leads the simulation to end at a value of the time coordinate close to $t^{\rm f}\approx5/4$. We choose this value, as our wave packages propagate with physical velocity $v=dx/dt=1$. Between $t^{\rm i}$ and $t^{\rm f}$ the wave will have reflected from the hard spatial system boundaries, characterized by $\phi=0$, and will have met again at the center of the simulation domain. This non-trivial dynamics allows us to showcase relevant aspects of our simulation approach.

\begin{figure}
    \includegraphics[scale=0.4]{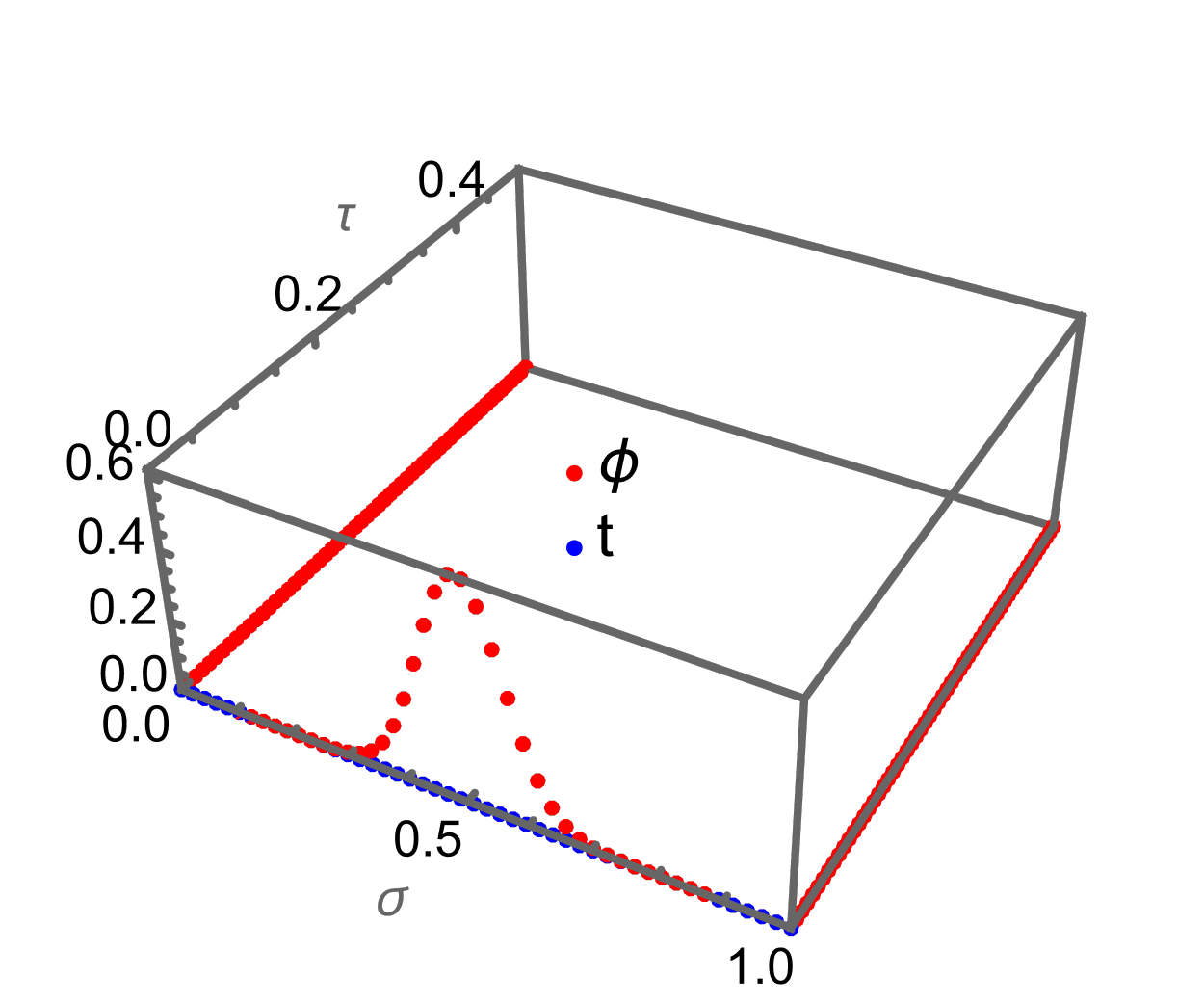}
    \caption{Initial and boundary values provided for the dynamic field $\phi$ as red points and initial values of the time mapping $t$.}
    \label{fig:IBconfPhit}
\end{figure}

We initialize the field $\phi$ with a localized bump, which is symmetrically located in the center of the spatial domain at $\sigma=1/2$ with functional form $\phi(0,\sigma)={\rm sin}[\pi\sigma]{\rm exp}[-100(\sigma-1/2)^2]$. On the spatial boundary we prescribe vanishing Dirichlet boundary conditions. In \cref{fig:IBconfPhit} we plot the initial and boundary values for the dynamical field $\phi$ as red points and the initial value of $t$ as blue points.

In order to mimic the continuum theory, we deploy SBP operators. The examples shown will be based on the \texttt{SBP121} approximation of the derivative. At the end of the section we will provide some comparisons with results based on \texttt{SBP242} operators. 
\begin{figure}
    \includegraphics[scale=0.3]{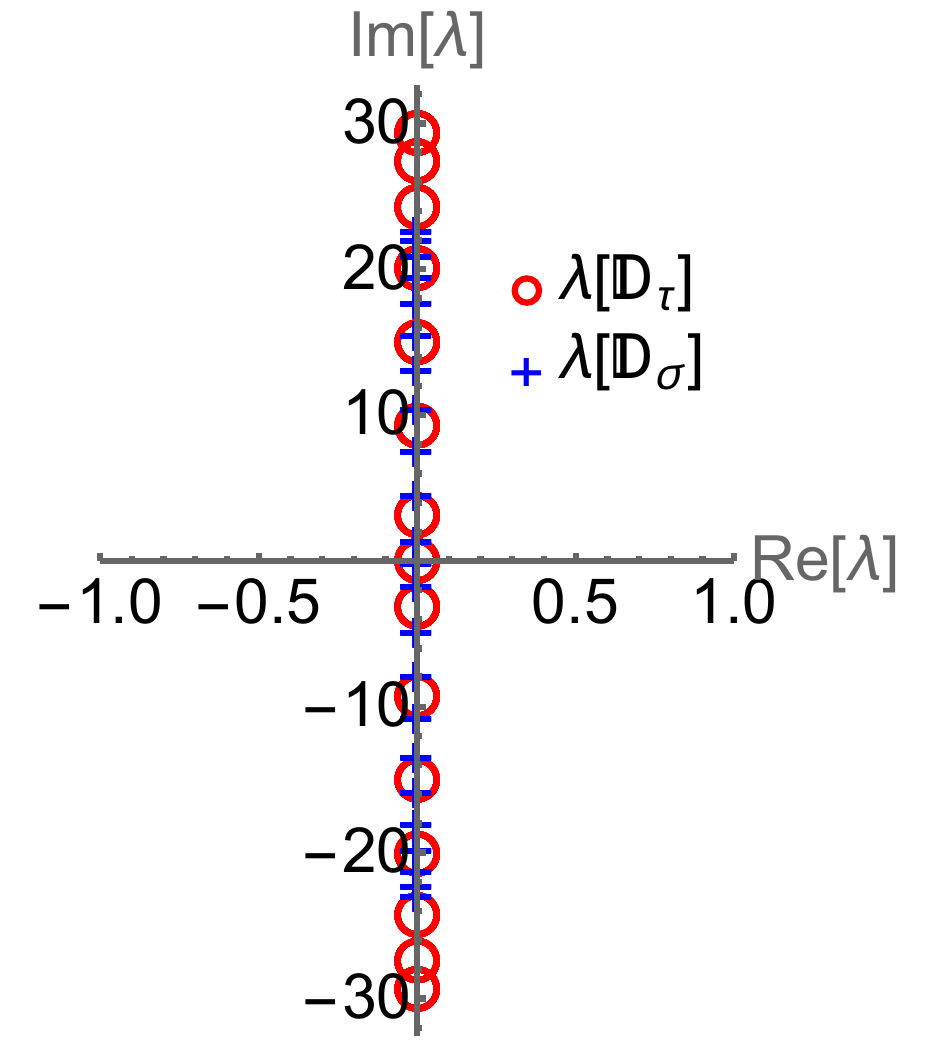}
    \includegraphics[scale=0.3]{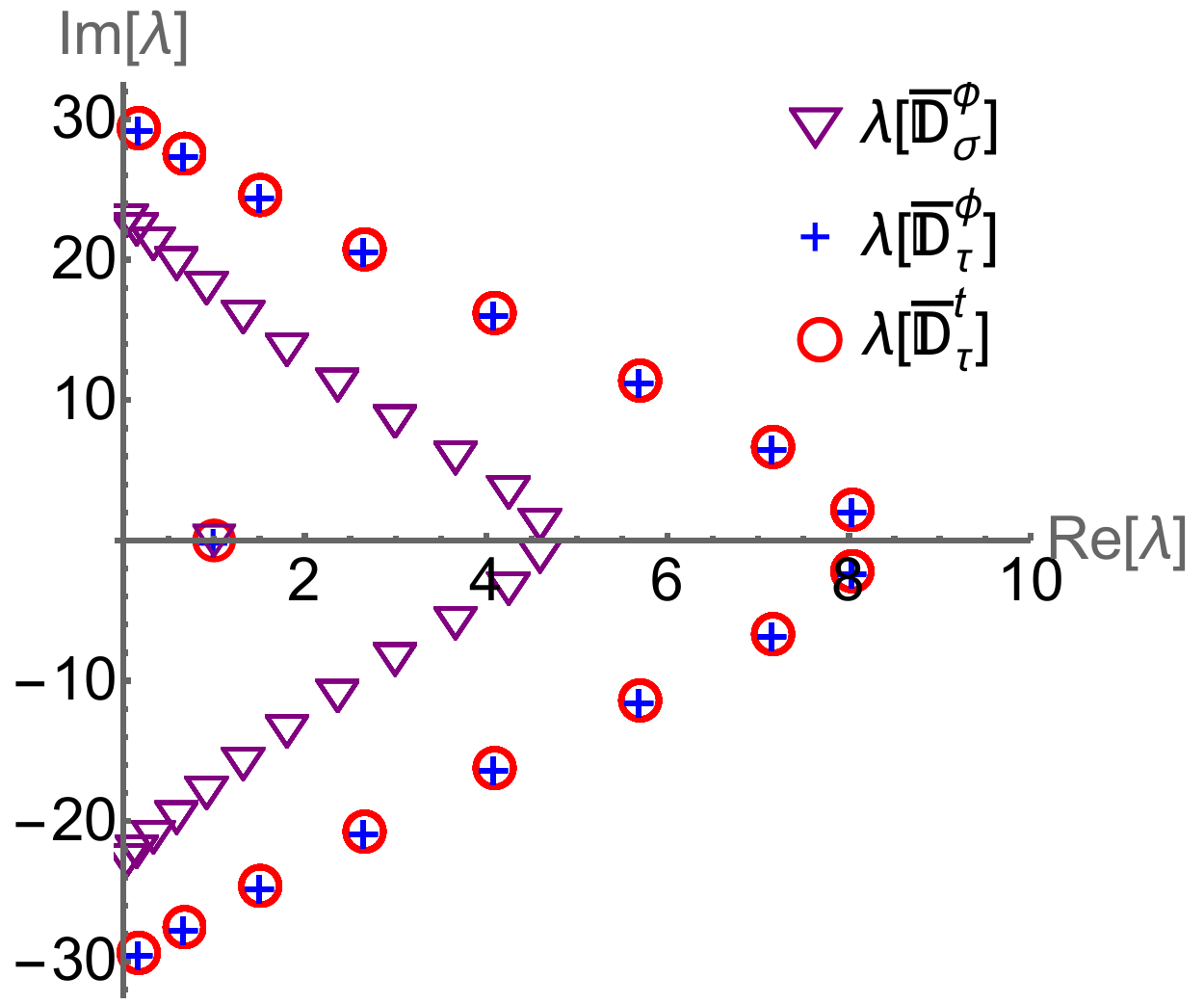}
    \caption{(left) Eigenvalue spectrum of the unregularized \texttt{SBP121} finite difference operators $\mathds{D}_\tau$ (red circles) and $\mathds{D}_\sigma$ (blue crosses) on a grid with $N_\sigma=24$ and $N_\tau=16$. Since $\Delta \tau<\Delta \sigma$ the purely imaginary eigenvalues of $\mathds{D}_\tau$ spread over a larger interval than for the spatial derivative. Each operator has exactly two zero modes depicted at the origin. (right) Eigenvalue spectrum of the regularized \texttt{SBP121} operators (using the initial and boundary data from \cref{fig:IBconfPhit}) in affine coordinates. The zero modes of the SBP operator are lifted and the original physical zero mode, i.e. the constant function, is now associated with the unit eigenvalue.}
    \label{fig:SBP121EValSpec}
\end{figure}
An important ingredient for use of SBP operators in the discretization of IBVP actions is their regularization, in order to avoid the $\pi$-mode. In  \cref{fig:SBP121EValSpec} we provide for completeness the actual eigenvalue spectrum of our \texttt{SBP121} operators for a smaller grid with $N_\sigma=24$ and $N_\tau=16$ points, which already allows us to showcase all relevant properties. The left panel shows the purely imaginary eigenvalues of $\mathds{D}_\tau$ (red circles) and $\mathds{D}_\sigma$ (blue crosses) with their two degenerate zero eigenvalues located at the origin. In the right panel, we instead show the spectrum of the three regularized operators deployed in \cref{eq:discrEsim} each formulated in affine coordinates with a structure as shown in \cref{eq:sketchregSBP}. We have two temporal ones $\bar{\mathds{D}}^t_\tau$ and $\bar{\mathds{D}}^\phi_\tau$, as well as the spatial one $\bar{\mathds{D}}^\phi_\sigma$. The data used for the regularization are the initial and spatial boundary values shown in \cref{fig:IBconfPhit}. Note that all zero modes have been successfully lifted and that the original physical zero mode, the constant function, is now associated with an eigenvalue of unit magnitude.

In order to determine the classical solution we carry out a numerical optimization procedure, locating the critical point of \cref{eq:discrEsim} %\footnote{Let us stress that when operating with a genuine IBVP action, such as $\mathds{E}^{\rm L}_{\rm IBVP}$, there is no need to derive the equations of motion from the action.} 
in the full set of dynamical degrees of freedom ${\bf d}=\{ {\bf t}_{1,2},{\bm \phi}_{1,2}, {\bm \lambda}^t,{\bm \lambda}^\phi,\tilde {\bm \lambda}^t,\tilde{\bm \lambda}^\phi,{\bm \gamma}^t,\tilde {\bm \gamma}^\phi,\tilde {\bm \gamma}^t,{\bm \gamma}^\phi, {\bm \kappa}^\phi, \tilde{\kappa}^\phi,{\bm \xi}^\phi,\tilde{\bm \xi}^\phi\}$.

 For our proof-of-principle study we use the \texttt{Mathematica 13.2}  software\footnote{Example \texttt{Mathematica 13.2} code, based on the \texttt{SBP121} operator is available at the Zenodo repository under open access \cite{rothkopf_2024_11082746}.} and the minimization routines it provides. Note that in the presence of Lagrange multipliers the critical point of the action functional is in general a saddle point and not an extremum. With most established optimization algorithms designed to locate extrema, we first convert the action functional into a suitable form, i.e. we optimize not on $\mathds{E}^{\rm L}_{\rm IBVP}$ but on the norm of its gradient $|{\bm \nabla}_{\bm d}\mathds{E}^{\rm L}_{\rm IBVP}|$. When $\mathds{E}^{\rm L}_{\rm IBVP}$ features a saddle point, then the norm of its gradient exhibits a genuine minimum, since its values must be positive definite. In this way we convert saddle points into extrema, which can be straight forwardly located.

In practice we use multi-step preconditioning to locate the solution. First, in order to obtain appropriate starting values, we optimize on an action where the coordinate maps and field degrees are decoupled, similar to the $T\to\infty$ limit of the full theory. To exploit optimally the different stability and convergence properties of optimization algorithms, we start with the \texttt{QuasiNewton} method of the \texttt{Mathematica} command \texttt{FindMinimum}, iterating for around one to two-thousand steps, followed by application of the interior-point optimization \cite{wright2005interior} based method \texttt{IPOPT}. With this solution in hand we turn to the full action and use the \texttt{IPOPT} method to obtain the final result.

Since the \texttt{IPOPT} method in \texttt{Mathematica} makes reference to an external compiled library, it is unable to benefit from the inbuilt high precision arithmetic, which limits us to find extrema up to residuals of $\Delta\sim 10^{-20}$. While at first this may appear more than satisfactory, it was observed in the study of the world-line formalism in \cite{Rothkopf:2023ljz} that even smaller global residuals are often needed to e.g. enforce the physical limit, i.e. the equality of the degrees of freedom on the forward and backward branch, down to zero in double precision. The reason for the occurrence of such small residuals lies in the intricate cancellations that take place between the degrees of freedom on the forward and backward branch. Hence, as we are limited to using \texttt{IPOT}, our overall precision will be around $\Delta\sim10^{-6}$. I.e. within the accuracy of our minimization %, two functions agree if they differ by less than $\Delta\sim10^{-6}$, in particular 
a numerical zero amounts to values smaller than $\Delta\sim10^{-6}$.

\begin{figure}
    \includegraphics[scale=0.4
]{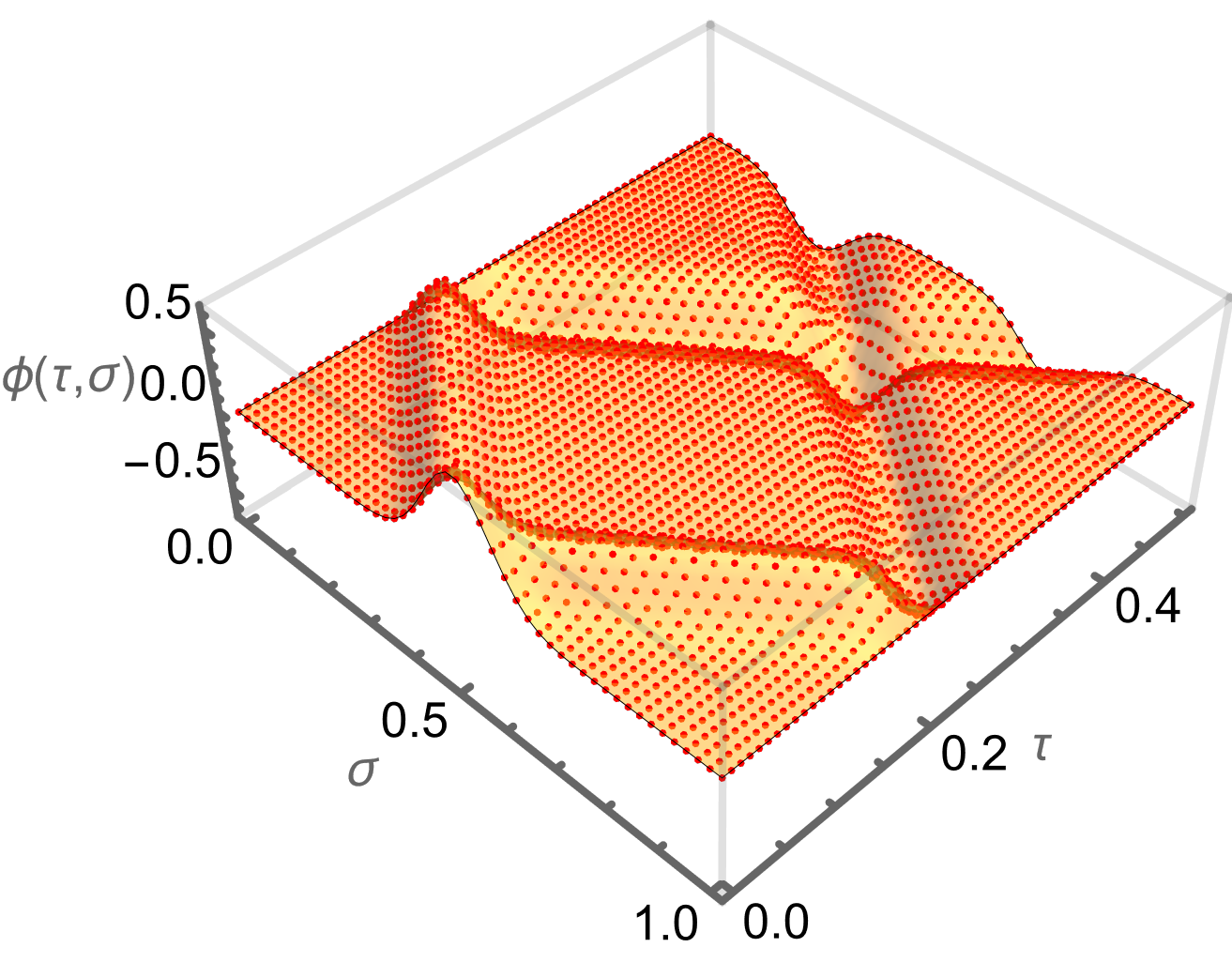}
%     \includegraphics[scale=0.3
% ]{img/PhiSolutionNs48Nt60tmax125SBP121_diffangle.pdf}
    \includegraphics[scale=0.4]{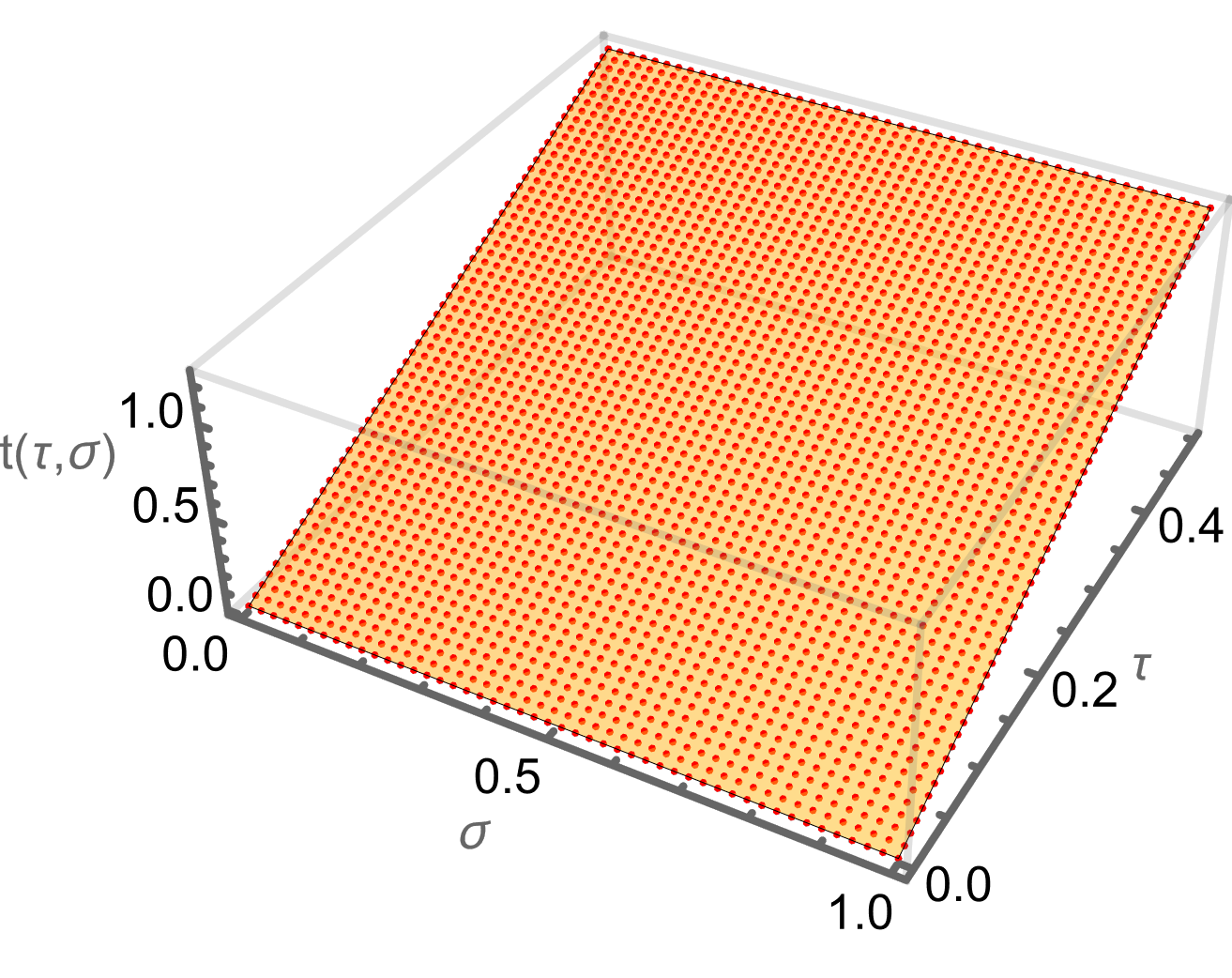}
    \caption{(top) Two  Classical solution for the field $\phi$ (top) and the time-mapping $t$ (bottom) from the critical point of \cref{eq:discrEsim} evaluated on a grid with $N_\sigma=48$ and $N_\tau=60$ points. Note that the initial and boundary conditions of \cref{fig:IBconfPhit} are exactly obeyed.}
    \label{fig:SolutionPhiAndt}
\end{figure}

We proceed by presenting the first result of our numerical study. The top panel of \cref{fig:SolutionPhiAndt} depicts the solution obtained for the field $\phi(\tau,\sigma)$ as a function of our abstract parameters $\tau$ and $\sigma$. The bottom panel on the other hand contains the values for the time-mapping $t(\tau,\sigma)$, which, due to our choice of $T=10^4$ at first sight resembles a plain sheet.

Our initial field configuration amounts to a  spatially symmetric bump of nearly Gaussian form. Based on intuition, we expect that this structure corresponds to an almost equal amount of left and right propagating modes. Indeed, we see two wave-packages emerge from the initial bump, travelling with opposite but equal speed towards the spatial boundaries of the simulation domain. 

The vanishing Dirichlet boundary conditions amount to the physical situation of the wave medium (such as e.g. an elastic sheet) being fastened at the boundary. As one would expect from Newton's third law, the amplitude of the wave-package changes sign after reflection. The wave attempts to push the boundary upward and the %(in this case infinitely heavy) 
boundary exerts an equal but opposite force onto the medium, leading instead to a downward movement. After reflection the wave-packages travel back towards the center of the spatial domain, where they interfere before heading towards the opposite end of the spatial domain. 

\begin{figure}
    \includegraphics[scale=0.4]{img/TmappingDerTauSolutionNs48Nt60tmax125SBP21.pdf}
    \includegraphics[scale=0.4]{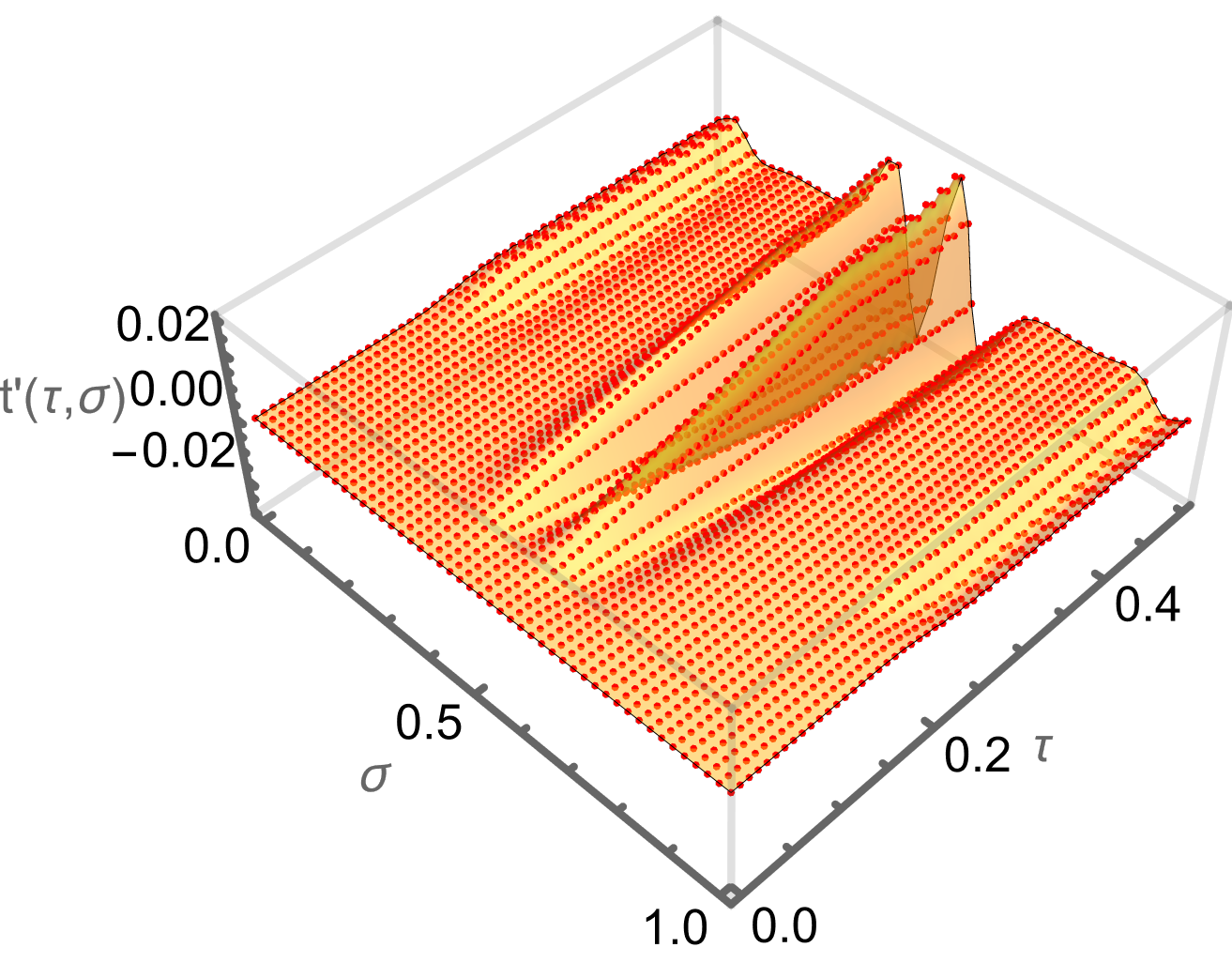}
    \caption{The temporal (top) and the spatial (bottom) derivatives of the time-mapping $t(\tau,\sigma)$ as obtained from the critical point of \cref{eq:discrEsim}, evaluated on a grid with $N_\sigma=48$ and $N_\tau=60$ points.}
    \label{fig:derivst}
\end{figure}

Next, we take a closer look at the behavior of the time mapping, whose $\tau$ derivative is plotted in the top panel of \cref{fig:derivst}, while its $\sigma$ derivative is shown in the bottom panel. We observe nontrivial structures intimately related to the dynamics of the propagating field.  
Starting with the $\tau$ derivative, if coordinate maps and fields evolved independently (as would be the case for $T\to\infty$), one would expect to find a constant value of $\Delta t/\Delta \tau=\dot t_{\rm IC}=5/2$ . This constant value amounts to a uniform temporal grid spacing throughout the simulation domain. For finite values of $T$ we observe deviations from such constant behavior, with larger values indicating a coarser and smaller values indicating a finer temporal grid in different regions of the domain spanned by $(\tau,\sigma)$. In particular two structures are noteworthy: close to the center of the $\sigma$ domain we find a double ridge structure stretching out over a finite interval in $\tau$ direction and then repeating as $\tau$ progresses. In addition we find significant deviations from constancy at the spatial boundaries of the simulation domain.

The extent of the central structure can be understood from a comparison with the propagation of the wave packages in the top panel of \cref{fig:SolutionPhiAndt}. As the wave-packages recede from their initial position in \cref{fig:SolutionPhiAndt}, no non-trivial dynamics occurs in the center of the simulation domain. Thus a coarser time spacing suffices in that region. Similarly, later, when the wave packages recombine to interfere, a field configuration very similarly to the initial conditions (except for a flipped sign) occurs and the temporal spacing reduces to resolve these non-trivial dynamics. With the wave-packages thereafter leaving the central region again to travel to the opposite end of the spatial domain there is also no more need for a high resolution and the ridge structure reappears. 

The structures close to the boundary also correlate with the propagation of the wave-packages. I.e. the deviations from the constant values occurs exactly where the wave-packages reflect from the spatial boundary. The temporal spacing is not just reduced but shows an intricate mix of refinement and coarsening. As the wave-package hits the wall, the adaption of the temporal grid emerges fully dynamically.
%We stress that the location where the wave-packages hit the wall is \textit{a priori} unknown to the algorithm and that the adaption of the temporal grid emerges fully dynamically here.

Continuing the investigation of the time mapping, we inspect its spatial derivative in the bottom panel of \cref{fig:SolutionPhiAndt} next. As expected from our discussion above, the central region and the boundaries show non-trivial modifications. Non-trivial gradients emerge in the central region, which appear accumulative. They build up as the wave-packages recede from their initial positions, while their growth is stopped over the short interval where the wave-packages recombine to interfere. Similarly, after the wave-packages have hit the spatial boundary, the time-mapping develops gradients, which persist at larger $\tau$ and we expect these effects to accumulate with each encounter of the boundary by the wave-packages.

In summary, the results clearly show that in our approach with dynamical coordinate maps the temporal grid adapts automatically to the dynamics of the propagating fields, constituting a genuine form of \textit{automatic adaptive mesh refinement}.

\begin{figure}
    \includegraphics[scale=0.27]{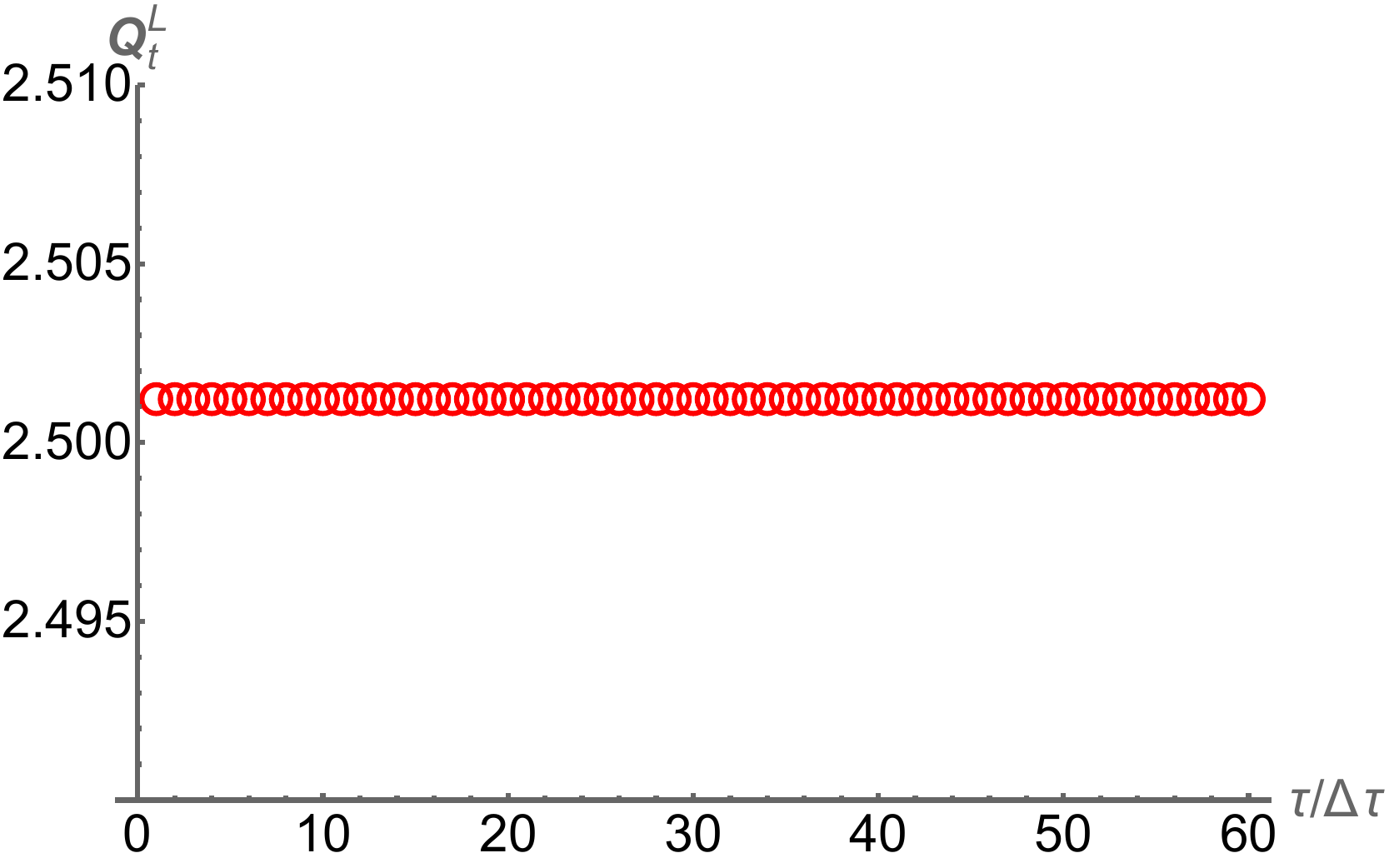}
    \includegraphics[scale=0.28]{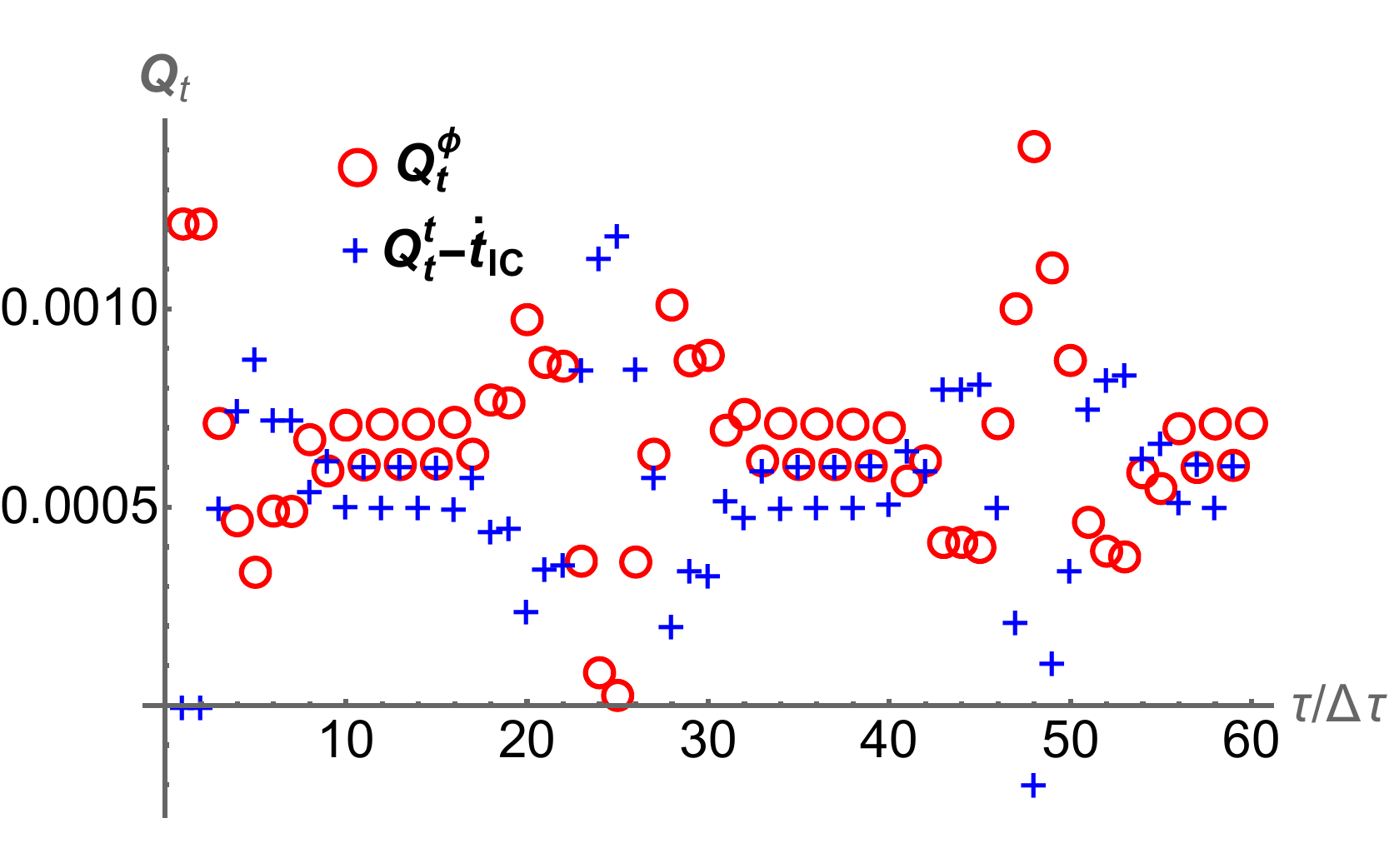}
    \caption{(top) Preservation of the Noether charge ${\bm Q}^{\rm L}_t$ (red open circles) over the whole simulation domain. (bottom) The two individual contributions to the Noether charge from the coordinates (blue crosses) and from the field dynamics (red circles). Note that we have subtracted the initial value $\dot t_0$ from the coordinate contribution for better visibility of the mutually complementary behavior. }
    \label{fig:NoetherChargeNt60}
\end{figure}

Let us continue with the second central result of our numerical study, the evaluation of the \textit{energy}, i.e. the Noether charge associated with time-translation symmetry in our approach. In the top panel of \cref{fig:NoetherChargeNt60} we show the values of ${\bm Q}^{\rm L}_t$, defined in \cref{eq:dsicrNoetherl}, which is associated with time-translations, a manifest symmetry of our discretized action \cref{eq:discrEsim}. Its values are given as red circles. A visual inspection reveals that ${\bm Q}^{\rm L}_t$ \textit{remains constant throughout the whole simulation with no deviation from its initial value}. Since ${\bm Q}^{\rm L}_t$ makes reference to first derivatives only, its value at initial $\tau^{\rm i}$ is determined solely by the initial conditions provided and that value persists throughout the simulation\footnote{While we can fix the value and $\tau$ derivatives of the degrees of freedom to their continuum values at $\tau^{\rm i}$, the values of the spatial derivatives follow self-consistently from an evaluation of $\mathds{D}_\sigma$ applied to the discretized initial values. I.e. when the initial energy is evaluated at different lattice spacings, the spatial derivative will encounter the same initial shape of e.g. the field bump at different resolutions, providing slightly different values for the energy there.}.  
%For completeness and to emphasize the important role that the explicit inclusion of boundary data via Lagrange multiplier plays, we also plot the values of the naive Noether charge ${\bm Q}_t$ as gray circles. That quantity receives correction terms due to the connection conditions enforced by $\tilde{\bm \gamma}^t$ at the final $\tau^{\rm f}$. Only when taking these contributions into account, as correctly done for ${\bm Q}^{\rm L}_t$, do we obtain a truly conserved Noether charge.

\begin{figure}
    \hspace{0.7cm}\includegraphics[scale=0.27]{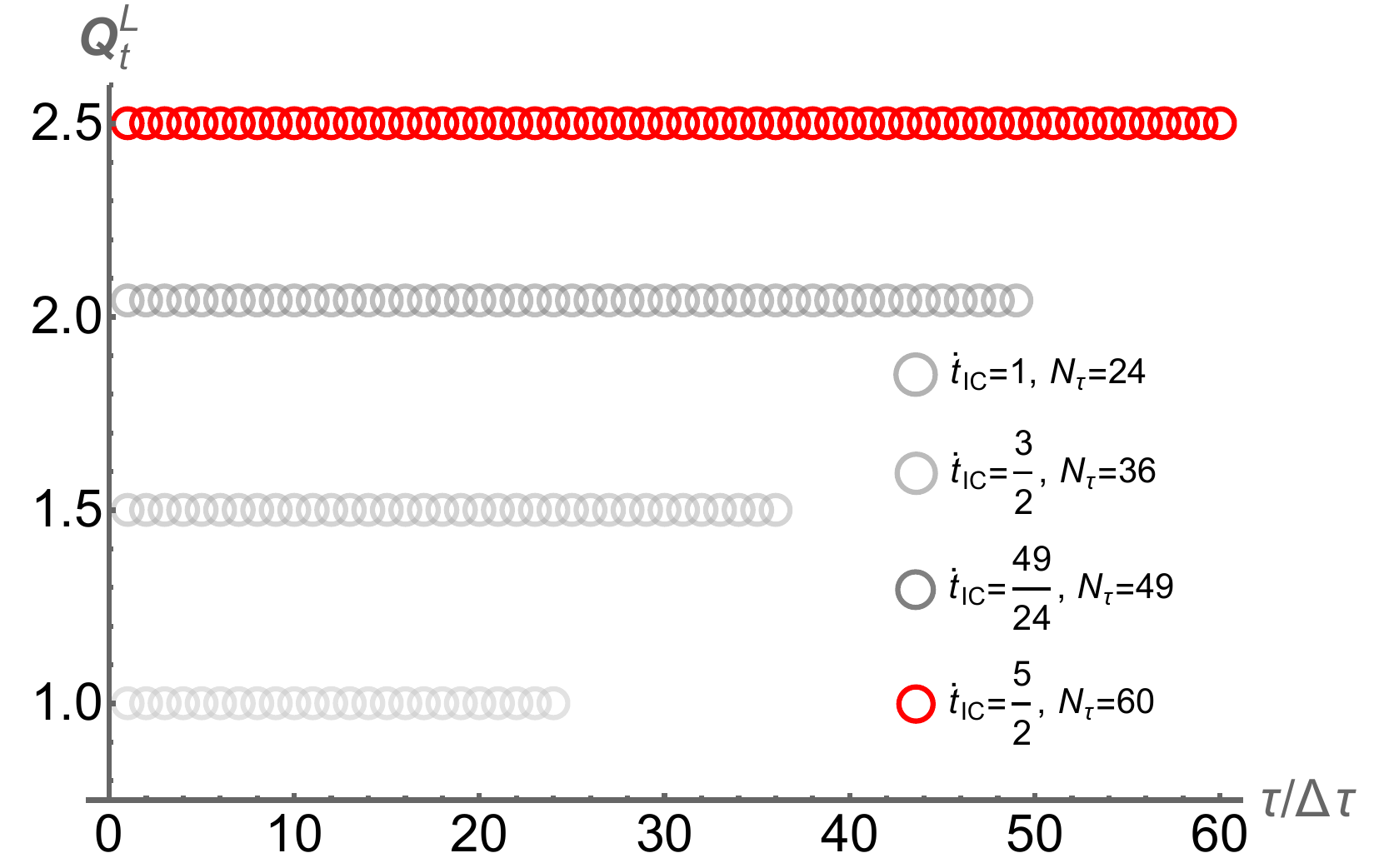}
    \caption{Values of the Noether charge ${\bm Q}^{\rm L}_t$ for different temporal extent of the simulation, implemented via different values of the initial $\tau$ derivative $\td$ of the time-mapping. Note that magnitude of the Noether charge depends on the initial value of $\dot t$.}
    \label{fig:DiffNoether}
\end{figure}

The total energy of our system is composed of a contribution from the field, and in addition from the time mapping, according to \cref{eq:1p1dnoethcont}. In the bottom panel of \cref{fig:NoetherChargeNt60} we plot the two contributions separately. To highlight the mutually complementary behavior we have subtracted the initial value $\dot t_{\rm IC}$ from the coordinate contribution. One finds that while the overall energy ${\bf Q}_t^{\rm L}$ is identically preserved, the coordinate and field contribution each show minute variation in time, which exactly cancel in their sum. This dependence of the energy on the coordinate mapping was expected, as it occurred in a similar fashion in the world-line formalism.

With its explicit coordinate mapping dependence, we showcase the values of ${\bm Q}^{\rm L}_t$ for different $\dot t_{\rm IC}$ in \cref{fig:DiffNoether}. The underlying resolution of the $\tau$ direction is kept constant and we increase $N_\tau$ such that approximately the same $t$ resolution is obtained. Due to our choice of a large value of $T=10^4$, the main contribution to the energy indeed comes from the time mapping itself, which is why the values of ${\bm Q}^{\rm L}_t$ increase with increasing $\td_{\rm IC}$. Note that the values are slightly larger than $\td$ and it is this difference that encodes the contributions from the field dynamics. We have explicitly checked that for all initial values $\dot t_{\rm IC}$ shown in \cref{fig:DiffNoether} the Noether charge remains at its initial value within the global numerical accuracy of our solution $\Delta \sim 10^{-6}$.

%The visual inspection of ${\bm Q}^{\rm L}_t$ in the top panel of \cref{fig:DiffNoether} indicated that its values remain constant. In the bottom panel of \cref{fig:DiffNoether} we provide quantitative confirmation, by plotting $\Delta E$ the difference between the discrete Noether charge at some $\tau$ vs. its value at $\tau^{\rm i}$, where it is specified solely by the initial conditions. We find values that are of the order of $10^{-7}$. As the global numerical accuracy of our solution obtained by minimization with the double precision arithmetic \texttt{IPOPT} library was $\Delta \sim 10^{-6}$, all the points shown for $\Delta E$ are consistent with zero.

Having established that the Noether charges are exactly preserved, we can now propose an explanation as to the guiding principle underlying the automatic adaptive mesh refinement observed in our simulations. We argue that it is \textit{based on the preservation of the space-time symmetries of our approach and the conservation of the associated Noether charge}. Only if the space-time grid adapts to the dynamics of the field, is it possible for a non-trivial Noether charge, such as ${\bm Q}_t^{\rm L}$ to be exactly preserved.

Note that while exact conservation of continuum Noether charges is a central feature of our novel approach, it is not synonymous with obtaining the continuum solution. For a finite grid spacing the solution will be affected by numerical error, which, as we demonstrate below, diminishes under grid refinement. The numerical errors however, by construction of our method, \textit{cannot affect the value of the Noether charges}.

Let us briefly return to the fact that the explicit inclusion of initial, boundary and connecting conditions in our discretized action \cref{eq:discrEsim} introduces non-trivial contributions to the values of the field and coordinate mappings at the spatial and temporal boundaries. 
%These contributions appear in the definition of the Noether charge in \cref{eq:dsicrNoetherl} and as shown in \cref{fig:NoetherChargeNt60} are crucial to take into account, in order obtain the correct expression for the conserved Noether charge ${\bf Q}^{\rm L}_t$ (compared to the naive expression ${\bf Q}_t$). Similarly 
We find e.g. that the $\tau$ derivative of the time mapping $\dot t$ in the top panel of \cref{fig:derivst} shows a jump at the last $\tau$ slice. This behavior is connected to the influence of the Lagrange multipliers there. The jump corresponds to a finite but small correction to the discretized equation of motion derived in the absence of Lagrange multipliers. Just as for the Noether charge it should be possible to define an expression for the time mapping in which the effect of the Lagrange multiplier is accounted for and for which the jump in the last time slice is absent. %Since the effects of the Lagrange multipliers in the Noether charge have been successfully captured in \cref{eq:dsicrNoetherl}.
We did not pursue the derivation of this correction term of $t$ further at this point. However, as we will show next, that the effects of the Lagrange multipliers do not affect the convergence of the solution under grid refinement.

\begin{figure}
    \includegraphics[scale=0.32]{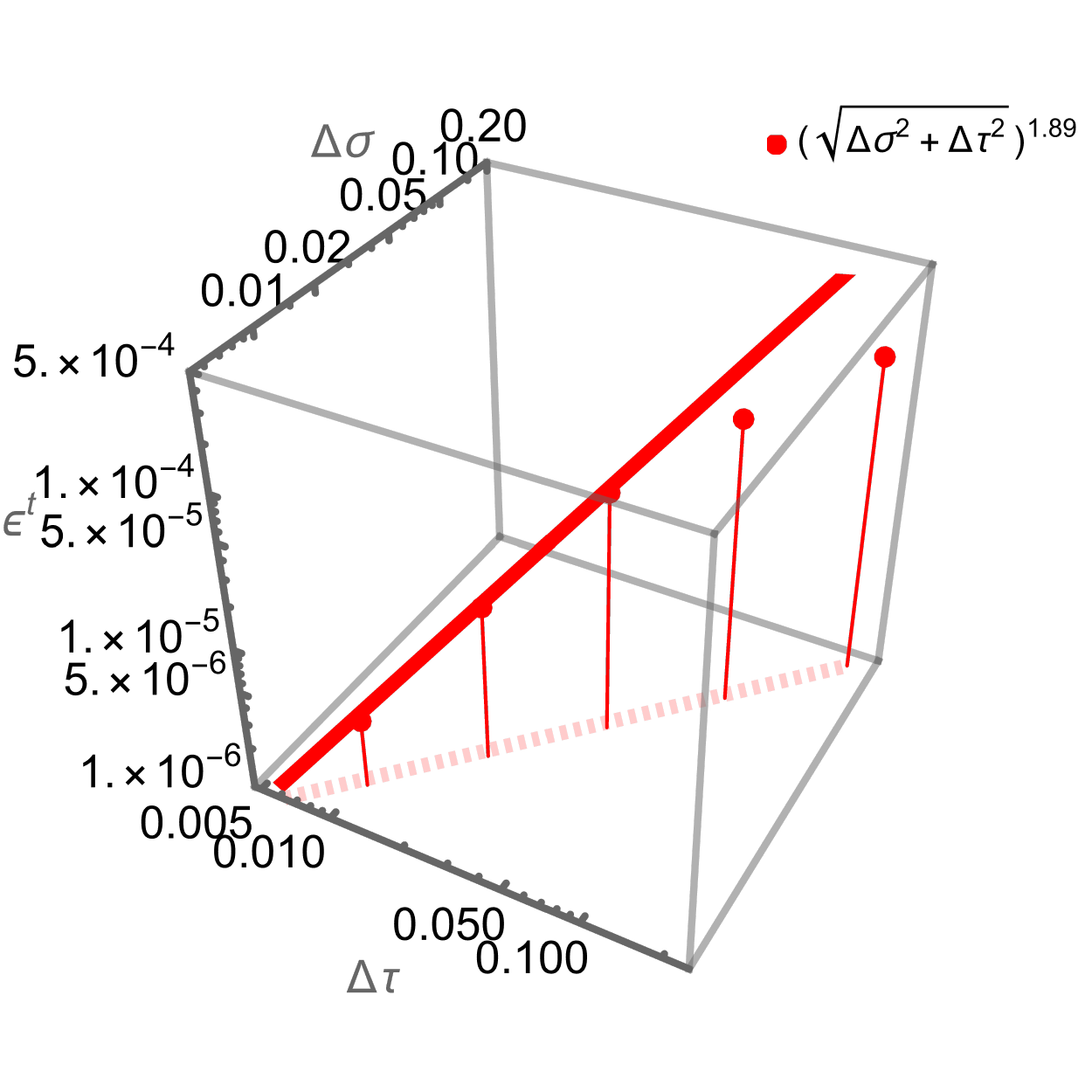}
    \includegraphics[scale=0.32]{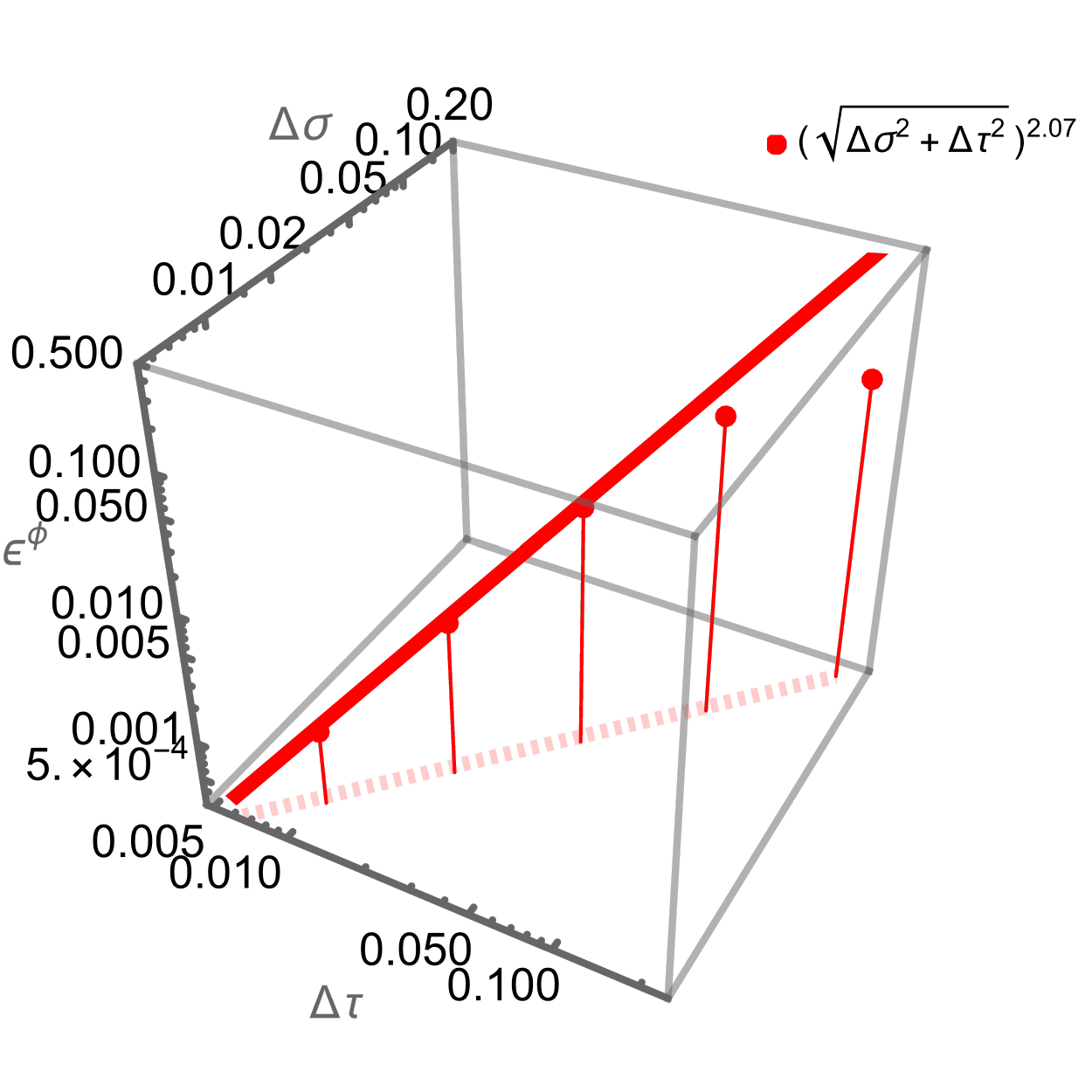}
    \caption{The global $L_2$-norm difference of the solutions for $t$ (top) and field $\phi$ (bottom) obtained via the direct determination of the critical point of $\mathds{E}^{\rm L}_{\rm IBVP}$ and by solving \cref{eq:1p1deomt,eq:1p1deomphi} on high resolution grids using the \texttt{MethodOfLines} method of \texttt{Mathematica}'s \texttt{NDSolve}.}
    \label{fig:Convergence}
\end{figure}

To demonstrate that we obtain the correct convergence rate under grid refinement, we compare the time mapping $t$ and field $\phi$ obtained from our approach, i.e. from finding the critical point of $\mathds{E}_{\rm IBVP}^{\rm L}$ to solutions of the continuum equations of motion \cref{eq:1p1deomt,eq:1p1deomphi} obtained via the \texttt{MethodOfLines} method of \texttt{Mathematica}'s \texttt{NDSolve} on high resolution grids, corresponding to the default \texttt{AccuracyGoal} of eight significant digits. The global $L_2$-norm differences obtained via integration over both $\tau$ and $\sigma$ 
\begin{align}
    &\epsilon^t = \sqrt{ ({\bm t}_1-{\bm t}_{\rm NDSolve})^T \mathds{H}\, ({\bm t}_1-{\bm t}_{\rm NDSolve})},\\
    &\epsilon^\phi = \sqrt{ ({\bm \phi}_1 - {\bm \phi}_{\rm NDSolve}) ^T \mathds{H} \,({\bm \phi}_1 - {\bm \phi}_{\rm NDSolve}) },
\end{align}
are shown in \cref{fig:Convergence}. The top panel shows $\epsilon^t$, the bottom panel $\epsilon^\phi$, each as red symbols. We concurrently reduce the spatial and temporal grid spacing and find that the deviation systematically reduces. The solid red line denotes the best power-law fit $\alpha \sqrt{\Delta \sigma^2+\Delta \tau^2}^\beta$ to the results. The fits tell us that for $\epsilon^t$ we have $\beta=1.89$ and for $\epsilon^\phi$ the scaling goes as $\beta=2.07$. Both of these values lie close to the ideal $\beta=2$, which corresponds to the ideal scaling of a \texttt{SBP121} operator. We note that due to the use of the pre-compiled \texttt{IPOPT} library and the associated limitation to global numerical accuracy to around $\Delta\sim10^{-6}$ we did not grid refine further.

\begin{figure}
    \includegraphics[scale=0.35]{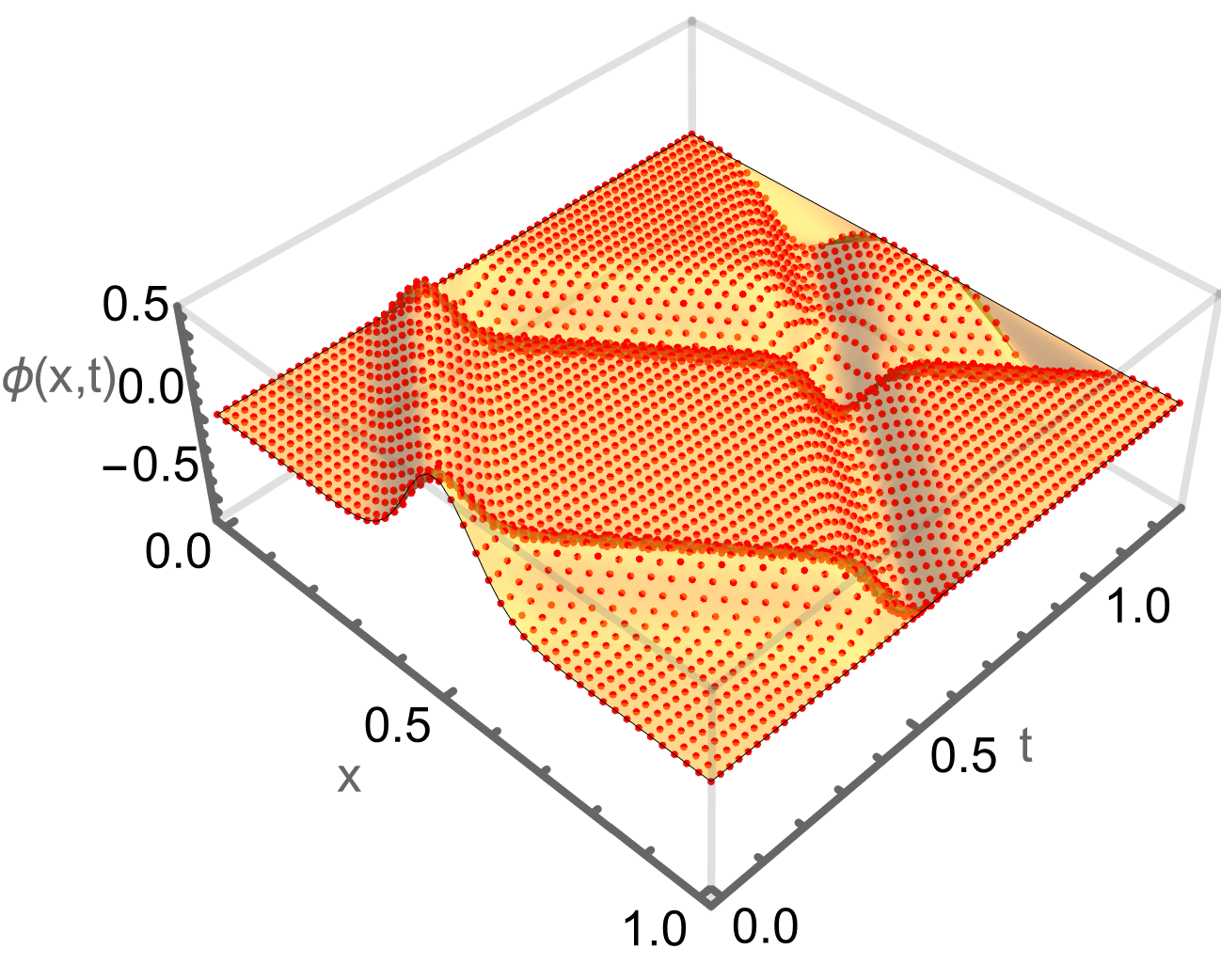}
    \includegraphics[scale=0.32]{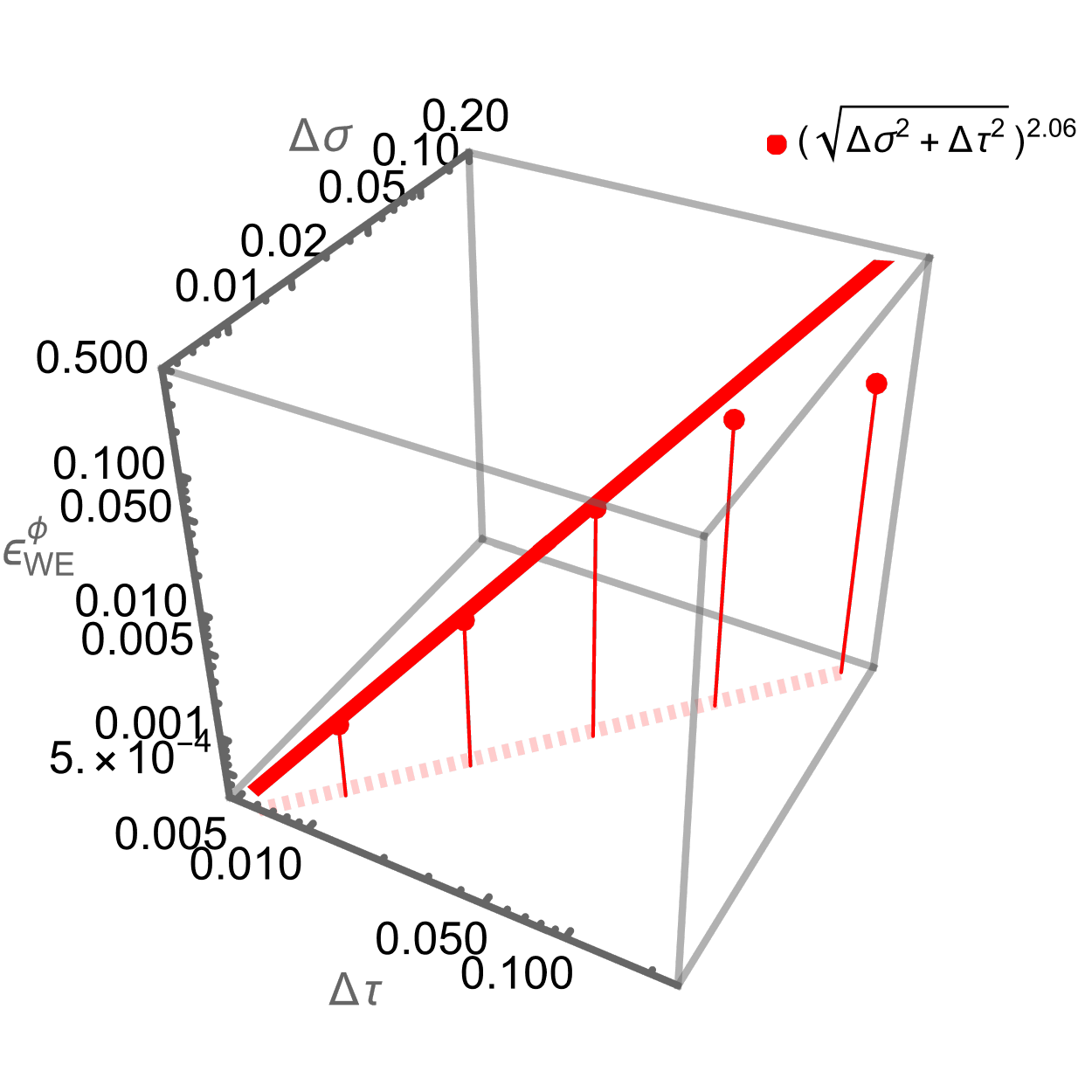}
    \caption{(top) The field $\phi(t,x)$ as function of the physical space-time coordinates $x$ and $t$. (bottom) The global $L_2$-norm difference between the solution for $\phi(x,t)$ obtained via the direct determination of the critical point of $\mathds{E}^{\rm L}_{\rm IBVP}$ and by solving the conventional wave equation on high resolution grids using the \texttt{MethodOfLines} method of \texttt{Mathematica}'s \texttt{NDSolve}. }
    \label{fig:physicalsol}
\end{figure}

\begin{figure}
    \includegraphics[scale=0.27]{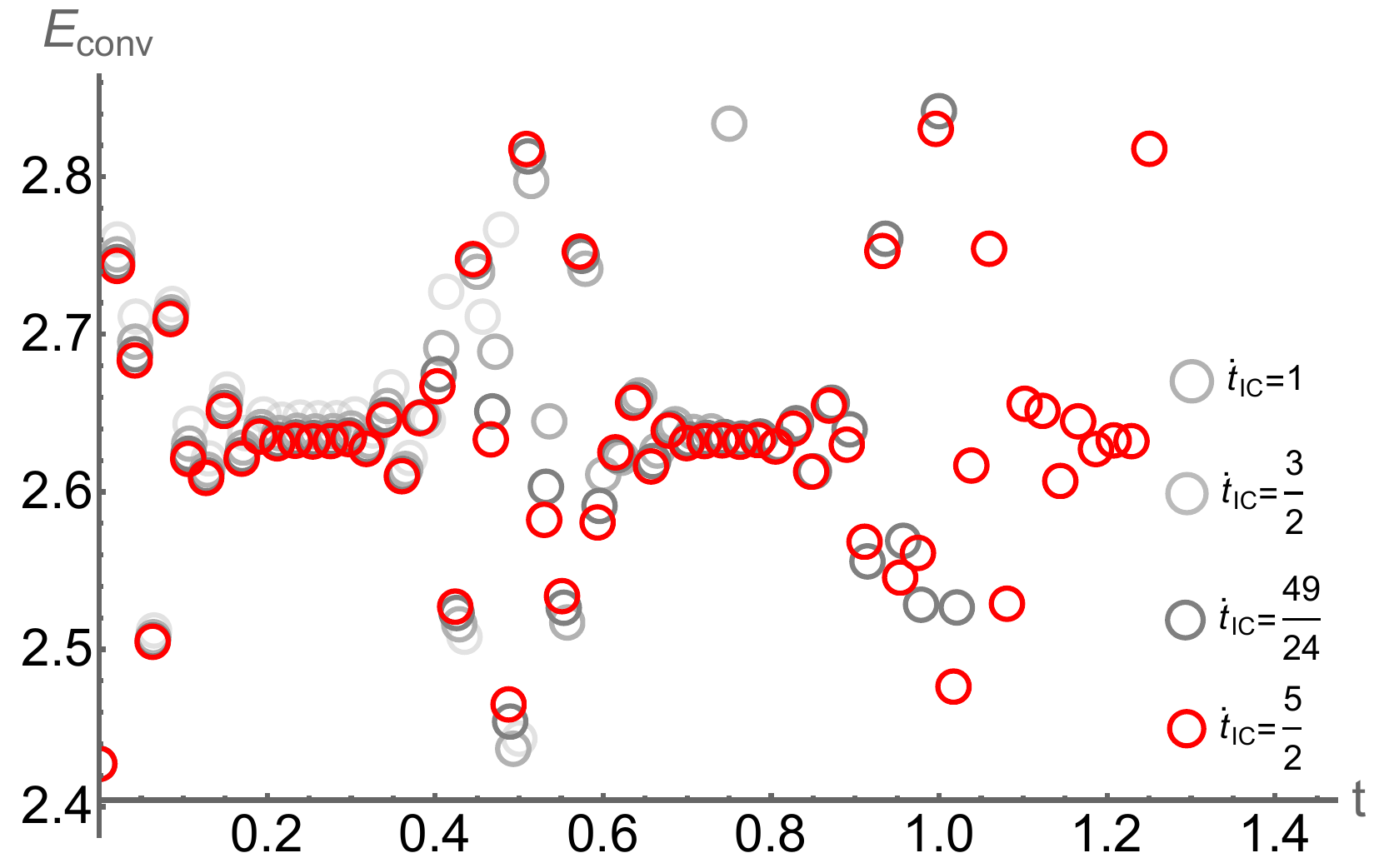}
    \caption{The conventional energy $E_{\rm conv}$ defined in \cref{eq:convnrg}, as a function of time $t$ on solutions of the field that propagate to different final times $t^{\rm f}$. Note that $E_{\rm conv}$ is insensitive to the initial $\tau$ derivative of the time mapping $\dot t$, which we use to select the final time of our simulation.}
    \label{fig:physicalsolnrg}
\end{figure}

For completeness let us also visualize the field solution in physical space-time coordinates, as shown in the top panel of \cref{fig:physicalsol}. While the spatial domain is fixed at all times to have the same extent, due to the use of our trivial mapping $x=\sigma$, the temporal simulation domain extends to slightly different values depending on the location along $x$. Due to T being large, these effects here are relatively small, as we saw that the deviation from a constant value of the derivatives of $t(\tau,\sigma)$ were at most on the level of a percent. 

We may ask how well the physical field configuration $\phi(x,t)$ obtained from our new approach agrees with the solution of the conventional wave equation (WE) in (\ref{eq:stdWPgoveq}). To this end we solve the wave equation via the \texttt{MethodsOfLines} method of \texttt{NDSolve} in \texttt{Mathematica} on high resolution grids, to achieve the \texttt{AccuracyGoal} of eight significant digits. This numerical solution is denoted as ${\bm \phi}_{\rm WE}$.

Due to the presence of the square root, the novel action contains additional terms, which are however suppressed by powers of $T$. For our choice of $T=10^4$ we find that for all the grid spacings considered in this study, convergence to the numerical wave equation solution ${\bf \phi}_{\rm WE}$ in the global $L_2$ norm
\begin{align}
    &\epsilon^\phi_{\rm WE} = \sqrt{ ({\bm \phi}_1 - {\bm \phi}_{\rm WE}) ^T \mathds{H} \,({\bm \phi}_1 - {\bm \phi}_{\rm WE}) }
\end{align}
can be found under grid refinement, as shown in the bottom panel of \cref{fig:physicalsol}. Of course, as one increases the grid resolution further one would eventually spot differences between our solution and that of the wave equation when the numerical error becomes smaller than the first correction term involving $T$.

It is also instructive to evaluate the conventional energy in the field $\phi$
\begin{align}
\nonumber E_{\rm conv}&=\int dx \frac{1}{2}\Big\{ \Big( \frac{\partial \phi}{\partial t}\Big)^2 + \Big( \frac{\partial \phi}{\partial x}\Big)^2 \Big\}\\
&\approx \frac{1}{2}\mathds{H}_\sigma \Big\{ \Big( {\mathds{D}_\tau}{\bm \phi}/  {\mathds{D}_\tau}{\bm t} \Big)^2 + \Big( {\mathds{D}_\sigma}{\bm \phi}/  {\mathds{D}_\sigma}{\bm x} \Big)^2\Big\}\label{eq:convnrg}
\end{align}
where division is understood in a point-wise fashion. The values of \cref{eq:convnrg} for different physical time extent are shown in \cref{fig:physicalsolnrg}. Note first, that in contrast to ${\bf Q}_t^{\rm L}$ this quantity at $t=0$ does not vary with $\dot t_{\rm IC}$, since the conventional energy only knows of the field itself. In addition we see that while the values of the conventional energy seem to stay close to a common value on average, the individual values of $E_{\rm conv}$ show significant oscillatory patterns (and a jump from the value at $t=0$). It is clearly not the actual conserved Noether charge associated with time translation symmetry.

\begin{figure}
    \includegraphics[scale=0.27]{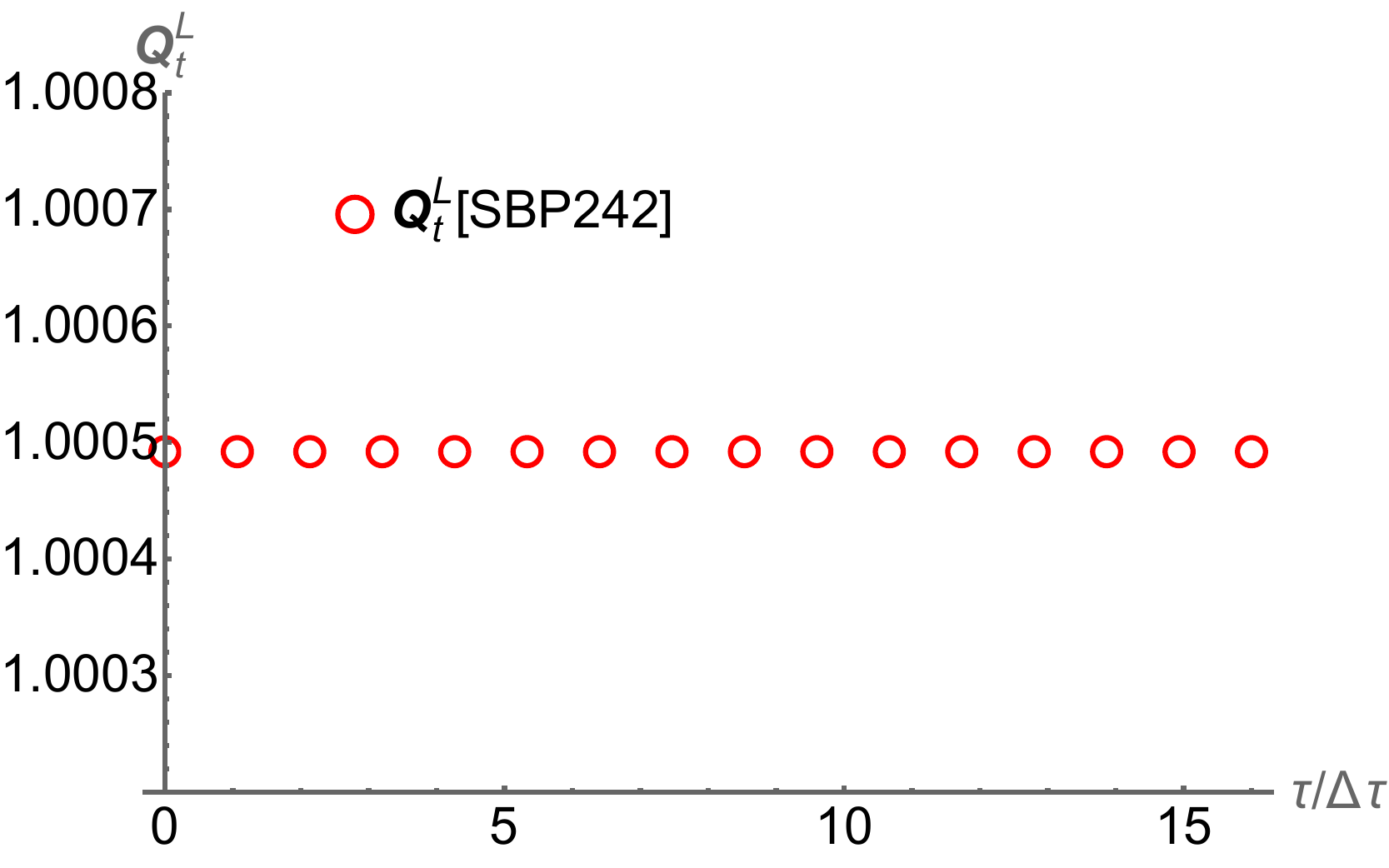}
    \caption{The values of the Langrange multiplier corrected ${\bm Q}^{\rm L}_t$ Noether charge from the \texttt{SBP242} operator on a grid with $N_\tau=16$ and $N_\sigma=24$. Also for the higher order operator we confirm that within our accuracy the Noether charge is exactly preserved.}
    \label{fig:Noether42}
\end{figure}

While a thorough investigation of the scaling properties of our approach for \texttt{SBP242} operators is planned for a future study, let us briefly confirm that the conservation of the Noether charge indeed holds even for the higher order SBP discretization scheme. I.e. the exact conservation observed in \cref{fig:NoetherChargeNt60} is not an accident of the lowest order with its very simple boundary modifications within $\mathds{D}$. Hence, in \cref{fig:Noether42} we show the values of ${\bm Q}^{\rm L}_t$ on a $N_\sigma=24$ and $N_\tau=16$ grid obtained from the \texttt{SBP242} scheme as red circles.  Note that also for \texttt{SBP242}, as long as the Lagrange multiplier contributions to the Noether charge are correctly incorporated according to \cref{eq:dsicrNoetherl}, we obtain \textit{exactly conserved values of the Noether charge}. 

\section{Summary and Outlook}

We have presented \textit{two novel conceptual contributions} to address central challenges of formulating IBVPs.

%a novel approach to the numerical solution of second order initial boundary value problems (IBVPs) introducing \textit{two novel conceptual contributions} to address the three central challenges of formulating IBVPs.

\subsection{Summary}

The \textit{first novel contribution} is the derivation of a reparameterization invariant continuum action, in which coordinate maps are incorporated as dynamical degrees of freedom. A new energy density scale $T$ determines when the dynamics of the field become inextricably intertwined with that of the coordinate maps.

In the novel action we do not integrate over space-time coordinates but instead over a set of underlying abstract parameters $\Sigma^a$. This fact leads to the \textit{second novel contribution}, the ability to discretize the action in terms of the abstract parameters, leaving the values of the coordinate maps continuous. In turn the discretized action retains manifest invariance under infinitesimal Poincar\'e transformations and Noether's theorem remains intact. Using summation-by-parts (SBP) operators, the Noether charges are exactly preserved in time.

To formulate a causal initial boundary value problem we deployed the doubled degrees-of-freedom construction of Galley, to obtain the IBVP action. The requires initial, boundary and connecting conditions are made explicit by introducing appropriate sets of Lagrange multiplier functions in the continuum.

Deploying appropriately regularized summation-by-parts finite difference operators formulated in affine coordinates, we outlined how to discretize the novel system action and arrive at the correct expression for the discrete Noether charge.

\Cref{sec:proofofprinciple11d} presented a proof-of-principle based on wave propagation in $1+1$ dimensions deploying a fully dynamic time mapping $t(\tau,\sigma)$ and a propagating field $\phi(\tau,\sigma)$. We showed how our approach offers new \textit{flexibility in boundary treatment}.

\noindent Our main numerical results are as follows:
\begin{itemize}
    \item The coordinate map $t$ adapts dynamically to the evolution of the field $\phi$, exhibiting region of coarser and finer resolution, correlated with the field dynamics. This constitutes a form of \textit{automatic adaptive mesh refinement}.
    \item The \textit{Noether charge} of time translation symmetry, is \textit{exactly conserved at its initial value}. This energy receives contributions from the field and the coordinate map. We find that the conventional energy density on the other hand shows oscillatory patterns, as it is not the correct conserved charge.
    \item The novel approach implemented with \texttt{SBP121} operators \textit{converges correctly under grid refinement}.
    \item  We have confirmed that the Noether charge is exactly conserved also for a \textit{higher order} \texttt{SBP242} scheme.
\end{itemize}

\subsection{Future work}

Our proof-of-principle explored the dynamics of the time-mapping $t(\tau,\sigma)$ but kept the spatial mapping trivial as $x=\sigma$. The inclusion of a fully dynamic spatial map $x(\tau,\sigma)$ is work in progress. While formally straight forward, one has to take into account that for fully dynamic coordinate maps the action contains additional redundancies. In the string-theory literature these are addressed by selecting e.g. conformal gauge.

As more flexible coordinate maps are introduced also the freedom in the boundary terms increases with more degrees of freedom contributing to their values. We are working on exploring whether there exist compatibility requirements for the boundary conditions, especially in the presence of a specific gauge choice.

Our proof of principle was implemented in $1+1$ dimensions and an extension to higher dimensions is called for. The numerical minimisation underlying the search for the critical point $\mathds{E}^{\rm L}_{\rm IBVP}$ will become more costly. This will eventually demand a departure from the \texttt{Mathematica} software and implementation in a more performant language, such as e.g. \texttt{C++} via the \texttt{Ceres} library, which is work in progress. 

In this study we have enforced initial, boundary, and connecting conditions in a strong manner with Lagrange multipliers. Weak imposition of boundary data offers many advantages, especially for numerical stability. It is therefore of interest to formulate a genuine weak version of our action based approach.

We believe that the preservation of the Noether charge is closely related to stability. In future work we will elaborate on this and will attempt a proof.

As explained but not elaborated on, the presence of $\tau$ and $\sigma$ derivatives on the coordinate fields opens up for new constructions of non-reflecting boundary conditions. We will explore that possibility for the wave equation in future work.

\subsection{Speculations}

It is natural to ask, whether a similar reparameterization invariant action as in the scalar case, can be constructed for other relevant theories, such as gauge fields (linear Maxwell electromagnetism, non-linear Yang-Mills theory). Treatment of gauge fields in string theory via the Born-Infeld construction suggests that such a construction is not unique but may be equivalent to a low order in the non-relativistic expansion.

So far the benefits of our novel action have been elucidated in the context of classical field theory. The preservation of space-time symmetries plays a central, if not more crucial role, for the accurate simulation of quantum fields. As breaking of space-time symmetries directly translates into a contamination of the particle spectrum observed in these simulations with unphysical modes, exact space-time symmetry preservation promises significant improvements in that context.

The novel action we propose motivates several lines of inquiry on the physics side. One may ask what the role of the new scale $T$ is, which was central to the construction of $\cal S$. In the world-line formalism, its counterpart $mc$ is associated with the rest energy of a point particle, a fundamental property of the propagating degree of freedom. Here $T$ seems to lend itself to the interpretation of a reference energy density for the coordinate dynamics and it is intriguing to investigate its possible relation to the physics of gravity and the cosmological constant through the general theory of relativity and string theory.

We may further ask about possible physics implications of the square root term in the fully reparameterization invariant action with dynamical coordinate maps. If $\cal S$ indeed were the correct high-energy action for second order field theory, then at low energies, small but finite correction terms to the conventional $S_{\rm ft}$ will ensue. This hypothesis should be straight forward to falsify, by exploring whether such corrections are compatible with the wealth of high precision collider based measurements available for the fundamental interactions on short scales and via astronomical observations on large scales.

%We believe that our approach to discrete IBVPs with its focus on the preservation of continuum space-time symmetries and the associated Noether charges provides the community with a valuable novel perspective, complementary to the extensive research activity on the subject. The demonstrated ability of our approach to provide an automatic adaptive mesh refinement guided by the symmetries of the system and the additional freedom in treatment of boundaries promises fertile ground for progress on front-line subjects related to industry relevant IBVPs. At the same time if the approach can be successfully transferred to the quantum domain it there too promises to significantly improve the accuracy of numerical simulations.

\section*{Acknowledgements}
A.\ R.\ and W.\ A.\ H.\ acknowledge support by the ERASMUS+ project 2023-1-NO01-KA171-HED-000132068. W.\ A.\ H.\ thanks the South African National Research Foundation and the SA-CERN Collaboration for financial support. J.\ N.\ was supported by the Swedish Research Council grant nr. 2021-05484 and the University of Johannesburg.

%%%%%%%%%%%%%%%%%%%%%%%%%%%%%%%%%%%%%%%%%
%%                                     									       %%
%%  So called Backmatter part starts here with acknowledgements		       %%
%%  funding information, and the bibliography.						       %%
%%                                     									       %%
%%%%%%%%%%%%%%%%%%%%%%%%%%%%%%%%%%%%%%%%%

\appendix    

\section{Derivation of the connecting conditions}
\label{sec:appconcond}

In the following we will derive how the action \cref{eq:novelactionEIVP} and connecting conditions of \cref{eq:SKconnecting} prevent the appearance of non-causal temporal boundary terms. To this end it is convenient to express ${\cal E}_{\rm IBVP}$ in linear combinations of $X_1,X_2$ and $\phi_1,\phi_2$. We consider the difference $X_-=X_1-X_2$, $\phi_-=\phi_1-\phi_2$ and the mean $X_+=\frac{1}{2}(X_1+X_2)$, $\phi_+=\frac{1}{2}(\phi_1+\phi_2)$ (known in the context of the Schwinger-Keldysh literature also as the quantum and classical components respectively). Let us inspect the variation of the action w.r.t. these new variables
\begin{align}
\delta {\cal E}_{\rm IBVP}[X_\pm,& \partial_a X_\pm, \phi_\pm, \partial_a \phi_\pm]\\
\nonumber=\int d^{(d+1)}\Sigma \Big\{ 
&\frac{\partial E_{\rm IBVP}}{\partial X_+^\mu}\delta X_+^\mu + \frac{\partial E_{\rm IBVP}}{\partial (\partial_a X_+^\mu)}\delta (\partial_a X_+^\mu) + \frac{\partial E_{\rm IBVP}}{\partial \phi_+}\delta \phi_+ + \frac{\partial E_{\rm IBVP}}{\partial (\partial_a \phi_+)}\delta (\partial_a \phi_+) \\
+& \frac{\partial E_{\rm IBVP}}{\partial X_-^\mu}\delta X_-^\mu + \frac{\partial E_{\rm IBVP}}{\partial (\partial_a X_-^\mu)}\delta (\partial_a X_-^\mu) + \frac{\partial E_{\rm IBVP}}{\partial \phi_-}\delta \phi_- + \frac{\partial E_{\rm IBVP}}{\partial (\partial_a \phi_-)}\delta (\partial_a \phi_-)\Big\}\label{eq:varEIVP}\\
\nonumber=\int d^{(d+1)}\Sigma \Big\{& 
 \Big( \frac{\partial E_{\rm IBVP}}{\partial X_+^\mu} -\partial_a \frac{\partial E_{\rm IBVP}}{\partial (\partial_a X_+^\mu)}\Big)\delta X_+^\mu + \Big(  \frac{\partial E_{\rm IBVP}}{\partial \phi_+} - \partial_a \frac{\partial E_{\rm IBVP}}{\partial (\partial_a \phi_+)} \Big)\delta \phi_+ \\
 +& 
 \Big( \frac{\partial E_{\rm IBVP}}{\partial X_-^\mu} -\partial_a \frac{\partial E_{\rm IBVP}}{\partial (\partial_a X_-^\mu)}\Big)\delta X_-^\mu + \Big(  \frac{\partial E_{\rm IBVP}}{\partial \phi_-} - \partial_a \frac{\partial E_{\rm IBVP}}{\partial (\partial_a \phi_-)} \Big)\delta \phi_-\Big\} \label{eq:varEeomIVP}\\
+\int &\prod_{a=1}^{d}d\Sigma_a\Big\{ \Big[ \frac{\partial E_{\rm IBVP}}{\partial (\partial_0 X_+^\mu)} \delta X_+^\mu \Big]_{\tau^{\rm i}}^{\tau^{\rm f}} + \Big[ \frac{\partial E_{\rm IBVP}}{\partial (\partial_0 X_-^\mu)} \delta X_-^\mu \Big]_{\tau^{\rm i}}^{\tau^{\rm f}} \label{eq:tauBCEvarXIVP} \\
&\qquad \quad   + \Big[ \frac{\partial E_{\rm IBVP}}{\partial (\partial_0 \phi_+)} \delta \phi_+ \Big]_{\tau^{\rm i}}^{\tau^{\rm f}} + \Big[ \frac{\partial E_{\rm IBVP}}{\partial (\partial_0 \phi_+)} \delta \phi_+ \Big]_{\tau^{\rm i}}^{\tau^{\rm f}} \Big\}\label{eq:tauBCEvarphiIVP}\\
\nonumber + {\rm spatial\, B.C.}
\end{align}
For the discussion of causality, the spatial boundary conditions do not play a role, which is why we do not list them here explicitly. Prescribing fixed initial conditions for the values of the d.o.f. in our system, we have $\delta X^\mu_\pm(\tau^{\rm i},\vec{\sigma})=0$ and  $\delta \phi_\pm(\tau^{\rm i},\vec{\sigma})=0$, making the lower boundaries of the terms in (\ref{eq:tauBCEvarXIVP}) and (\ref{eq:tauBCEvarphiIVP}) vanish.

At the end of the forward branch at $\tau^{\rm f}$ we can achieve the vanishing of some of the boundary terms by requiring that the values of the degrees of freedom on the forward and backward branch agree
\begin{align}
X_1^\mu(\tau=\tau^{\rm f},\vec{\sigma})=X_2^\mu(\tau=\tau^{\rm f},\vec{\sigma}),\quad \phi_1(\tau=\tau^{\rm f},\vec{\sigma})=\phi_2(\tau=\tau^{\rm f},\vec{\sigma}).\label{eq:connval}
\end{align}
This leads to $\delta X^\mu_-(\tau^{\rm f},\vec{\sigma})=0$ and  $\delta \phi_-(\tau^{\rm f},\vec{\sigma})=0$. The only remaining terms are those involving $X_+^\mu$ and $\phi_+$ at the final $\tau^{\rm f}$.

Since for a causal IBVP the values of the d.o.f. are \textit{not a priori} known at the final time $\tau^{\rm f}$, the corresponding variation in (\ref{eq:tauBCEvarXIVP}) and (\ref{eq:tauBCEvarphiIVP}) cannot be made to vanish there. Let us instead inspect the terms which multiply the variations, which in this case are $\partial E_{\rm IBVP}/\partial (\partial_0 X_+^\mu)$ and $\partial E_{\rm IBVP}/\partial (\partial_0 \phi)$. Since the variations for each $\mu$ component of $X$ are independent, we can consider just one of them, e.g. $\mu=d$, to establish how the boundary contribution for each and all of them can be made to vanish. 
By exploiting the inverse relations $X^\mu_1=X^\mu_+ + \frac{1}{2} X^\mu_-$ and $X^\mu_2=X^\mu_+ - \frac{1}{2} X^\mu_-$ we find that
\begin{align}
    \nonumber\frac{\partial E_{\rm IBVP}}{\partial (\partial_0 X_+^d)} &=  \frac{\partial E_{\rm IBVP}}{\partial (\partial_0 X_1^d)} \frac{\partial_0 X_1^d}{\partial_0 X_+^d}+\frac{\partial E_{\rm IBVP}}{\partial (\partial_0 X_2^d)} \frac{\partial_0 X_2^d}{\partial_0 X_+^d}\\
    &=  \frac{\partial E_{\rm BVP}[X_1,\partial_a X_1,\phi_1,\partial_a\phi_1]}{\partial (\partial_0 X_1^d)} -\frac{\partial E_{\rm BVP}[X_2,\partial_a X_2,\phi_2,\partial_a\phi_2]}{\partial (\partial_0 X_2^d)}\label{eq:tbnd0}
\end{align}
As one can confirm in \cref{eq:novelactionE}, $E_{\rm BVP}$ makes reference to $X$ in the terms involving ${\rm det}[g]$ and ${\rm adj}[g]_{ab}$. 

It is important to realize that in \cref{eq:tbnd0} each of the two terms refers to a different induced metric for the forward and backward branch respectively. Thus the full adjugate of the induced metric needs to be the same between the two branches for cancellation to occur. This in turn requires that at final $\tau^{\rm f}$ all the different derivatives of $X^\mu$ must agree on the two branches.
We thus conclude that in order for the boundary terms related to the $X_+$ variables to vanish, we must require 
\begin{align}
\partial_a X_1^\mu|_{\tau=\tau^{\rm f}}= \partial_a X_2^\mu|_{\tau=\tau^{\rm f}}.\label{eq:coordSKderiv}
\end{align}

For the field degrees of freedom we have to consider instead 
\begin{align}
    \nonumber\frac{\partial E_{\rm IBVP}}{\partial (\partial_0 \phi_+)} &=  \frac{\partial E_{\rm IBVP}}{\partial (\partial_0 \phi_1)} \frac{\partial_0 \phi_1}{\partial_0 \phi_+}+\frac{\partial E_{\rm IBVP}}{\partial (\partial_0 \phi_2)} \frac{\partial_0 \phi_2}{\partial_0 \phi_+}\\
    &=  \frac{\partial E_{\rm BVP}[X_1,\partial_a X_1,\phi_1,\partial_a\phi_1]}{\partial (\partial_0 \phi_1)} -\frac{\partial E_{\rm BVP}[X_2,\partial_a X_2,\phi_2,\partial_a\phi_2]}{\partial (\partial_0 \phi_2)}.
\end{align}
In turn from the field dependent terms in \cref{eq:novelactionE} one is led to expressions such as
\begin{align}
 \frac{\partial}{\partial (\partial_0 \phi_1)} \partial_a\phi_1\partial_b\phi_1 {\rm adj}[g]_{ab}= 2\, {\rm adj}[g]_{0b}\partial_b\phi_1
\end{align}
Since we have made sure by \eqref{eq:coordSKderiv} that the adjugate of the induced metric agrees betweeen the forward and backward branch at final $\tau$, we can thus write
\begin{align}
\frac{\partial E_{\rm BVP}}{\partial (\partial_0 \phi_1)} -\frac{\partial E_{\rm BVP}}{\partial (\partial_0 \phi_2)}= 2 {\rm adj}[g]_{0b}\Big( \partial_b\phi_1 - \partial_b\phi_2\Big).
\end{align}
The adjugate of the induced metric is in general a dense matrix. In order for the above expression to vanish, we find that it is necessary to identify the derivatives of the field with respect to all $\Sigma^a$ at final $\tau$
\begin{align}
\partial_a \phi_1|_{\tau=\tau^{\rm f}}= \partial_a \phi_2|_{\tau=\tau^{\rm f}}.\label{eq:fieldSKderiv}
\end{align}

The outcome of the above derivation entails that we need to fix all $(d+1)$ derivatives where the forward and backwards path meet. However in the continuum formalism, the fact, that we identify the values of our degrees of freedom at $\tau^{\rm f}$ already implies that all spatial derivatives too are the same\footnote{In the discrete setting, depending on the choice of stencil used for the derivative approximation, it may become necessary to introduce individual constraints for the spatial derivatives}. Then it is only the derivatives in $\tau$ that we need to enforce to be equal. 

The combination of \cref{eq:connval,eq:coordSKderiv,eq:fieldSKderiv} constitutes the connecting conditions summarized in the main text in \cref{eq:SKconnecting}.

\FloatBarrier

\begin{backmatter}

\section*{Competing interests}
  The authors declare that they have no competing interests.

\section*{Author's contributions}
    \begin{itemize}
         \item A.\ Rothkopf: construction of the novel action (contribution based on the world-line formalism and field theory), numerical implementation, writing, editing
         \item W.\ A.\ Horowitz: construction of the novel action (contribution based on string theory and field theory), editing
         \item J.\ Nordstr\"om: cross-validating of the novel approach (contribution based on applied analysis), editing
    \end{itemize}

%%%%%%%%%%%%%%%%%%%%%%%%%%%%%%%%%%%%%%%%%
%%                                     									       %%
%%  Bibliography part starts here								       %%
%%                                     									       %%
%%%%%%%%%%%%%%%%%%%%%%%%%%%%%%%%%%%%%%%%%

\bibliographystyle{stavanger-mathphys}

%%%%%%%%%%%%%%%%%%%%%%%%%%%%%%%%%%%%%%%%%
%%                                     									       %%
%%  Specify your BibTeX bibliography file here or manually insert references  %%
%%                                     									       %%
%%%%%%%%%%%%%%%%%%%%%%%%%%%%%%%%%%%%%%%%%

\bibliography{references}

% or include bibliography directly:
% \begin{thebibliography}
% \bibitem{b1}
% \end{thebibliography}

\end{backmatter}

%%%%%%%%%%%%%%%%%%%%%%%%%%%%%%%%%%%%%%%%%
%%                                     									       %%
%%  End of the document										       %%
%%                                     									       %%
%%%%%%%%%%%%%%%%%%%%%%%%%%%%%%%%%%%%%%%%%

\end{document}